\documentclass[11pt,leqno,english]{amsart}
\makeatletter
\def\paragraph{\@startsection{paragraph}{4}%
  \z@{.5\linespacing}{.5\linespacing}%
  {\bfseries}}
\def\@secnumfont{\ifnum\@toclevel>1\bfseries\else\mdseries\fi}
\makeatother

\usepackage[left=1.3in,top=1.3in,right=1.3in,nofoot]{geometry}
\usepackage{amsmath}
\usepackage{amsfonts}
\usepackage{amssymb}
\usepackage{amscd}
\usepackage[titletoc,title]{appendix}
\usepackage{hyperref}

\newcommand{\Q}{\mathbb{Q}}
\newcommand{\R}{\mathbb{R}}

\newcommand{\Z}{\mathbb{Z}}

\newcommand{\ca}{\mathcal{A}}
\newcommand{\Cb}{\mathcal{B}}

\newcommand{\Cd}{\mathcal{D}}

\newcommand{\Hh}{\mathcal{H}}
\newcommand{\Ii}{\mathcal{I}}

\newcommand{\cn}{\mathcal{N}}
\newcommand{\Oo}{\mathcal{O}}
\newcommand{\cp}{\mathcal{P}}

\newcommand{\cs}{\mathcal{S}}
\newcommand{\cx}{\mathcal{X}}
\newcommand{\cz}{\mathcal{Z}}

\newcommand{\superc}{supercuspidal}
\newcommand{\supers}{supersingular}
\newcommand{\repr}{representation}
\newcommand{\resp}{\mathrm{resp}}

\newcommand{\bb}{\mathbf{B}}

\newcommand{\bg}{\mathbf{G}}
\newcommand{\bgl}{\mathbf{\GL}}
\newcommand{\bm}{\mathbf{M}}
\newcommand{\bs}{\mathbf{S}}
\newcommand{\bx}{\mathbf{x}}
\newcommand{\bh}{\mathbf{H}}
\newcommand{\bn}{\mathbf{N}}
\newcommand{\bp}{\mathbf{P}}
\newcommand{\gs}{\mathbf{S}}
\newcommand{\bu}{\mathbf{U}}
\newcommand{\bz}{\mathbf{Z}}

\newcommand{\ba}{\backslash}
\newcommand{\rg}{\rightarrow}
\newcommand{\lgr}{\longrightarrow}

\DeclareMathOperator{\Supp}{{Supp}}
\DeclareMathOperator{\id}{{id}}
\DeclareMathOperator{\St}{{St}}
\DeclareMathOperator{\op}{{op}}
\DeclareMathOperator{\is}{{is}}
\DeclareMathOperator{\der}{{der}}
\DeclareMathOperator{\sep}{{sep}}
\DeclareMathOperator{\simplyc}{{sc}}
\DeclareMathOperator{\GL}{GL}
\DeclareMathOperator{\PGL}{PGL}
\DeclareMathOperator{\SL}{SL}
\DeclareMathOperator{\val}{{val}}
\DeclareMathOperator{\Hom}{{Hom}}
\DeclareMathOperator{\End}{{End}}
\DeclareMathOperator{\Ker}{{Ker}}
\DeclareMathOperator{\Ima}{{Im}}
\DeclareMathOperator{\Gal}{{Gal}}
\DeclareMathOperator{\Ind}{{Ind}}
\DeclareMathOperator{\ind}{{ind}}

\DeclareMathOperator{\ad}{{ad}}

\newcommand{\ikg}{\ind_K^GV}
\newcommand{\Cci}{C_c^\infty}
\newcommand{\Bs}{\bar{\sigma}}

\newtheorem{montheo}{Theorem}

\usepackage{xcolor}
\definecolor{teal}{rgb}{0.0, 0.5, 0.5}
\definecolor{forest}{rgb}{0.13, 0.55, 0.13}

\begin{document}

\title[A classification of irreducible mod $p$ representations of $p$-adic groups]
{A classification of irreducible admissible\\ mod $p$ representations of\\ $p$-adic reductive groups}

\author{N. Abe}
\address[N. Abe]{Creative Research Institution (CRIS), Hokkaido University, N21, W10, Kita-ku, Sapporo, Hokkaido 001-0021, Japan}
\thanks{The first-named author was supported by JSPS KAKENHI Grant Number 26707001.}
\email{abenori@math.sci.hokudai.ac.jp}

\author{G. Henniart} 
\address[G. Henniart]{Universit\'e de Paris-Sud, Laboratoire de Math\'ematiques d'Orsay, Orsay cedex F-91405 France;
CNRS, Orsay cedex F-91405 France}
\email{Guy.Henniart@math.u-psud.fr}

\author{F. Herzig} 
\address[F. Herzig]{Department of Mathematics, University of Toronto,
  40 St.\ George Street, Room 6290, Toronto, ON M5S 2E4, Canada}
\thanks{The third-named author was partially supported by a Sloan Fellowship and an NSERC grant.}
\email{herzig@math.toronto.edu}

\author{M.-F. Vign\'eras}
\address[M.-F. Vign\'eras]{Institut de Math\'ematiques de Jussieu, 175 rue du Chevaleret, Paris 75013 France}
\email{vigneras@math.jussieu.fr}

\date{}

\begin{abstract}
Let $F$ be a locally compact non-archimedean field, $p$ its residue characteristic, and $\textbf{G}$ a connected reductive group over $F$. Let $C$ be an algebraically closed field of characteristic $p$. We give a complete classification of irreducible admissible $C$-representations of $G=\bg(F)$, in terms of supercuspidal $C$-representations of the Levi subgroups of $G$, and parabolic induction. Thus we push to their natural conclusion the ideas of the third-named author, who treated the case $\bg=\GL_m$, as further expanded by the first-named author, who treated split groups $\bg$. As in the split case, we first get a classification in terms of supersingular representations of Levi subgroups, and as a consequence show that supersingularity is the same as supercuspidality.
\end{abstract}
\maketitle

\tableofcontents

\section{Introduction}\label{I}

\subsection{}\label{I.1} The study of congruences between classical modular forms has met considerable success in the past decades. When interpreted in the natural framework of automorphic forms and representations, such congruences naturally lead to representations over fields of positive characteristic, rather than complex representations. In our local setting, where the base field is a locally compact non-archimedean field $F$, this means studying representations of $G=\bg(F)$, where $\bg$ is a connected reductive group over $F$, on vector spaces over a field $C$ of positive characteristic $p$, which we assume algebraically closed. As $C$ is fixed throughout, we usually say representation instead of representation on a $C$-vector space or $C$-representation.

Our representations satisfy natural requirements: they are smooth, in that every vector has open stabilizer in $G$ (smoothness is always understood for   representations of $G$  or its subgroups), and most of the time they are \textbf{admissible}: a representation of $G$ on a $C$-vector space $W$ is admissible if it is smooth and for every open subgroup $J$ in $G$, the space $W^J$ of vectors fixed under $J$ has finite dimension. The overall goal is to understand irreducible admissible representations of~$G$.

Here we consider only the case where the residue characteristic of $F$ is $p$.

\subsection{}\label{I.2} In this paper we classify irreducible admissible representations of $G$ in terms of parabolic induction and supercuspidal representations of Levi subgroups of~$G$.
Such a classification was obtained for $\bg=\GL_2$ in the pioneering work of L.~Barthel and R. Livn\'e \cite{BL1,BL2} -- see also some recent work \cite{Abd,Che,Ko,KX,Ly2} on situations where, mostly, $\bg$ has relative semisimple rank~$1$.  

New ideas towards the general case were set forth by the third-named author \cite{He1,He2}, who gave the classification for $\bg=\GL_n$ over a $p$-adic field $F$; his ideas were further expanded by the first-named author \cite{Abe} to treat the case of a split group $\bg$, still over a $p$-adic field $F$. T.\ Ly extended the arguments of \cite{He1,He2} to treat 
$\bg=\GL_{3/D}$ where $D$ is a division algebra over $F$, allowing $F$ to have characteristic $p$.
Here, using the first steps accomplished in \cite{HV1,HV2} (see \ref{I.5}, \ref{I.6}), we treat  general $\bg$ and $F$.

\subsection{}\label{I.3} To express our classification, we recall parabolic induction. If $P$ is a parabolic subgroup of $G$ and $\tau$ a  representation of $P$ on a $C$-vector space $W$, we write $\Ind_P^G \tau$ for the natural representation of $G$, by right translation, on the space $\Ind_P^GW$ of smooth functions $f:G\rg W$ such that $f(pg)=\tau(p)f(g)$ for $p$ in $P$, $g$ in $G$. The functor  $\Ind_P^G$ is exact. In fact we use $\Ind_P^G \tau$ only when $\tau$ comes via inflation from a  representation $\sigma$ of the Levi quotient of $P$, and we write $\Ind_P^G\sigma$ instead of $\Ind_P^G\tau$.
A representation of $G$ is said to be \textbf{supercuspidal} if it is  irreducible, admissible, and does not appear as a  subquotient of a parabolically induced representation $\Ind_P^G\sigma$, where $P$ is a  proper parabolic subgroup of $G$ and $\sigma$ an irreducible admissible representation of the Levi quotient of~$P$.
First we construct irreducible admissible representations of $G$. The construction uses the ``generalized Steinberg'' representations investigated by E.~Gro\ss e-Kl\"onne \cite{GK} and the third-named author \cite{He2} when $\bg$ is split, and by T.~Ly \cite{Ly1} in general: for any pair of parabolic subgroups $Q\subset P$ in $G$, $\St_Q^P$ is the natural representation of $P$ in the quotient of $\Ind_Q^P 1$ by the sum of the subspaces $\Ind_{Q'}^P1$, for parabolic subgroups $Q'$ with $Q\subsetneq Q'\subset P$; the representation $\St_Q^P$ factors through the unipotent radical $U_P$ of $P$ and gives the representation $\St_{Q/U_P}^{P/U_P}$ of its reductive quotient, so $\St_Q^P$ is irreducible and admissible \cite{GK,Ly1}.

Start with a parabolic subgroup $P$ of $G$, with Levi quotient  $M$, and a representation $\sigma$ of $M$. Then there is a largest parabolic subgroup $P(\sigma)$ of $G$, containing $P$, such that $\sigma$ inflated to $P$ extends to $P(\sigma)$ (see \ref{II.7}). That extension is unique, we write it $^e\sigma$; it is  trivial on the unipotent radical of $P(\sigma)$. It is irreducible and admissible if $\sigma$ is.  We consider triples $(P,\sigma,Q)$: a triple consists of a parabolic subgroup $P$ of $G$, a  representation $\sigma$ of the Levi quotient $M$ of $P$, and a parabolic subgroup $Q$ of $G$ with $P\subset Q \subset P(\sigma)$; we say that the triple is \textbf{supercuspidal} if $\sigma$ is a supercuspidal representation of $M$. To a triple $(P,\sigma,Q)$ we associate the representation $I(P,\sigma,Q)=\Ind_{P(\sigma)}^G(^e\sigma\otimes \St_Q^{P(\sigma)})$.

\begin{montheo}
For a supercuspidal triple $(P,\sigma,Q)$, $I(P,\sigma,Q)$ is irreducible and admissible.
\end{montheo}

\begin{montheo}
Let $(P,\sigma,Q)$ and $(P',\sigma',Q')$ be supercuspidal triples.
Then $I(P,\sigma,Q)$ and $I(P',\sigma',Q')$ are isomorphic if and only if there is an element $g$ of $G$ such that $P'=gPg^{-1}$, $Q'=g\,Qg^{-1}$ and $\sigma'$ is equivalent to $p'\mapsto\sigma(g^{-1}p'g)$.
\end{montheo}

\begin{montheo}
Any irreducible  admissible representation of $G$ is isomorphic to $I(P,\sigma,Q)$ for some supercuspidal triple $(P,\sigma,Q)$.
\end{montheo}

Hopefully the classification expressed by these theorems will be useful in extending the mod $p$ local Langlands correspondence beyond $\GL_2(\Q_p)$.

\subsection{}\label{I.4} Using the classification results above, it is possible to describe the irreducible components of $\Ind_P^G\sigma$ where $P$ is a parabolic subgroup of $G$ and $\sigma$ an irreducible admissible representation of the Levi quotient $M$ of $P$; in particular we show that $\Ind_P^G\sigma$ has finite length and that all its irreducible subquotients are admissible and occur with multiplicity one. 
Also we have a notion of \textbf{supercuspidal support}: if $(P,\sigma,Q)$ is a supercuspidal triple, then $\pi=I(P,\sigma,Q)$ occurs as a subquotient of $\Ind_P^G \sigma$ and if $\pi$ occurs as a subquotient of $\Ind_{P'}^G\sigma'$ for a supercuspidal representation $\sigma'$ of (the Levi quotient of) a parabolic subgroup $P'$ of $G$ then $(P',\sigma')$ is conjugate to $(P,\sigma)$ in $G$ as in Theorem 2.  All that is proved in \ref{VI.3}. It is the conjugacy class of $(P,\sigma)$ that we call the supercuspidal support of~$\pi$.  

\vskip2mm
\noindent\textbf{Remark} Even in the case of $\GL_n(F)$ (for which we refer to the introduction of \cite{He2}), the classification and its consequences are rather simpler than for complex representations:  intertwining operators do not exist in our context; this ``explains'' the multiplicity one result above, which does not hold for complex representations \cite{Ze}. By contrast, supercuspidal  mod $p$ representations remain a complete mystery, apart from the case of $\GL_2(\Q_p)$ \cite{Br} and groups closely related to it \cite{Abd,Che, Ko, KX}.

The existence of a supercuspidal support for complex irreducible representations is a classical
result \cite[2.9 Theorem]{BZ1}; for mod $\ell$ representations with $\ell \ne p$ it is unknown (even for finite reductive
groups of characteristic $p$ outside the case of general linear groups), except for inner forms of
$\GL_n(F)$ where, as above, it is not proved directly but is a consequence of the classification of
irreducible representations  \cite[Th\'eor\`eme A]{MS}.

\subsection{}\label{I.5} As in \cite{He2,Abe} our classification is not established directly using supercuspidality. Rather we get a classification in terms of  supersingular representations of Levi subgroups of $G$ -- the term was first used by Barthel and Livn\'e for $G = \GL_2(F)$ -- and deduce Theorems 1 to 3 from it.  To define supersingularity, we need to make some  choices, and a priori the notion depends on these choices.

So we fix a maximal $F$-split torus $\gs$ in $\bg$ and a special point $\mathbf{x}_0$ in the apartment corresponding to $\gs$ in the semisimple Bruhat-Tits building of $G$; we let $K$ be the special parahoric subgroup of $G$ corresponding to $\mathbf{x}_0$. We also fix a minimal parabolic subgroup $B$ of $G$ with Levi subgroup $Z$, the $F$-points of the centralizer   of $\gs$, and we write $U$ for the unipotent radical of $B$.

Let $V$ be an irreducible representation of $K$ -- it has finite dimension. If $(\pi,W)$ is an admissible representation of $G$, then $\Hom_K(V,W)$ is a finite-dimensional $C$-vector space; by Frobenius reciprocity $\Hom_K(V,W)$ is identified with $\Hom_G(\ind_K^G V,W)$, where $\ind_K^G$ denotes compact induction, so that $\Hom_K(V,W)$ is a right-module over the intertwining algebra $\Hh_G(V)=\End_G(\ind_K^G V)$ of $V$ in $G$.
If $\Hom_K(V,W)$ is not zero we say that $V$ is a \textbf{weight} of $\pi$; in that case the centre\footnote{Note that $\Hh_G(V)$ is commutative in many cases, for example when $G$ is split, but not in general \cite{HV1}.} $\cz_G(V)$ of $\Hh_G(V)$ has eigenvectors in $\Hom_K(V,W)$, and we focus on the corresponding characters of $\cz_G(V)$, which we call the (Hecke) \textbf{eigenvalues} of $\cz_G(V)$ in~$\pi$.

For any parabolic subgroup $P$ of $G$ containing $B$, with Levi component $M$ containing $Z$ and unipotent radical $N$, the space of coinvariants $V_{N\cap K}$ of $N\cap K$ in $V$ provides an  irreducible representation of $M\cap K$ and by \cite{He1,He2, HV2} there is a natural injective algebra homomorphism 
$$ \cs _M^G: \Hh_G(V) \rg \Hh_M(V_{N\cap K}) 
$$   with explicit image  (see \ref{III.3}). It induces a homomorphism between centres $\cz_G(V)\rg\cz_M(V_{N\cap K})$. Both homomorphisms are localizations at a central element.
A character $\chi:\cz_G(V)\to C$ is said to be \textbf{\supers} if, in the above situation, it can be extended to a character of $\cz_M(V_{N\cap K})$ only when $P=G$ (see Chapter \ref{III}, part A) for details). A \textbf{supersingular} representation of $G$ is an irreducible admissible representation $(\pi,W)$ such that for all weights $V$ of $\pi$, all eigenvalues of $\cz_G(V)$ in $\pi$ are supersingular\footnote{That is consistent with the definition in \cite{He2,Abe}; but the reader should be aware that the definition in \cite{HV2} is slightly different, maybe not equivalent.}.

A triple $(P,\sigma,Q)$  as in \ref{I.3} is a \textbf{$B$-triple} if $P$  contains $B$; it is said to be \textbf{supersingular} if  
it is a $B$-triple  and $\sigma$ is a supersingular representation of the Levi quotient of~$P$.

Theorems 1 to 3 are consequences of the following results.

\begin{montheo}
For each supersingular triple $(P,\sigma,Q)$, the representation $I(P,\sigma,Q)$ is irreducible and admissible. If $\pi$ is an irreducible admissible representation of $G$, there is a supersingular triple $(P,\sigma,Q)$ such that $\pi$ is isomorphic to $I(P,\sigma,Q)$; moreover $P$ and $Q$ are unique and $\sigma$ is unique up to isomorphism.
\end{montheo}

\begin{montheo}
Let $\pi$ be an irreducible admissible \repr\ of $G$. Then $\pi$ is \superc\  if and only if it is \supers.
\end{montheo}

(For $\bg=\GL_2$ this was discovered by Barthel and Livn\'e.)

\vskip2mm
Note that Theorem 5 implies, in particular, that the notion of supersingularity does not depend on the choices of $\gs$, $K$, $B$ necessary for the definition -- beware that in general two choices of $K$ will not even be conjugate under the adjoint group of $G$.

\vskip2mm

\noindent\textbf{Remarks} 1) We also show that, if $\pi$ is an irreducible admissible \repr\ of $G$, and for some weight $V$ of $\pi$ there is an eigenvalue of $\cz_G(V)$ in $\pi$ which is \supers, then $\pi$ is \supers/\superc.

2) Let $(P,\sigma,Q)$ be a supersingular (or supercuspidal) $B$-triple. Then $I(P,\sigma,Q)$ is finite dimensional if and only if $P=B$ and $Q=G$.

\subsection{}\label{I.6} As in \cite{He2} and \cite{Abe}, a lot of our arguments bear on the relation between parabolic induction  $\Ind_P^G$ in $G$ and compact induction  $\ind_K^G$ from $K$ to $G$.

Let $V$ be an irreducible  \repr\ of $K$, and let $P$ be a parabolic subgroup of $G$ containing $B$, with Levi component $M$ containing $Z$, and unipotent radical $N$. In \cite{HV2}, inspired by \cite{He1}, \cite{He2}, a canonical intertwiner
$$
\Ii: \ind_K^G V \lgr \Ind_P^G(\ind_{M\cap K}^M V_{N\cap K})
$$
was investigated. In fact the morphism $ \cs_M^G$ of \ref{I.5} is such that for $f$ in $\ind_K^GV$ and $\Phi$ in $\Hh_G(V)$ we have
$$
\Ii(\Phi(f)) = \cs_M^G(\Phi)(\Ii(f)),
$$
where the action of $\cs_M^G(\Phi)$ on $\Ii(f)$ is via its natural action on $\ind_{M\cap K}^MV_{N\cap K}$.
Under a suitable regularity condition of $V$ with respect to $P$ \cite{HV2}, cf.\ \ref{III.14} Theorem, $\Ii$ induces an isomorphism 
$$\chi\otimes \ind_K^G V \overset{\sim}{\underset{}{\lgr}}\Ind_P^G(\chi \otimes \ind_{M\cap K}^M V_{N\cap K})$$ for any character $\chi$ of $\cz_G(V)$ which extends to $\cz_M(V_{N\cap K})$: such an extension is unique, we still denote it by $\chi$; the first tensor product is over $\cz_G(V)$, the second one over $\cz_M(V_{N\cap K})$. Here we obtain a generalization of that result, which we now proceed to explain.

We consider an irreducible  \repr\ $V$ of $K$, and a character $\chi:\cz_G(V)\to C$. There is a smallest parabolic subgroup $P$ containing $B$ -- we write $P=MN$ as above -- such that $\chi$ extends to a character, still written $\chi$, of $\cz_M(V_{N\cap K})$; there is a  natural parabolic subgroup $P_e$, containing $P$, such that the  \repr\ $\chi\otimes (\ind_{M\cap K}^M V_{N\cap K})$ of $M$, inflated to $P$, extends to a \repr\ of $P_e$ --  write $^e(\chi \otimes \ind_{M\cap K}^M V_{N\cap K})$ for that extension. Using similar notation as in \ref{I.3}, we write  $I_e(P,\chi\otimes \ind_{M\cap K}^M V_{N\cap K},Q)$ for $\Ind_{P_e}^G(^e(\chi\otimes \ind_{M\cap K}^M V_{N\cap K}) \otimes \St_Q^{P_e})$ when $Q$ is a parabolic subgroup between $P$ and $P_e$.

\begin{montheo}[Filtration Theorem]
With the previous notation, $\tau=\chi \otimes \ind_K^G V$ has a natural filtration by sub\repr s $\tau_Q$, where $Q$ runs through parabolic subgroups of $G$ with $P\subset Q\subset P_e$ and $\tau_{Q'}\subset\tau_{Q}$ if $Q'\subset Q$.
 The quotient $\tau_Q/\sum\limits_{Q'\subsetneq Q}\tau_{Q'}$ is isomorphic to $I_e(P,\chi\otimes \ind_{M\cap K}^M V_{N\cap K},Q)$.
\end{montheo}

This last theorem should be compared to the following (the proof, in Chapter~\ref{V}, explains that comparison). Let $\pi = \Ind_P^G(\chi\otimes\ind_{M\cap K}^M V_{N\cap K})$. It also has a natural filtration by sub\repr s $\pi_Q$ for $Q$ as above, but this time $\pi_{Q'} \subset \pi_{Q}$ if $Q'\supset Q$, and the quotient $\pi_{Q}/\sum\limits_{Q'\supsetneq Q}\pi_{Q'}$ is isomorphic to $I_e(P_e,\chi \otimes\ind_{M\cap K}^M V_{N\cap K},Q)$.
In particular the filtrations on $\tau$ and $\pi$ give rise to the same  subquotients, but in reserve order, so to say. (We note that the representation $\pi_Q$ above corresponds to the representation $I_Q$ in Chapter~\ref{V}.) 

A striking example is when $V$ is trivial character of $K$ and $\chi$ is the ``trivial'' character of $\cz_G(V)=\Hh_G(V)$: in that case $P=B=ZU$,   $P_e=G$, and  $\chi\otimes \ind_{Z\cap K}^Z V_{U\cap K}$ is the trivial character of $Z$. In $\pi =\Ind_B^G 1$, the trivial character of $G$ occurs as a sub\repr\ and the Steinberg \repr\ $\St_B^G$ as a quotient, whereas the reverse is true in $\chi\otimes \ind_K^G1$.

Theorem 6 is new even for $\GL_n$ ($n > 2$). A weaker version of this theorem is proved in
\cite[Proposition 4.7]{Abe} when $\bg$ is split with simply connected derived subgroup and $P = B$
(and in \cite{BL2} in the further special case when $\bg = \GL_2$). On the way, following the ideas
of \cite{Abe}, we prove the freeness of $R_M \otimes_{\cz_G(V)} \ind_K^G V$ as $R_M$-module, where
$R_M$ denotes the ``supersingular quotient'' of $\cz_M(V_{N \cap K})$. This may
be of independent interest. Again this result was established for  $\bg = \GL_2$ in \cite{BL1}, but see also
the recent paper \cite{GK2}.
\subsection{}\label{I.7} To prove Theorem~4 we follow the same strategy as in \cite{He2,Abe} (see the introduction of \cite{He2} for an outline). If
$(P,\sigma,Q)$ is a supersingular triple, we need to prove that $\pi=I(P,\sigma,Q)$ is irreducible;
that is done by showing that for any weight $V$ of $\pi$ and any eigenvector $\varphi$ for
$\cz_G(V)$ in $\Hom_K(V,\pi)$ with corresponding eigenvalue $\chi$, $\pi$ is generated as a \repr\
of $G$ by the image of $\varphi$. When $V$ is suitably regular, that is seen as a consequence of the
isomorphism $\chi\otimes \ind_K^G V \simeq \Ind_P^G(\chi\otimes \ind_{M\cap K}^M V_{N\cap K})$
recalled in \ref{I.6} above (see \ref{III.14}).  We reduce to that suitably
regular case by using a \emph{change of weight theorem}, which gives explicit sufficient conditions
on $V$, $V'$, and $\chi$ for having an isomorphism $\chi\otimes \ind_K^G V \simeq \chi\otimes
\ind_K^G V'$.  (Here, $V'$ is an irreducible representation of $K$ that is ``slightly less regular''
than $V$ and such that $(V')_{U\cap K} \simeq V_{U\cap K}$.)  We refer the reader to Sections
\ref{IV.2}, \ref{III.23} for the precise statement and its use in the proof of Theorem~4.

The main novelty in our methods is our proof of the change of weight
theorem. It is also the hardest and most subtle part of our arguments.
Previously, for split groups, a version of this theorem was established in \cite[\S 6]{He2}
and \cite[\S 4]{Abe} by computing the composition of two intertwining operators and applying the
Lusztig--Kato formula. We do not know if this approach can be generalized.  Our new proof does not
involve Kazhdan--Lusztig polynomials, but rather proceeds by embedding $\ind_K^G V$, $\ind_K^G V'$
into the parabolically induced representation $\cx = \Ind_B^G(\ind_{Z \cap K}^Z \psi_V)$ using the
intertwiner $\Ii$ of \ref{I.6}, where $\psi_V : Z\cap K \to C^\times$ describes the action of $Z\cap
K$ on $V_{U\cap K} \simeq (V')_{U\cap K}$. The representation $\ind_K^G V$
contains a canonical (up to scalar) fixed vector under a pro-$p$ Iwahori subgroup $I \subset K$
which generates $\ind_K^G V$ as a representation of $G$, and similarly for $\ind_K^G V'$. 
Our proof then studies the action of the pro-$p$-Iwahori
Hecke algebra $\End_G(\ind_I^G 1)$ on $\cx^I$ to relate the two compact inductions inside $\cx$.  We
crucially rely on the description of the pro-$p$-Iwahori Hecke algebra recently given for general $G$ by
the fourth-named author in \cite{Vig4}, in particular the Bernstein relations in this algebra.

We arrive at a \emph{dichotomy} in \ref{IV.1} Theorem and \ref{IV.2} Corollary, namely
our change of weight results depend on whether or not $\psi_V$ is trivial on a certain subgroup of
$Z \cap K$. When $G$ is split, the triviality is always guaranteed, but that is not
always so for inner forms of $\GL_n$ \cite[Lemme 3.10.1]{Ly3} and even
for unramified unitary groups in 3 variables. This dichotomy may explain why
we did not find an easy generalization of the previous proofs for split $G$.

\subsection{}\label{I.8}  Let $\pi$ be an irreducible admissible representation of $G$, $P=MN$ a parabolic subgroup of $G$,  and $\tau$ an irreducible admissible representation of $M$ inflated to $P$.  In a sequel to this article we will apply our classification to tackle natural questions as   the computation of the $N$-coinvariants or  the $P$-ordinary part of $\pi$,  
the description of the lattice of subrepresentations of $\Ind_P^G \tau$,   the generic    irreducibility of the representations $\Ind_P^G \tau \chi$  where $\chi$ runs over the unramified characters of $M$ (this question was raised by J.-F. Dat).

\subsection{}\label{I.9} We end this introduction with some comments on the organization of the paper. In Chapter \ref{II} we fix notation and we examine  when a \repr\ of a parabolic subgroup of $G$, trivial on its unipotent radical, can be extended to a larger parabolic subgroup. For a triple $(P,\sigma,Q)$ as in \ref{I.3}, we construct $I(P,\sigma,Q)$ and show that  it is admissible if $\sigma$ is. In Chapter \ref{III} we give most of the proof of Theorem~4. The irreducibility proof was outlined in \ref{I.7}. The proof that $\pi = I(P,\sigma,Q)$ determines $P$, $Q$, and $\sigma$ up to isomorphism comes from examining the possible weights and Hecke eigenvalues for $\pi$ (\ref{III.24}). Finally, to prove that every irreducible admissible \repr\ $\pi$ of $G$ has the form $I(P,\sigma,Q)$ we use the filtration theorem (Theorem~6). The proof of the change of weight theorem 
is given in Chapter \ref{IV}; this is the technical heart of our paper. In Chapter \ref{V} we deduce the 
filtration theorem from the change of weight theorem. We trust that the reader will see easily that there is no loop in our arguments. Finally, Chapter \ref{VI} gives the proof of Theorems 1, 2, 3, 5 and other consequences of the classification, already stated in \ref{I.4}. That section can essentially be read independently, taking Theorem~4 for granted.

\medskip \noindent \textbf{Acknowledgments} 
We thank the following institutions, where part of our work was carried out: IHES,  IMJ Paris 7, IMS Singapore, MSRI, Paris 11.
We thank the referee for helpful comments.

\section{Extension to a larger parabolic subgroup}\label{II}

\subsection{}\label{II.1} Let us first fix notation, valid throughout the paper. As stated in the introduction, our base field $F$ is locally compact and non-archimedean, of residue characteristic $p$; its ring of integers is $\Oo$, its residue field $k$, and $q$ is the cardinality of $k$; we write $|\ |$ for the normalized absolute value of $F$.

A linear algebraic group over $F$ will be written with a boldface letter like $\bh$, and its group of $F$-points will be denoted by the corresponding ordinary letter $H=\bh(F)$\footnote{We shall use a similar convention for groups over $k$.}.

We fix our connected reductive $F$-group $\bg$\footnote{$\bg$ is fixed, but otherwise arbitrary, so the results we establish for $\bg$ can be applied to other reductive groups over~$F$.}, and a maximal $F$-split torus $\gs$ in $\bg$; we write $\bz$ for the centralizer of $\gs$ in $\bg$, ${\boldsymbol{\cn}}$ for its  normalizer, and $W_0=W(\bg,\gs)$ for the Weyl group ${\boldsymbol{\cn}}/\bz$; we recall that $W_0=\cn/Z$ \cite[21.2 Theorem]{Bo}. We also fix a minimal $F$-parabolic subgroup $\bb$ of $\bg$ with Levi subgroup $\bz$, and write $\bu$ for its unipotent radical. As is customary, we say that $P$ is a parabolic subgroup of $G$ to mean that $P=\bp(F)$, where $\bp$ is an $F$-parabolic subgroup of $\bg$. If $P$ contains $B$, we usually write $P=MN$ to mean that $M$ is the Levi component of $P$ containing $Z$, and $N$ the unipotent radical of $P$; we then write $P_{\op}=M N_{\op}$ for the parabolic subgroup opposite to $P$ with  respect to~$M$; in particular $B_{\op}=ZU_{\op}$.

We let $\Phi$ be the set of roots of $\gs$ in $\bg$, so $\Phi$ is a subset of the group $X^*(\gs)$ of characters of $\gs$; we let $\Phi^+$ be the subset of roots of $\gs$ in $\bu$, called positive roots, and $\Delta$ for the set of simple roots of $\gs$ in $\bu$. If $X_*(\gs)$ is the group of cocharacters of $\gs$ we write $\langle\; , \, \rangle$ for the natural pairing $X^*(\gs)\times X_*(\gs) \rg \Z$; for $\alpha$ in $\Phi$, the corresponding coroot \cite[expos\'e XXVI, \S7]{SGA3} is written $\alpha^\vee$ and for $I\subset \Phi$ we put $I^\vee =\{\alpha^\vee\mid \alpha\in I\}$. We choose a positive definite symmetric bilinear form on $X^*(\gs) \otimes_\Z \R$, invariant under $W_0$, which induces a notion of orthogonality between roots; for roots $\alpha$, $\beta$ we have $\alpha\perp \beta$ if and only if $\langle\alpha,\beta^\vee\rangle=0$.

For $\alpha$ in $\Phi$ we write $U_\alpha=\bu_\alpha(F)$, for the corresponding root subgroup ($\bu_\alpha$ is written $\bu_{(\alpha)}$ in \cite[\S21]{Bo}), and $s_\alpha\in W_0$ for the corresponding reflection. For $I\subset \Delta$ we let $W_I$ be the subgroup generated by $\{s_\alpha\mid \alpha\in I\}$, $\cn_I$ for the inverse image of $W_I$ in $\cn$, $P_I$ for the parabolic subgroup $U \cn_IU$ (it contains $B$), $P_I=M_IN_I$ for its Levi decomposition, $M_I$ containing $Z$; if $I$ is a singleton $\{\alpha\}$ we rather write $P_\alpha=M_\alpha N_\alpha$. We set $\Delta_P=I$ if $P=P_I$.  We note that for $I$, $J\subset\Delta$, $P_{I\cap J}=P_I\cap P_J$, $M_{I\cap J}=M_I \cap M_J$.

\subsection{}\label{II.2} As announced in the introduction, we tackle here a preliminary question: if $P$ is a parabolic subgroup of $G$ and $\sigma$ a representation of $P$ trivial on its unipotent radical $N$, when can $\sigma$ be extended to a larger parabolic subgroup $Q$ of $G$? Dividing by the unipotent radical of $Q$, which is contained in $N$, we loose no generality in assuming that $Q=G$. If $\sigma$ extends to $G$, then any extension has to be trivial on the normal subgroup $\langle\, {}^G N\rangle$ of $G$ generated by $N$, so that $\sigma$ has to be trivial on $P\cap\langle \, {}^G N\rangle$. So we need to understand what $\langle\, {}^G N\rangle$ is. That question, which involves no representation theory, will be dealt with presently.

 \subsection{}\label{II.3} Of particular importance in our setting will be the subgroup $G'$ of $G$ generated by $U$ and $U_{\op}$. Beware that the notation, which will be applied to other reductive groups (like the Levi subgroups of $\bg$), is unusual, and that $G'$ is not generally the group of points over $F$ of a reductive subgroup of $\bg$: this occurs already for $\bg=\PGL_2$. Since $G$ is generated by $U$, $U_{\op}$ and $Z$, see e.g.\ \cite[Proposition 6.25]{BoT}, $G'$ is normal in $G$ so is also the subgroup of $G$ generated by the unipotent radicals of the parabolic subgroups of $G$, and we have $G=ZG'$. Sometimes we have $G'=G$, though.
 
 \vskip2mm
\noindent\textbf{Proposition} \textit{Assume that $\bg$ is semisimple, simply connected, almost $F$-simple and isotropic. Then $G'=G$, and $G$ has no non-central proper normal subgroup. Moreover, $Z$ is generated by the $Z\cap M_\alpha'$, $\alpha$ running through} $\Delta$.
\vskip2mm
 
 \noindent \textbf{Proof} The first assertion is due to Platonov \cite[Theorem 7.6]{PlR} and the second one then follows from work of Tits \cite[Theorem 7.1]{PlR}. The final assertion is due to  Prasad and Raghunathan \cite{PrR} -- actually their result is valid over any field. $\square$
 
 \vskip2mm
 
 \noindent\textbf{Remark} Let $\bg$ be as in the proposition, let $\alpha\in \Delta$ and $\bg_\alpha$ the subgroup of $\bg$ generated by $\bu_\alpha$ and $\bu _{-\alpha}$; since $\bg_\alpha$ satisfies the hypotheses of the proposition, we have $M_\alpha'=G_\alpha'=G_\alpha$.
 
 \subsection{}\label{II.4} In the following sections (\ref{II.5}--\ref{II.8}) our strategy is to reduce statements for $G$ to a much        
simpler group $G^{\is}$ via a homomorphism $G^{\is} \to G$ whose image is $G'$. The group $G^{\is}$ has the property that it is a product of groups of the      
form considered in \ref{II.3} Proposition, and the homomorphism $G^{\is} \to G$ restricts to an isomorphism on unipotent radicals of parabolic          
subgroups.

Let  $\bg^{\simplyc}$ be the simply connected covering  of the derived group $\bg^{\der}$ of $\bg$. 
Recall that $\bg^{\simplyc}$ is the direct product of its almost $F$-simple components. We let $\Cb$ be an indexing set for the \textbf{isotropic} almost $F$-simple components of $\bg^{\simplyc}$
 and for $b\in \Cb$ we write $\tilde{\bg}_b$ for the corresponding component. We put $\bg^{\is} = \prod\limits_{b\in \Cb} \tilde{\bg}_b$, and denote by $\iota$ the natural homomorphism ${\bg}^{\is}\rg \bg$, factoring through $\bg^{\is}\rg \bg^{\simplyc} \rg \bg^{\der} \rg \bg$.
 
 We need to understand the relation between parabolic subgroups of $\bg$ and parabolic subgroups of $\bg^{\is}$. The following comes from \cite[\S21, \S22]{Bo}, going through the factorization of $\iota$.
 
 The connected component of $\iota^{-1}(\gs)$ is a maximal $F$-split torus $\tilde{\gs}$ of $\bg^{\is}$, and $\gs$ is the product of $\iota(\tilde{\gs})$ and the maximal $F$-split torus in the centre of $\bg$. The centralizer of  $\tilde{\gs}$ in $\bg^{\is}$ is $\tilde{\bz}= \iota^{-1}(\bz)$, its normalizer $\tilde{\boldsymbol{\cn}}=\iota^{-1}(\boldsymbol{\cn})$, and $\iota$ induces an isomorphism $W(\bg^{\is},\tilde{\gs}) = W(\bg,\gs)$ (see in particular \cite[22.6 Theorem]{Bo}); in particular $W_0$ has representatives in $\iota(G^{\is})$. As $\bg^{\is}$ is a direct product $\prod \tilde{\bg}_b$ (over $b\in \Cb$) we have corresponding natural decompositions $\tilde{\gs}=\prod \tilde{\gs}_b$, $\tilde{\bz}=\prod \tilde{\bz}_b$, $\tilde{\boldsymbol{\cn}}=\prod \tilde{\boldsymbol{\cn}}_b$ and $W(\bg^{\is},\gs)= \prod W(\tilde{\bg}_b,\tilde{\gs}_b)$. Note that $\iota(\tilde{\bg}_b)$ is normal in $\bg$ for each $b\in \Cb$.
 
 The map $\bp \mapsto\tilde{\bp}= \iota^{-1}(\bp)$ is a bijection between $F$-parabolic subgroups of $\bg$ and $F$-parabolic subgroups of $\bg^{\is}$, and $\iota$ induces an isomorphism (cf.\ \cite[22.6 Theorem]{Bo}) of the unipotent radical $\tilde{\bn}$ of $\tilde{\bp}$ onto the unipotent radical $\bn$ of $\bp$. Also, $\tilde{\bm}= \iota^{-1}(\bm)$ is the Levi component of $\tilde{\bp}$ containing $\tilde{\bz}$. In particular $\tilde{\bb}=\iota^{-1}(\bb)$ is a minimal $F$-parabolic subgroup of $\bg^{\is}$; it is the direct product of minimal parabolic subgroups $\tilde{\bb}_b$ of $\tilde{\bg}_b$, and its unipotent radical $\tilde{\bu}$ is the direct product of the $\tilde{\bu}_b$, with $\tilde{\bu}_b$ the unipotent radical of $\tilde{\bb}_b$. Via $\iota$ we get an identification\footnote{More precisely the natural map $\tilde{\gs}\rg \gs$ induces a group homomorphism $X^*(\gs) \rg X^*(\tilde{\gs})$ through which the roots of $\gs$ in $\bu$ are identified with the roots of $\tilde{\gs}$ in $\tilde{\bu}$. By \cite[Exp.\ XXVI, 7.4]{SGA3} if $\alpha$ is a root of $\gs$ in $\bu$ and $\tilde{\alpha}$ the corresponding root of $\tilde{\gs}$ in $\tilde{\bu}$, then $\tilde{\alpha}^\vee$ goes to $\alpha^\vee$ via the transposed morphism $X_*(\tilde{\gs})\rg X_*(\gs)$. In the sequel we make no distinction between $\alpha$ and $\tilde{\alpha}$, $\alpha^\vee$ and $\tilde{\alpha}^\vee$.} of the roots of $\gs$ in $\bu$ with the roots of $\tilde{\gs}$ in $\tilde{\bu}$, so that $\Delta$, in particular, also appears as the set of simple roots of $\tilde{\gs}$ in $\tilde{\bu}$; as such $\Delta$ is a disjoint union of the sets $\Delta_b$, $b\in\Cb$, where $\Delta_b$ is the set of roots of $\tilde{\gs}$ (or $\tilde{\gs}_b$) in $\tilde{\bu}_b$; that partition of $\Delta$ is the finest partition into mutually orthogonal subsets. Those subsets are the connected components of the Dynkin diagram of $\bg$ (with set of vertices $\Delta$) so we can view $\Cb$ as the set of such components.

\vskip2mm

\noindent\textbf{Proposition} 
$G'=\iota(G^{\is}).
$
\vskip2mm

\noindent\textbf{Proof} By \ref{II.3} Proposition we have $\tilde{G}_b'=\tilde{G}_b$ for each $b\in \Cb$ so $(G^{\is})'=G^{\is}$; since $\iota$ induces an isomorphism of $\tilde{U}$ onto $U$ and $\tilde{U}_{\op}$ onto $U_{\op}$, we get $G'=\iota(G^{\is})$. $\square$

\vskip2mm

Note that the proposition implies that $Z \cap G'=\iota(\tilde Z)$.

\subsection{}\label{II.5} \noindent\textbf{Notation} For  $I\subset\Delta$, set $\Cb(I)=\{b\in \Cb \mid I\cap \Delta_b\not= \Delta_b\}$.

\vskip2mm

\noindent\textbf{Proposition} 
\textit{Let $I\subset \Delta$. Then the normal subgroup $\langle\, ^G N_I\rangle$ of $G$ generated by $N_I$ is} $\iota(\prod\limits_{b\in \Cb(I)}\tilde{G}_b)$.
\vskip2mm

\noindent\textbf{Proof} We have $\tilde{N}_I=\prod\limits_{b\in \Cb}( \tilde{N}_I \cap \tilde{G}_b)$ and $\tilde{N}_I\cap \tilde{G}_b$ is the unipotent subgroup of $\tilde{G}_b$ corresponding to $I\cap \Delta_b\subset \Delta_b$. For $b\in \Cb -\Cb(I)$, $\tilde{N}_{I}\cap \tilde{G}_b$ is trivial; for $b$ in $\Cb(I)$, $\tilde{N}_{I}\cap \tilde{G}_b$ is non-trivial, and provides a non-central subgroup of $\tilde{G}_b$ so by \ref{II.3} Proposition the normal subgroup of $G^{\is}$ generated by $\tilde{N}_{I}\cap \tilde{G}_b$ is $\tilde{G}_b$; the proposition follows. $\square$

\vskip2mm

\noindent \textbf{Corollary}
\begin{enumerate}
\item $\displaystyle P_I \cap \langle\, ^G N_I\rangle = \iota\big(\prod_{b\in\Cb(I)} (\tilde{P}_I \cap \tilde{G}_b)\big),$ 
\item $\displaystyle 
M_I \cap \langle\, ^G N_I\rangle = \iota\big(\prod_{b\in\Cb(I)} (\tilde{M}_I \cap \tilde{G}_b)\big), $
\item $\displaystyle 
M_I\langle\, ^G N_I\rangle = G, $
\item $\displaystyle \langle\, ^G N_I\rangle$   {\it contains}  $N_{I,\op}$.
\end{enumerate}
\vskip2mm

\noindent\textbf{Proof} Parts (i) and (ii) are immediate consequences of the previous considerations. Let us prove (iii). From the proposition $\langle\, ^G N_I\rangle$  contains $\iota(\tilde{G}_b)$ for $b\in \Cb(I)$, but for $b\in \Cb - \Cb(I)$, $M_I$ contains $\iota(\tilde{G}_b)$, so finally $M_I\langle\, ^G N_I\rangle$ contains $\iota(G^{\is})=G'$. Since $M_I$ contains $Z$ and $G=ZG'$, we get (iii). Part (iv) follows from the proposition because $N_{I,\op}$ is $\iota(\prod\limits_{b\in\Cb(I)} (\tilde{N}_{I,\op}\cap\tilde{G}_b))$. $\square$
\medskip

\noindent\textbf{Remark} For $b\in \Cb$, $\tilde{M}_I \cap \tilde{G}_b$ can be also described as the product $\prod\limits_c \tilde{M}_{\Delta_c}$ over the connected components $c$ of the Dynkin diagram obtained from that of $\tilde{G}_b$ by deleting vertices outside~$I$. (We note that the product is not direct.)

\subsection{}\label{II.6} There is another useful characterization of $M_I \cap\langle\, ^G N_I\rangle$.

\vskip2mm

\noindent\textbf{Proposition} 
\textit{Let $I\subset \Delta$. Then $M_I \cap\langle\, ^G N_I\rangle$ is the normal subgroup of $M_I$ generated by $Z\cap M_\alpha'$, for $\alpha$ running through} $\Delta-I$.
\vskip2mm

\noindent\textbf{Proof} Let $\alpha\in \Delta-I$ and let $b\in \Cb$ be such that $\alpha\in\Delta_b$, so that $M_\alpha'\subset \iota(\tilde{G}_b)$. As $\alpha\notin I$, $b$ belongs to $\Cb(I)$ so $\iota(\tilde{G}_b)$ is included in $\langle\, ^G N_I\rangle$ by \ref{II.5} Proposition, and consequently $Z\cap M_\alpha'\subset M_I\cap\langle\, ^G N_I\rangle$. To prove that $M_I\cap\langle\, ^G N_I\rangle$ is the normal subgroup of $M_I$ generated by the $Z\cap M_\alpha'$, $\alpha\in \Delta-I$, it is enough, by \ref{II.5}, to work within $\tilde{G}_b$. So we now assume that $\bg=\bg^{\is}$ and $\bg$ is almost $F$-simple. If $I=\Delta$, $N_I$ is trivial so there is nothing to prove. So let us assume $I\not=\Delta$, so that $\langle\, ^G N_I\rangle=G$ by \ref{II.3} since $N_I$ is not trivial. We can apply to $\bm_I$ all the considerations applied to $\bg$ in the current chapter, so we see that $M_I=Z\prod H_J$ where $J$ runs through connected components of the Dynkin diagram with set of vertices $I$ associated to $M_I$, and $\bh_J$ is the corresponding semisimple simply connected almost $F$-simple subgroup of $\bm_I$. Let $J$ be such a connected component. As the Dynkin diagram attached to $G$ is by assumption connected, there is $\alpha$ in $\Delta-I$ with $\langle J,\alpha^\vee\rangle\not=0$. Choose $\alpha'$ in $J$ with $\langle\alpha',\alpha^\vee\rangle\not=0$ and $x\in F^\times$ with $\alpha'(\alpha^\vee(x))^2\not=1$. We have $\alpha^\vee(x)\in Z\cap M_\alpha'$, $U_{\alpha'} \subset H_J\subset M_I$, and the map from $U_{\alpha'}$ to itself given by $u\mapsto \alpha^\vee(x) u\alpha^\vee (x)^{-1}u^{-1}$ is onto\footnote{If $2\alpha'$ is not a root, then $\alpha^\vee(x)$ acts on $U_{\alpha'}$ (a vector group) via multiplication by $\alpha'(\alpha^\vee(x))$. If $2\alpha'$ is a root, then $\alpha^\vee(x)$ acts on $U_{2\alpha'}$ via $\alpha'(\alpha^\vee(x))^2$ and on $U_{\alpha'}/U_{2\alpha'}$ via $\alpha'(\alpha^\vee(x))$.}. The normal subgroup of $M_I$ generated by $Z\cap M_\alpha'$ contains $\alpha^\vee(x)$\footnote{It follows from \ref{II.4}, footnote 5, that $\alpha^\vee(x)$ belongs to $M_\alpha'$; on the other hand it belongs to $S\subset M_I$.}
and $u\alpha^\vee(x)^{-1}u^{-1}$ for $u\in U_{\alpha'}$, so it contains $U_{\alpha'}$. By \ref{II.3} Proposition it contains $H_J$ and in particular $Z\cap M_{\alpha''}'$ for all $\alpha''\in J$. We conclude that the normal subgroup of $M_I$ generated by the $Z\cap M_\alpha'$, $\alpha\in \Delta-I$, contains $Z\cap M_\alpha'$ for all $\alpha\in\Delta$. By \ref{II.3} Proposition it contains $Z$; since we have seen that it contains each $H_J$, it is equal to $M_I=Z\prod H_J$.~$\square$

\subsection{}\label{II.7} Keeping the same notation, we can now derive consequences for \repr s.

\vskip2mm

\noindent\textbf{Proposition} 
\textit{Let $I\subset \Delta$, and let $\sigma$ be a \repr\ of $M_I$. Then the following conditions are equivalent:}

\textit{(i) $\sigma$ extends to a \repr\ of $G$ trivial on $N_I$,}

\textit{(ii) for each $b\in \Cb(I)$, $\sigma$ is trivial on $\iota(\tilde{M}_I\cap \tilde{G}_b)$,
}

\textit{(iii) for each $\alpha\in \Delta-I$, $\sigma$ is trivial on $Z\cap M_\alpha'$.
}

\textit{\noindent When these conditions are satisfied, there exists a unique  extension $^e\sigma$  of $\sigma$ to $G$  which is trivial on $N_I$, and it  is smooth, admissible or irreducible if and only if $\sigma$~is.}

\vskip2mm 

\noindent\textbf{Proof} As already said in \ref{II.2}, if $\sigma$ extends to a representation of $G$ trivial on $N_I$, the extension is trivial on $\langle\, ^{G}  N_I\rangle$ so $\sigma$ is certainly trivial on $M_I \cap \langle\, ^{G} N_I\rangle$. Consequently, (i) implies (ii) and (iii) by \ref{II.5}, \ref{II.6}. Conversely, under assumptions (ii) or (iii), $\sigma$ is trivial on $M_I\cap \langle\, ^{G}N_I\rangle$ hence extends, trivially on $\langle\, ^{G} N_I\rangle$, to a representation of $M_I\langle\, ^{G} N_I\rangle$, which is $G$ by \ref{II.5} Corollary (iii). The extension $^{e}\sigma$ is necessarily unique. 
Assume that $\sigma$ extends to a representation $^{e}\sigma$ of $G$ trivial on $N_I$. Since $\sigma$ and $^{e}\sigma$ have the same image, $\sigma$ is irreducible if and only if $^{e}\sigma$ is. As $P_I$ is a topological subgroup of $G$, $\sigma$ is smooth if $^{e}\sigma$ is. Conversely, assume that $\sigma$ is smooth and let $x$ be a vector in the space of $\sigma$, $J$ its stabilizer in $P_I$;  by \ref{II.5} Corollary (iv), $N_{I,\op}$ acts trivially on $^{e}\sigma$ and the stabilizer of $x$ in $G$, which contains $N_{I,\op}J$, is open in $G$, so $^{e}\sigma$ is smooth.

 As $P_I$ is a topological subgroup of $G$, $^{e}\sigma$ is admissible if $\sigma$ is. Conversely assume $^{e}\sigma$ is admissible; for each open subgroup $J$ of $M_I$, a vector in $\sigma$ fixed by $J$ is also fixed by the subgroup generated by $J$, $N_I$ and $N_{I,\op}$ which is open in $G$, so $\sigma$ is admissible. $\square$

\vskip2mm

\noindent \textbf{Remark 1}  By \ref{II.5} Remark, condition (ii) illustrates that $\sigma$ can extend to $G$ (trivially on $N_I$) only for very strong reasons: for any connected component $\Delta_b$ of the Dynkin diagram of $G$ meeting $\Delta-I$, $\sigma$ has to be trivial on $M_{\Delta_c}'$ for any connected component $\Delta_c$ of the Dynkin diagram of $M_I$ included in $\Delta_b$. By \ref{II.3} Proposition applied to $M_{\Delta_c}^{\is}$ that last condition is also equivalent to $\sigma$ being trivial on $U_\beta$ for some, or any, $\beta\in \Delta_c$.

\vskip2mm

\noindent \textbf{Remark 2} The coefficient field plays no role here. Properties (i), (ii) and (iii) are equivalent for a representation of $M_I$ over a commutative ring. 
The last assertion of the proposition remains also true for representations over a commutative ring (for admissibility, suppose as usual that the ring is noetherian).

\vskip2mm
\noindent\textbf{Notation} Let $P=MN$ be a parabolic subgroup of $G$ containing $B$, and let $\sigma$ be a representation of $M$. We let $\Delta(\sigma)$ be the set of $\alpha\in \Delta-{\Delta_P}$ such that $\sigma$ is trivial on $Z\cap M_\alpha'$. We let $P(\sigma)$ be the parabolic subgroup corresponding to $\Delta(\sigma)\sqcup\Delta_P$.

\vskip2mm

\noindent \textbf{Corollary 1} \textit{Let $P=MN$ be a parabolic subgroup of $G$ containing $B$, and let $\sigma$ be a representation of $M$. Then the parabolic subgroups of $G$ containing $P$ to which $\sigma$ extends, trivially on $N$, are those contained in $P(\sigma)$. In that case the extension is unique and is smooth, admissible or irreducible if $\sigma$ is. 
}
\vskip2mm

The corollary is immediate from the proposition applied to Levi components of parabolic subgroups of $G$ containing $P$.

\vskip2mm

\noindent \textbf{Remark 2} Since any parabolic subgroup $P$ of $G$ is conjugate to one containing $B$, it follows, as stated in the introduction, that if $\sigma$ is a representation of $P$ trivial on its unipotent radical, there is a maximal parabolic subgroup $P(\sigma)$ of $G$ to which $\sigma$ can be extended, and the extension is smooth, admissible or irreducible if (and only if) $\sigma$ is.

\vskip2mm

\noindent\textbf{Corollary 2} \textit{Keep the assumptions and notation of Corollary $1$, and assume further that $\Delta(\sigma)$ is not orthogonal to $\Delta_M$. Then there is a proper parabolic subgroup $Q$ of $M$, containing $M\cap B$, such that $\sigma$ is trivial on the unipotent radical of $Q$; moreover $\sigma$ is a subrepresentation of $\Ind_Q^{M}(\sigma_{|Q})$, and $\sigma_{|Q}$ is irreducible or admissible if $\sigma$ is. In particular, $\sigma$ cannot be supercuspidal.}

\vskip2mm

\noindent\textbf{Proof} We may assume that $G=P(\sigma)$. Let $\alpha\in \Delta(\sigma)$ not orthogonal to $\Delta_M$, and let $b\in \Cb$ such that $\alpha\in \Delta_b$. Then $\Delta_b\cap \Delta_M \ne \Delta_b$, so $\sigma$ is trivial on $\iota(\tilde{M}\cap \tilde{G}_b)$ by the proposition.
As $\alpha$ is not orthogonal to $\Delta_M$, $\Delta_b\cap \Delta_M$ is not empty. If $Q$ is the (proper) parabolic subgroup of $M$ corresponding to $\Delta_M-\Delta_b$, then $\iota(\tilde{M} \cap \tilde{G}_b)$ contains the unipotent radical $N_Q$ of $Q$ and $\sigma$ is trivial on $N_Q$. Then, obviously, $\sigma$ is a subrepresentation of $\Ind_Q^{M}(\sigma_{|Q})$ and by the proposition, applied to $M$ instead of $G$, if $\sigma$ is irreducible or admissible, so is its restriction to the Levi component of $Q$. By the definition of supercuspidality, $\sigma$ cannot be supercuspidal. $\square$

\vskip2mm

\noindent \textbf{Remark 3} The last assertion of Corollary 2 explains why the case of interest to us is when $\Delta_M$ and $\Delta(\sigma)$ are orthogonal -- an analogous result will be obtained when $\sigma$ is assumed supersingular instead of supercuspidal (\ref{III.17} Corollary). 
As a special case, assume that the (relative) Dynkin diagram of $\bf G$ is connected, and $\sigma$ is a supercuspidal representation of $M$ extending to $G$. Then either $M=G$ or $M=Z$; in the latter case, $\sigma$ is trivial on $Z\cap G'$ and finite dimensional.

\vskip2mm

\noindent \textbf{Remark 4} For the record, let us state a few useful facts when $\Delta$ is the disjoint union of two subsets $I$ and $J$, orthogonal to each other. Then $M_I'$ and $M_J'$ are normal subgroups of $G$, commuting with each other. We have $G'=M_I'M_J'$, $M_I=ZM_I'$, $M_J=ZM_J'$, $M_I\cap M_J=Z$ and in particular $M_I\cap M_\alpha'=Z\cap M_\alpha'$ for $\alpha\in J$. Also, $M_I'\cap M_J'$ is finite and central in $G$: indeed, decomposing $\tilde{G}$ as $\tilde{G}_I\times \tilde{G}_J$, $M_I' \cap M_J'$ is simply the image under $(g_1,g_2) \mapsto \iota(g_1)$ of $\Ker \iota \subset \tilde{G}_I\times \tilde{G}_J$. The inclusion of $M_I$ in $G$ induces an isomorphism $M_I/(M_I \cap M_J' )\simeq G/M_J'$ (and similarly for $M_J$).

\vskip2mm

\noindent \textbf{Remark 5} Let $\alpha\in \Delta$ belong to the component $\Delta_b$. The normal subgroup of $G$ generated by $Z\cap M_\alpha'$ is $\iota(\tilde{G}_b)$ because $Z\cap M_\alpha'$ is not central in $M_\alpha'$. If $\sigma$ is a representation of $G$ which is trivial on $Z\cap M_\alpha'$, it is then trivial on $\iota(\tilde{G}_b)$ and the conclusions of Corollary 2 hold (with $M=G$).

\subsection{}\label{II.8} To go further we need the generalized Steinberg representations already recalled in the introduction.

\vskip2mm

\noindent\textbf{Lemma} \textit{Let $Q$ be a parabolic subgroup of $G$. Then lifting functions on $G$ to functions on $G^{\is}$ via $\iota$ gives an isomorphism of $\Ind_Q^G1$ with $\Ind_{\tilde{Q}}^{G^{\is}}1$. The \repr\ $\St_Q^G \circ \iota$ of $G^{\is}$ is isomorphic to $\St_{\tilde{Q}}^{G^{\is}}$; the restriction of $\St_Q^G$ to $G'$ is irreducible and admissible.}

\vskip2mm

\noindent\textbf{Proof} We have $ZG'=G$ and $Q$ contains $Z$, so $G=QG'$. Besides $Q\cap G'=\iota(\tilde{Q})$. It follows that $\iota$ induces a bijection of $\tilde{Q}\ba G^{\is}$ onto $Q\ba G$; that bijection is continuous hence is a homeomorphism by Arens' theorem \cite[p.\ 65]{MZ}. The first assertion follows and the others are immediate consequences. $\square$

\medskip

Now let $P=MN$ be a parabolic subgroup of $G$, let $\sigma$ be a representation of $M$, inflated to $P$.  Then by \ref{II.7} Corollary 1, $\sigma$ extends (uniquely) to a  \repr\ $^e\sigma$ of $P(\sigma)$. For each parabolic subgroup $Q$ with $P\subset Q\subset P(\sigma)$ we can form the \repr\ $^e\sigma \otimes \St_{Q}^{P (\sigma)}$ of $P(\sigma)$.

\vskip2mm

\noindent\textbf{Proposition} \textit{ $\sigma$ is irreducible ($\resp.$ admissible) if and only if $^e\sigma \otimes \St_{Q}^{P (\sigma)}$ is irreducible ($\resp.$ admissible).} 

\vskip2mm

From this, we get (see for instance \cite[Lemma 4.7]{Vig3}):

\vskip2mm

\noindent \textbf{Corollary} \textit{ $\sigma$ is admissible if and only if $\Ind_{P (\sigma)}^G(^e\sigma \otimes \St_{Q}^{P (\sigma)})$ is admissible.}  

\vskip2mm

\noindent\textbf{Proof of the proposition} The unipotent radical of $P(\sigma)$ acts trivially on both $^e\sigma$ and $\St_{Q}^{P (\sigma)}$. Therefore  we  may assume $P(\sigma)=G$. 

By the lemma above $\St_{Q}^G \circ \iota$ is the generalized Steinberg \repr\ $\St_{\tilde Q}^{G^{\is}}$. For $b$ in $\Cb-\Cb(\Delta_Q)$, $\Delta_Q\cap \Delta_b=\Delta_b$  so that by construction $\St_{\tilde Q}^{ G^{\is}}$ is trivial on $\tilde{G}_b$; consequently,  its restriction to $H=\prod\limits_{b\in \Cb(\Delta_Q)}$ $\tilde{G}_b$ is irreducible. On the other hand by \ref{II.5}, $^e\sigma$ is trivial on the normal subgroup $\iota(H)$. If $\sigma$ is irreducible, the irreducibility of $^e\sigma\otimes \St_{Q}^{{G}}$ comes then from Clifford theory  as  in \cite[Lemma~5.3]{Abe}\footnote{To apply that lemma, note that Schur's lemma is valid for the restriction of $\St_{Q}^{{G}}$ to $\iota(H)$.}.

  Assume that $\sigma$ is admissible, so $^e\sigma$ is admissible too. As above $\iota(H)$ acts trivially on $^e\sigma$ and the restriction of $\St_{Q}^{{G}}$ to $\iota(H)$ is admissible. If $L$ is an open subgroup of $G$, the vectors in $\St_{Q}^{{G}}$ fixed under $L\cap \iota(H)$ form a finite dimensional vector space $X$. The vectors fixed by $L$ in $^e\sigma\otimes \St_{Q}^{{G}}$ are in $^e\sigma \otimes X$. There is an open subgroup $L'$ of $L$ acting trivially on $X$ and $(^e\sigma \otimes X)^{L'}= {}^e \sigma^{L'}\otimes X$ is finite dimensional. Consequently, $^e \sigma \otimes \St_{Q}^{{G}}$ is admissible.
  
 Conversely, if  $^e\sigma \otimes \St_Q^G$ is irreducible,  obviously $\sigma $ is irreducible. If $^e\sigma \otimes \St_Q^G$ is admissible so is $\sigma $. Indeed, if $J$ is an open subgroup of $G$ then $(^e\sigma )^J \otimes (\St_Q^G)^J$ is contained in $(^e\sigma  \otimes \St_Q^G)^J$, so if $J$ is small enough for $(\St_Q^G)^J$ to be non-zero, we deduce that $(^e\sigma)^J$ is finite-dimensional; thus $^e\sigma $ is admissible and so is $\sigma $ by \ref{II.7} Proposition.  $\square$

\medskip
  
\noindent \textbf{Remark} \textit{ Assume that $\Delta_M$ is orthogonal to $\Delta-\Delta_M$.    Let $\sigma $ be a representation of $M$  which extends to $G$ trivially on $N $, and let $Q $  be a parabolic subgroup  of $G$ containing $P$.  }

1)  \textit{The representation $^e{\sigma } \otimes \St_{Q }^G$  of $G$ determines $\sigma$ and $Q$.   }

 2)  \textit{Any subquotient $\pi$ of  $^e{\sigma } \otimes \St_{Q }^G$ is of the form $^e{\sigma_\pi } \otimes \St_{Q }^G$ for some representation $\sigma_\pi $ of $M$  which extends to $G$ trivially on $N $.}

\vskip2mm

\noindent\textbf{Proof} 1) We put $J=\Delta - \Delta_M$.   As  $Q$ contains $M $, $\St_Q^G$ is trivial on the normal subgroup $M'$, and restricting to $M_J$ functions on $G$ gives an isomorphism of $\St_Q^G$ onto $\St_{Q\cap M_J}^{M_J}$.
The restriction of $^e{\sigma} \otimes \St_{Q}^{G}$ to $M_J'$ is a direct sum of irreducible \repr s $\St_{Q}^{G}|_{M_J'}$, and that representation determines $Q$ (\ref{II.8} Lemma).
  Seen as  a representation of $G$, $\Hom_{M_J'} (\St_{Q}^{G}$, $^e{\sigma} \otimes \St_{Q}^{G})$ is isomorphic to $^e\sigma$ (use for example  \cite[Lemma~5.3]{Abe}),  and $^e\sigma $ determines $\sigma $. 

2)  The restriction of $\pi$ to $M_J'$ is a sum of copies of the irreducible \repr\ $\St_Q^G|_{M_J'}$. By Clifford theory \cite[Lemma~5.3]{Abe}, $\pi$ is isomorphic to $\Hom_{M_J'}(\St_Q^G,\pi) \otimes \St_Q^G$. Moreover, $ \Hom_{M_J'}(\St_Q^G,\pi)$ is a  \repr\ of $G$ trivial on $M_J'$  hence determines a representation $\sigma_\pi$ of  $M$ via the map $M \twoheadrightarrow G/M_J'$ and $^e{\sigma _\pi} \simeq \Hom_{M_J'}(\St_Q^G,\pi)$ as a representation of $G$. $\square$

\section{Supersingularity and classification}\label{III}

\subsection{}\label{III.1} This chapter is devoted to the proof of  \ref{I.5} Theorem 4, and is rather long. It is divided into parts A) to H). In  part A) we give some more detail on supersingularity, and in part B) we describe a parametrization for the irreducible \repr s of $K$. The next step in part C) is to determine the weights and eigenvalues of parabolically induced \repr s. We then proceed to the analysis of the \repr s $I(P,\sigma,Q)$: we first determine $P(\sigma)$ in part D), and after that we compute the weights and eigenvalues of $I(P,\sigma,Q)$ for a \supers\ triple $(P,\sigma,Q)$ in part E). The subsequent proof of the irreducibility of $I(P,\sigma,Q)$ in part F) uses a change of weight theorem proved in Chapter~\ref{IV}. From the knowledge of weights and eigenvalues, we easily deduce in part G) when  $I(P_1,\sigma_1,Q_1)$ is isomorphic to $I(P_2,\sigma_2,Q_2)$  for supersingular triples $(P_1,\sigma_1,Q_1)$ and $(P_2,\sigma_2,Q_2)$. In part H) we finally prove exhaustion, i.e.\ that every irreducible admissible representation of $G$ has the form $I(P,\sigma,Q)$ for some supersingular triple $(P,\sigma,Q)$: that uses a result established only in Chapter \ref{V} as a further consequence of the change of weight theorem.

\vskip2mm

\noindent\textbf{Notation}   The special maximal parahoric subgroup $K\subset G$ is fixed throughout; we write $K(1)$ for its pro-$p$-radical and  $H^0$ for $H\cap K$, when $H$ is a subgroup of $G$.  Note that $Z^0$ is the unique parahoric subgroup of $Z$ and that  $Z(1)= Z\cap K(1)$ is   the  unique pro-$p$ Sylow subgroup of $Z^0$.

\vskip2mm
\paragraph{{\large A) Supersingularity}}
\addcontentsline{toc}{section}{\ \ \ A) Supersingularity}

\subsection{}\label{III.2} Consider an irreducible \repr\ $(\rho,V)$ of $K$; it is finite-dimensional and trivial on $K(1)$.  The  classification of such objects  will be recalled in part~B).

We view the intertwining algebra $\Hh_G(V)$ as a Hecke algebra,  the convolution algebra of compactly supported functions $\Phi:G \rightarrow \End_C(V)$ satisfying
$$
\Phi(kgk') = \rho(k)\Phi(g)\rho(k')\quad \mathrm{for}\ g\ \mathrm{in}\ G,\ k\ \mathrm{and}\ k'\ \mathrm{in}\ K.
$$
The convolution operation is given by
$$
(\Phi * \Psi) (g) = \sum_{h\in G/K} \Phi(h)\Psi(h^{-1}g)\quad \mathrm{for}\ \Phi,\Psi\ \mathrm{in} \ \Hh_G(V). \leqno(\mathrm{III}.2.1)
$$
The action on $\ind_K^GV$ is also given by convolution:
$$
(\Phi * f) (g) = \sum_{h\in G/K} \Phi(h)f(h^{-1}g)\quad \mathrm{for}\ f\in \ind_K^GV,\ \Phi\in \Hh_G(V).
 \leqno(\mathrm{III}.2.2)
$$

\subsection{}\label{III.3} We need to recall the structure of $\Hh_G(V)$ and its centre $\cz_G(V)$, as elucidated in 
\cite{HV1}, building on \cite{He1,He2}; note that $\Hh_G(V)$ is commutative in the context of \cite{He1,He2, Abe}.

Let $P=MN$ be a parabolic subgroup of $G$ containing $B$. Then  the space of coinvariants $V_{N^0}$ of $N^0$ in $V$ affords an irreducible \repr\ of $M^0$ (which is the special parahoric subgroup of $M$ corresponding to the special point $\textbf{{x}}_0$). For each   \repr\ $\sigma$ of $M$ on a vector space $W$, Frobenius reciprocity and the equalities $G=KP=PK$, $P^0=M^0N^0$, give a canonical isomorphism:
$$
\Hom_G(\ind_K^G V, \Ind_P^GW) \stackrel{\sim}{\longrightarrow} \Hom_M(\ind_{M^0}^M V_{N^0},W)
\leqno(\mathrm{III}.3.1)
$$
The natural algebra homomorphism $\cs_M^G: \Hh_G(V) \lgr \Hh_M(V_{N^0})$ of \ref{I.5} is given concretely~by
$$
[\cs_M^G(\Psi)(m)] \bar{v} = \sum_{n\in N^0\ba N} \overline{\Psi(nm)(v)}
\ \mathrm{for}\ m\ \mathrm{in}\ M,\ v\ \mathrm{in}\ V,
\leqno(\mathrm{III}.3.2)
$$
where a bar indicates the image in $V_{N^0}$ of a vector in $V$ \cite[Proposition~2.2]{HV2}. Recall that (III.3.1) is $\Hh_G(V)$-linear if we let $\Hh_G(V)$ act on the right-hand side via $\cs_M^G$. Recall also that $\cs_M^G$ is injective \cite[Proposition 4.1]{HV2}.

For varying $P=MN$, the homomorphisms $\cs_M^G$ satisfy obvious transitivity properties, and $\cs_Z^G$ identifies $\Hh_G(V)$ with a  subalgebra of $\Hh_Z(V_{U^0})$ which we now describe.
For a root $\alpha$ in $\Phi=\Phi(\bg,\gs)$, the group homomorphism $|\alpha|: x \mapsto |\alpha(x)|$ from $S$ to $\R_+^\times$ extends uniquely to a group homomorphism $Z \rg \R_+^\times$ trivial on $Z^0$, and we still write $|\alpha|$ for that extension. We write $Z^+$ for the set of $z$ in $Z$ such that $|\alpha|(z) \le 1$ for all $\alpha\in \Delta$. Then by \cite[Proposition 4.2]{HV2} $\Hh_G(V)$ is identified via $\cs_Z^G$ with the subalgebra of $\Hh_Z(V_{U^0})$ consisting of elements supported on $Z^+$. By \cite[1.8 Theorem]{HV1}, $\cz_G(V)$ is the subalgebra $\Hh_G(V) \cap \cz_Z(V_{U^0})$ of $\cz_Z(V_{U^0})$ consisting of elements supported on $Z^+$.

\subsection{}\label{III.4} The group $Z$ normalizes $Z^0$ and its pro-$p$ radical $Z(1)$ and the quotient $Z/Z^0$ is a finitely generated abelian  group. The coinvariant space $V_{U^0}$ is in fact a line, and $Z^0$ acts on it via a character $\psi_V: Z^0 \rg C^\times$ trivial on $Z(1)$: see  part~B), for the difference between the notation $\psi_V$ here and in \cite{HV2}. For $z\in Z$, the coset $Z^0z$ supports a non-zero function in $\Hh_Z(V_{U^0})$ if and only if $z$ normalizes $\psi_V$, and such a function is in $\cz_Z(V_{U^0})$ if and only if $\psi_V(zz'z^{-1}z'^{-1})=1$ for all $z'\in Z$ normalizing $\psi_V$.

\medskip 

\noindent \textbf{Notation}  We let $Z_{\psi_V}$ be the subgroup of $Z$ defined by this last condition. It contains $S$ and~$Z^0$. 

\medskip 

For $z\in Z$ normalizing $\psi_V$ we write $\tau_z\in \Hh_Z(V_{U^0})$ for the function with support $Z^0z$ and value $\id_{V_{U^0}}$ at $z$; we have
$$
\tau_z * \tau_{z'} =\tau_{zz'}\ \mathrm{for}\ z,\ z'\ \mathrm{in}\ Z\ \mathrm{normalizing}\ \psi_V.
$$
Identifying $\Hh_G(V)$ and $\Hh_M(V_{N^0})$ with subalgebras of $\Hh_Z(V_{U^0})$ via $\cs_Z^G$ and $\cs_Z^M$, we can now describe $\Hh_M(V_{N^0})$ as the localization of $\Hh_G(V)$ at some central element \cite[Proposition~4.5]{HV2} (so that $\cz_M(V_{N^0})$ is the localization of $\cz_G(V)$ at the same element).

\vskip2mm

\noindent \textbf{Proposition} \textit{Let $M=M_I$ for some $I\subset \Delta$, and let $s\in S$ satisfy $|\alpha|(s)<1$ for $\alpha\in \Delta-I$, $|\alpha|(s)=1$ for $\alpha\in I$. Then $\Hh_M(V_{N^0})$ is the localization of $\Hh_G(V)$ at $\tau_s$, and $\cz_M(V_{N^0})$ the localization of $\cz_G(V)$ at~$\tau_s$.}

\vskip2mm

\noindent \textbf{Notation}  For each $\alpha\in \Delta$, we choose $z_\alpha$ in $S$ such that $|\alpha|(z_\alpha)<1$ and  $|\alpha'|(z_\alpha)=1$ for $\alpha'\in \Delta-\{\alpha\}$.
For  a character  $\chi$ of $\cz_G(V)$, we let $\Delta_0(\chi)=\{\alpha\in \Delta\mid \chi(\tau_{z_\alpha})=0\}$. 
\vskip2mm

In the above proposition, we can take $s=\prod\limits_{\alpha\in \Delta-I}
z_\alpha$; then  $\tau_s$ is the product $\tau_s=\prod\limits_{\alpha\in \Delta-I} \tau_{z_\alpha}$ in any order.
\vskip2mm

\noindent \textbf{Lemma} \textit{Let $\chi$ be a character of $\cz_G(V)$. Then $I = \Delta_0(\chi)$  is the smallest subset  of $\Delta$ such that $\chi$ extends to a character of $\cz_{M_I}(V_{N_I^0})$. For $z$ in $Z^+ \cap Z_{\psi_V}$ we have $\chi(\tau_z)\not=0$ if and only if $|\alpha|(z)=1$ for all $\alpha\in \Delta_0(\chi)$. In particular, $\Delta_0(\chi)$ does not depend on $\{z_\alpha\}$.}

\vskip2mm

\noindent\textbf{Proof} As $\cz_{M_I}(V_{N_I^0})$ is the localization of $\cz_G(V)$ at $\prod\limits_{\alpha\in \Delta-I}\tau_{z_\alpha}$, $\chi$ extends to a character of $\cz_{M_I}(V_{N_I^0})$ if and only if $\chi(\tau_{z_\alpha})\not=0$ for $\alpha\in \Delta-I$.  The first assertion follows. Let $z\in Z^+\cap Z_{\psi_v}$; if for some $\alpha\in \Delta_0(\chi)$ we have $|\alpha|(z)<1$, then for some positive integer $r$, $z^r=z_\alpha t$ with $t\in Z^+ \cap Z_{\psi_V}$, and $\chi(\tau_z)^r=\chi(\tau_{z_\alpha})\chi(t)=0$, so $\chi(\tau_z)=0$; if $|\alpha|(z)=1$ for all $\alpha\in\Delta_0(\chi)$ then with $s=\prod\limits_{\alpha\in \Delta-\Delta_0(\chi)}z_\alpha$ there is a positive integer $x$ such that $s^r=zt$ for some $t\in Z^+\cap Z_{\psi_V}$ and similarly $\chi(\tau_z)\not=0$ since $\chi(\tau_s)\not=0$. $\square$

\vskip2mm  We write $Z_{\Delta}^\bot $ for the set of $z\in Z$ with $|\alpha|(z)=1$ for all $\alpha\in \Delta$. 
Using the lemma, we can restate the definition of \supers ity (\ref{I.5}) for a character of $\cz_G(V)$.

\vskip2mm

\noindent \textbf{Corollary} \textit{For a  character $\chi$ of $\cz_G(V)$, the following conditions are equivalent:}
 
 (i) \textit{$\chi$  is \supers ,}

(ii)  $\Delta_0(\chi)=\Delta$,

(iii) \textit{$\chi(\tau_z)=0$ for all $z$ in $Z^+ \cap Z_{\psi_V}$ not in $Z_{\Delta}^\bot$.}
  
\vskip2mm

\paragraph{{\large B) Irreducible representations of $K$}}
\addcontentsline{toc}{section}{\ \ \ B)  Irreducible representations of \texorpdfstring{$K$}{K}}

\subsection{}\label{III.5} For a subgroup $H\subset G$ we put $\overline{H}=(H\cap K)/(H\cap K(1))$.
As recalled above, irreducible  representations of $K$ factor through $\overline{K}= K/K(1)$. Information about $\overline{K}$ comes from \cite{BT1,BT2}, see also \cite{Ti}. The group $\overline{K}$ is naturally the group of points (over the residue field $k$ of $F$), of a connected reductive group, which we write $\bg_k$, so that $\overline{G}= \overline{K}=\bg_k(k)$\footnote{We warn the reader that when $\bg$ is semisimple, $\bg_k$ is not necessarily semisimple. If $\bh_k$ is an algebraic group over $k$, we put $H_k=\bh_k(k)$, so that for many algebraic subgroups $\bh$ of $\bg$ in the current chapter, we can use indifferently the notations $\overline{H}$ or $H_k$ for $(H\cap K)/(H\cap K(1))$ -- we mostly  use $\overline{H}$.}. 
We also have $\overline{S}=\gs_k(k)$, where $\gs_k$ is a maximal split torus in $\bg_k$, with a natural identification of $X^*(\gs_k)$ and $X^*(\gs)$; if $\bz_k$ is the centralizer of $\gs_k$ in $\bg_k$ then $\overline{Z}=\bz_k(k)$, and similarly for the normalizer $\boldsymbol{\cn}_k$ of $\gs_k$ in $\bg_k$. As $K$ is a special parahoric subgroup, every element of $W_0$ has a  representative in $K$ so that $W_0=\cn^0/Z^0$, and reduction mod $K(1)$ yields an identification of $W_0$ with $W(\bg_k,\gs_k)=\overline{\cn}/\overline{Z}$.

Similarly $\overline{B}=\bb_k(k)$ for a minimal parabolic subgroup $\bb_k$ of $\bg_k$ with Levi component $\bz_k$ (which is a torus since $k$ is finite) and unipotent radical $\bu_k$ such that $\overline{U}=\bu_k(k)$.

\subsection{}\label{III.6} The root system of $\gs_k$ in $\bg_k$ is a sub-root system of the root system of $\gs$ in $\bg$, using the above-mentioned identification of $X^*(\gs_k)$ and $X^*(\gs)$. We write $\Phi_k$ for the set of roots of $\gs_k$ in $\bg_k$; we have $\Phi_k\subset \Phi$. A reduced root $\alpha\in \Phi$ belongs to $\Phi_k$ if $2\alpha$ is not a root in $\Phi$; if $\alpha$ and $2\alpha$ are roots in $\Phi$, then $\alpha$ or $2\alpha$ or both are in $\Phi_k$ -- all three cases can occur.

So we get a natural bijection $\alpha\mapsto \overline{\alpha}$ from reduced roots in $\Phi$ to reduced roots in $\Phi_k$, which sends positive roots to positive roots, and the set $\Delta$ of simple roots in $\Phi$ to the set $\Delta_k$ of simple roots in  $\Phi_k$. When  $\alpha\in \Phi$ is reduced, we have $\overline{U}_\alpha=\bu_{k,\overline{\alpha}}(k)$. Henceforward we \textbf{identify} the reduced roots of $\Phi_k$ with those of $\Phi$, hence
 $\Phi_k$ with $\Phi$, $\Delta_k$ with $\Delta$, via $\alpha\mapsto\overline{\alpha}$. Then for $I\subset \Delta$ the parabolic subgroup $P_I=M_IN_I$ is such that $\overline{P_I}=\bp_{I,k}(k)$, $\overline{M_I}=\bm_{I,k}(k)$, $\overline{N_I}=\bn_{I,k}(k)$.

\subsection{}\label{III.7} Let $\bb_{\op}$ be the parabolic subgroup of $\bg$ opposite to $\bb$\footnote{When convenient, we put the index $\op$ on top.} (with respect to $\bz$) and $\bu_{\op}$ its unipotent radical; then $\overline B_{\op}= \bb_{k,\op}(k)$ where $ \bb_{k,\op}$ is the parabolic subgroup of $\bg_k$ opposite to $\bb_k$. Similarly we have $\overline U_{\op}= \bu_{k,\op}(k)$ for their unipotent radicals.

From \cite[Proposition 6.25]{BoT} we get that $\overline{G}$ is generated by the union of $\overline{Z}$, $\overline{U}$, $\overline U_{\op}$. The subgroup $\overline{G}'$ of $\overline{G}$ generated by the union of $\overline{U}$ and $\overline U_{\op}$ is normal in $\overline{G}$; it is the image in $\overline{G}$ of $\bg_{k,\simplyc}(k)$ where $\bg_{k,\simplyc}$ is the simply connected covering of the derived group of $\bg_k$. Note\footnote{Recall $G'$ is the subgroup of $G$ generated by $U$ and $U_{\op}$.} that $G'{}^{0}$ certainly contains $U^0$ and $(U_{\op})^0$ so that its image in $\overline{G}$ contains $\overline{G}'$. But it can be larger, so we need to distinguish $\overline{G}'$ and $\overline{G'}$\footnote{To avoid confusion, we sometimes write $G_k'$ rather than $\overline{G}'$.}; the discrepancy is actually quite important in Chapter~\ref{IV}.
\vskip2mm

\noindent \textbf{Lemma} \textit{(i) The map $(U\cap K(1))\times Z(1)\times (U_{\op}\cap K(1)) \rg K(1)$ given by the product law is bijective, and similarly for any order of the factors.}

\noindent\textit{(ii) $K$ is generated by the union of $U^0$, $Z^0$ and $(U_{\op})^0$.}

\vskip2mm

\noindent\textbf{Proof} Assertion (i) is due to Bruhat and Tits \cite[4.6.8 Corollaire]{BT2}. Since $\overline{G}$ is generated by the union of $\overline{Z}$, $\overline{U}$ and $\overline U_{\op}$, $K$ is generated by the union of $Z^0$, $U^0$, $(U_{\op})^0$ and the normal subgroup $K(1)$; then (ii) follows from~(i). $\square$

\vskip2mm

The lemma has a consequence which will be useful later. As in \ref{III.4} we write $Z_{\Delta}^\bot$ for the set of $z\in Z$ such that $|\alpha|(z)=1$ for all $\alpha\in \Delta$. Equivalently, $Z_{\Delta}^\bot = \Ker v_Z$ in the notation of  \cite[3.2]{HV1}. (We have in fact that $|\alpha|(z)=q^{-\langle\alpha,v_Z(z)\rangle}$ for $\alpha\in \Delta$ and $z \in Z$.)

\vskip2mm

\noindent \textbf{Corollary} \textit{$Z_{\Delta}^\bot$ is the normalizer of $K$ in $Z$.}

\vskip2mm

\noindent\textbf{Proof} If $z\in Z$ normalizes $K$ it also normalizes $U_\alpha^0$ for all $\alpha\in \Phi$. Given the action of $z$ on the filtration of $U_\alpha$ \cite{Ti}, that is equivalent to $|\alpha|(z)=1$ for $\alpha\in \Phi$. Conversely if $|\alpha|(z)=1$ for $\alpha$ in $\Delta$ then $|\alpha|(z)=1$ for all $\alpha$ in $\Phi$ and $z$ normalizes $U_\alpha^0$ for all $\alpha\in \Phi$; it then normalizes $U^0$ and $(U_{\op})^0$, so it normalizes~$K$. That proves that $Z_{\Delta}^\bot$ is the normalizer of $K$ in $Z$. $\square$

\vskip2mm

\noindent \textbf{Remark} By the Cartan decomposition the normalizer of $K$ in $G$ is $Z_{\Delta}^\bot K$.

\subsection{}\label{III.8} We can now recall (see \cite{HV1,HV2} and the references therein) the parametrization of the irreducible representations of $\overline{G}$, up to isomorphism. 

If $(\rho,V)$ is an irreducible representation of $\overline{G}$, then $V^{\overline{U}}$ is a line, on which $\overline{Z}$ acts via a character, say $\eta: \overline{Z} \rg C^\times$. Let $\Delta(\eta)$ be the set of simple roots $\alpha\in \Delta$ such that $\eta$ is trivial on $\overline{Z} \cap {M}_{\alpha,k}'$ (where ${M}_{\alpha,k}$ is the Levi subgroup of $\overline{G}$ corresponding to $\{\alpha\}$), and as in \ref{III.7} ${M}_{\alpha,k}'$ is the subgroup of ${M}_{\alpha,k}$ generated by (the union of) $\overline{U}_\alpha$ and $\overline{U}_{-\alpha}$. The stabilizer of the line $V^{\overline{U}}$ in $\overline{G}$ is a parabolic subgroup containing $\overline B$ corresponding to a subset $\Delta_V$ of $\Delta(\eta)$, and $V$ is characterized up to isomorphism by the pair $(\eta,\Delta_V)$; all such pairs do occur. In \cite{HV2}, $(\eta,\Delta_V)$ is called the \textbf{standard parameter} of~$V$.

\subsection{}\label{III.9} In this paper, we are interested in coinvariants rather than invariants, so we use different parameters. Let $V$ be an irreducible representation of $\overline{G}$ with standard parameter $(\eta,\Delta_V)$.  
\vskip2mm

\noindent \textbf{Lemma} \textit{The group $\overline{Z}$ acts on the line  $V_{\overline{U}}$  via the character $\eta \circ w_0$ where $w_0$ is the longest element in $W_0$. Moreover the stabilizer of the kernel of $V\rg V_{\overline{U}}$ is the parabolic subgroup of $\overline{G}$ corresponding to the subset $-w_0\Delta_V$ of~$\Delta$.}

\vskip2mm

\noindent\textbf{Proof} By \cite[Proposition 3.14]{HV2} the projection $V\rg V_{\overline{U}}$ induces a $\overline{Z}$-equivariant isomorphism of $V^{\overline{U}_{\op}}$ onto $V_{\overline{U}}$; the first assertion  comes from \cite[3.11]{HV2}. The stabilizer we look at is also the stabilizer of the line $(V^*)^{\overline{U}}$ in the contragredient representation $V^*$ of $V$; the second assertion follows from  by \cite[3.12]{HV2}. $\square$

\vskip2mm

\noindent \textbf{Definition} \textit{ The \textbf{parameter} of $V$ is the pair $(\psi_V,\Delta(V))$ where $\overline{Z}$ acts on $V_{\overline{U}}$ via $\psi_V$ and the stabilizer in $\overline{G}$ of the kernel of $V\rg V_{\overline{U}}$ is $\overline{P}_{\Delta(V)}$.}

\vskip2mm

\noindent\textbf{Remarks} 1)  We have  $\psi_V=\eta \circ w_0$ and $\Delta(V)=-w_0\Delta_V$.

\noindent 2) The antistandard parameter of $V$ \cite[3.11]{HV2} is $(\psi_V,-\Delta(V))$.

\noindent 3) $V$ is determined up to isomorphism by its parameter. One has $\Delta(V) \subset \Delta(\psi_V)$, and all pairs $(\psi,I)$ with $I\subset \Delta(\psi)$ occur as parameters.

\subsection{}\label{III.10}

\noindent \textbf{Lemma} \textit{Let $V$ be an irreducible representation of $K$, and let  $ P= MN$ be a  parabolic subgroup  of $G$ containing $B$}.
\begin{enumerate}
\item $V_{\overline{N}}$ \textit{is an irreducible \repr\ of $\overline{M}$ with parameter} $(\psi_V,\Delta_M\cap \Delta(V))$.
\item \textit{$V$ is $\overline P_{\op}$-regular in the sense of} \cite[Def.~3.6]{HV2} \textit{if and only if} $\Delta(V) \subset \Delta_M$.
\end{enumerate}

Here, $\overline P_{\op}=\overline{M}\,\overline N_{\op}$ is the parabolic subgroup of $\overline{G}$ opposite to $\overline{P}$ (relative to $\overline{M}$).

\vskip2mm

\noindent\textbf{Proof} By \cite[3.11]{HV2} $V^{\overline N_{\op}}$ is an irreducible \repr\ of $\overline{M}$ and its antistandard parameters are $(\psi_V,-(\Delta_M\cap \Delta(V)))$. On the other hand, the projection $V\rg V_{\overline{N}}$ induces an $\overline{M}$-equivariant isomorphism of $V^{\overline N_{\op}}$ onto $V_{\overline{N}}$, so (i) comes from Remark 2) above. By \cite[Def.~3.6]{HV2} $V$ is $\overline P_{\op}$-regular if and only if $-\Delta(V) \subset- \Delta_M$ i.e.\ $\Delta(V)\subset \Delta_M$, whence~(ii). $\square$

\vskip2mm

\noindent\textbf{Remarks} 1) Since $\overline P_{\Delta(V)}$ is the stabilizer of the kernel of the projection $V\rg V_{\overline{U}}$, $V$ is one-dimensional if and only if $\Delta(V)=\Delta$. It follows  from  part (i) of the lemma  that $V_{\overline{N}}$ is one-dimensional if and only if $\Delta_M\subset \Delta(V)$. That provides a useful characterization of~$\Delta(V)$.

2) In this paper we will not use the terminology of a weight $V$ being $\overline P_{\op}$-regular. We will phrase everything in terms of the equivalent
condition $\Delta(V) \subset \Delta_M$ of the above lemma.

\medskip

\noindent \textbf{Examples} 1) Consider the case where $V$ is the trivial representation of $\overline{G}$. Then $\psi_V=1$ and $\Delta(V)=\Delta$. Representations $V$ with parameter $(1,I)$ for $I\subset \Delta$ are particularly important to us (cf.~\ref{III.18} below).

2) Let $\eta$ be a character of $\overline{Z}$; then $\eta$ extends to a character of $\overline M_{\Delta(\eta)}$: indeed, that extension is the irreducible representation of $\overline M_{\Delta(\eta)}$ with parameter $(\eta, {\Delta(\eta)})$.

\subsection{}\label{III.11} Consider the simply connected covering $\bg_{k,\simplyc}$ of the derived group $\bg_{k,\der}$ of $\bg_k$  and write $j:\bg_{k,\simplyc}\rg \bg_k$ for the natural morphism. Put ${G}_{k,\simplyc}=\bg_{k,\simplyc}(k)$. We can repeat exactly the same considerations as in \ref{II.4} in this context of finite reductive groups, and we use the analogous notation -- note however that since $k$ is finite, every   almost $k$-simple component of $\bg_{k,\simplyc}$ is isotropic. In particular $j$ induces an isomorphism between $\tilde{\bu}_k$ and $\bu_k$, and $\Delta_k$ also appears as the set of simple roots of $\tilde{\gs}_k$ in~$\tilde{\bu}_k$.

From \ref{III.7}, recall that 
$${G}_k'= j({G}_{k,\simplyc}).$$

\noindent \textbf{Proposition} \textit{Let $(\rho,V)$ be an irreducible representation of $G_k$ with parameter $(\psi_V,\Delta(V))$. Then $(\rho\circ j,V)$ is an irreducible \repr\ of $G_{k,\simplyc}$ with parameter $(\psi_V \circ j_{|\tilde{Z}_k},\Delta(V))$.}

\vskip2mm

Here, $\tilde{Z}_k = \tilde{\bz}_k(k)$, where $\tilde{\bz}_k$ is the centralizer of $\tilde{\gs}_k$ in $\bg_{k,\simplyc}$; we use similarly abbreviated notation below.
By the fact above  and  the inclusion  ${G}_k'=\overline{G}' \subset  \overline{G'}$ (\ref{III.7}), we get:

\vskip2mm

\noindent \textbf{Corollary} \textit{The restriction of $\rho$ to ${G}_k' $, and a fortiori to $\overline{G'}$, is irreducible.}

\vskip2mm

\noindent \textbf{Proof of the proposition} Since $V_{\tilde{U}_k}$, equal to $V_{\overline{U}}$, is one-dimensional, the cosocle of $\rho\circ j$ is irreducible. Similarly $V^{\tilde{U}_{k,\op}}$ equal to $V^{\overline U_{\op}}$ is one dimensional, so the socle of $\rho\circ j$ is irreducible too. As the projection of $V^{\overline U_{\op}}$ to $V_{\overline{U}}$ is non-zero, the map from the socle of $\rho\circ j$ to its cosocle is non-zero, and $\rho\circ j$ is indeed irreducible. Clearly $\tilde{Z}_k$ acts on $V_{\tilde{U}_k}=V_{\overline{U}}$ by $z \mapsto \psi_V \circ j(z)$, and $\tilde{P}_{\Delta(V),k}= j^{-1}(\overline P_{\Delta(V)})$ stabilizes the kernel of $V\rg V_{\tilde{U}_k}$. But for $I\subset \Delta$, we have $\overline{P_I}=\overline{Z} j(\tilde{P}_{I,k})$, so if $\tilde{P}_{I,k}$ stabilizes that kernel, $I\subset \Delta(V)$. $\square$

\medskip

\paragraph{{\large C) Weights of parabolically induced representations}}
\addcontentsline{toc}{section}{\ \ \ C)  Weights of parabolically induced representations}

\subsection{}\label{III.12} Let $P=MN$ be a parabolic subgroup of $G$ containing $B$, and $(\tau,W)$ a \repr\ of $M$. We investigate the weights of $\Ind_P^GW$ and the corresponding Hecke eigenvalues. From now on, we identify the irreducible \repr s of $K$ and those of $\overline{G}= K/K(1)$. 
\vskip2mm

In this part C) we let $(\rho,V)$ be an irreducible \repr\ of $K$, with parameter $(\psi_V,\Delta(V))$. Recall that if $(\pi,X)$ is a \repr\ of $G$, for example $X=\Ind_P^GW$, then $\Hom_K(V,X)$ is a right $\Hh_G(V)$-module  via Frobenius reciprocity. The formula for the action~is
$$
(\varphi \Phi)(v)= \sum_{g\in G/K} g\varphi(\Phi(g^{-1})v)\quad \mathrm{for}\ v\in V,\ \varphi \in \Hom_K(V,X),
\leqno(\mathrm{III.12.1})
$$
and $\Phi\in \Hh_G(V)$.

\vskip2mm

\noindent \textbf{Proposition}\textit{
(i) The natural isomorphism }
$$
\Hom_K(V,\Ind_P^GW) \stackrel{\mathrm{can}\atop\sim}{\lgr} \Hom_{M^0}(V_{N^0},W)
$$
\textit{is $\Hh_G(V)$-linear, where $\Hh_G(V)$ acts on the right-hand side via $\cs_M^G$.}

\textit{(ii) $V$ is a weight for $\Ind_P^GW$ if and only if $V_{N^0}$ is a weight for~$W$.}

\textit{(iii) The map $\cs_M^G$ identifies the eigenvalues of $V$ in $\Ind_P^GW$  and the   eigenvalues  of $V_{N^0}$ in $W$.}

\vskip2mm

\noindent\textbf{Proof} (i) comes from \ref{III.3} and (ii) is an immediate consequence.
We have seen that $\cz_M(V_{N^0})$ is the localization of $\cz_G(V)$ at some element $\tau_s$. Clearly $\tau_s$ acts invertibly on $\Hom_{M^0}(V_{N^0},W)$; as the canonical isomorphism is $\Hh_G(V)$-linear, $\tau_s$ also acts invertibly on  $\Hom_K(V,\Ind_P^GW)$, which gives (iii). $\square$
\vskip2mm
A useful consequence of  (III.12.1) is the following lemma. Recall  that for $z\in Z^+\cap Z_{\psi_V}$, $\cz_G(V)$ contains a unique element $T_z$ such that $\Supp T_z=KzK$ and $T_z(z)\in \End_C(V)$ induces the identity on $V^{\overline{U}_{\op}}$ 
\cite[7.3, 2.9]{HV1}.
\vskip2mm

\noindent \textbf{Lemma} \textit{Let $(\pi,X)$ be a \repr\ of $G$ and $\varphi \in \Hom_K(V,X)$. Let $z\in Z_{\psi_V}$. Assume $z\in Z_{\Delta}^\bot$, i.e.\ that $z$ normalizes $K$. Then $\cs_Z^G(T_z)=\tau_z$ and $(\varphi\tau_z)(v)=z^{-1}\varphi(v)$ for $v$ in $V^{\overline{U}_{\op}}$. If $\varphi$ is an eigenvector for $\cz_G(V)$ with eigenvalue $\chi$, then $z^{-1}$ acts on $\varphi(V^{\overline{U}_{\op}})$  by $\chi(\tau_z)$. }

\vskip2mm

\noindent\textbf{Proof} By assumption $zK=Kz$, and  the endomorphism  $T_z(z)$   satisfies $\rho(k)T_z(z)=T_z(z)\rho(z^{-1}kz)$ for $k\in K$ \cite[7.3]{HV1}. As $z$ normalizes $U^0$ and $(U_{\op})^0$,  $T_z(z)$ induces endomorphisms of  $V_{\overline U}$ and $V^{\overline U_{\op}}$; since the natural map $V^{\overline U_{\op}}\to V_{\overline U}$ is an isomorphism, $T_z(z)$ induces the identity on $V_{\overline{U}}$. From (III.3.2) we get $\cs_Z^G(T_z)=\tau_z$, and (III.12.1) gives
$$
(\varphi T_z)(v)=z^{-1}\varphi(T_z(z)v)\quad \mathrm{for}\ v\in V,
$$  hence the result. $\square$

\subsection{}\label{III.13} Let $\varphi \in \Hom_K(V,\Ind_P^GW)$ and 
$\varphi_M \in \Hom_{M^0}(V_{N^0},W)$ correspond via (III.3.1). 
Then $\varphi$ gives rise to a $G$-morphism, again written $\varphi$, 
 from $\ind_K^GV$ to $\Ind_P^GW$, and similarly we get an $M$-morphism 
 $\varphi_M:\ind_{M^0}^M V_{N^0}\rg W$. 

Consider the following diagram, where horizontal maps are canonical isomorphisms
\begin{equation*}
\begin{CD}
\Hom_G(\ind_K^G V,\ \Ind_P^G(\ind_{M^0}^M V_{N^0})) @>{\mathrm{can}\atop\sim} >> \Hom_M(\ind_{M^0}^M V_{N^0},\ind_{M^0}^M V_{N^0})\\
@VV{\Ind_P^G \varphi_M}V @VV{\varphi_M}V\\
\Hom_G(\ind_K^G V,\ \Ind_P^GW)@>{\mathrm{can}\atop\sim} >> \Hom_M(\ind_{M^0}^M V_{N^0},W)
\end{CD}
\end{equation*}

By naturality, the vertical maps obtained by composing with $\Ind_P^G\varphi_M$ and $\varphi_M$, as indicated, make the diagram commutative. 
  The identity map of $\ind_{M^0}^M V_{N^0}$ yields the canonical intertwiner
$$
\Ii: \ind_K^G V \longrightarrow \Ind_P^G( \ind_{M^0}^M V_{N^0})
\leqno(\mathrm{III}.13.1)
$$
mentioned in \ref{I.6}.  We get:

 \vskip2mm
\noindent\textbf{Lemma}  $\Ind_P^G \varphi_M \circ \Ii =\varphi$.
\vskip2mm

By \cite[Proposition 4.1]{HV2}, $\Ii$ is injective. 
As $\Ii$ is $\Hh_G(V)$-linear, it factors as follows:
\begin{eqnarray*}
\ind_K^G V \lgr \cz_M(V_{N^0})\otimes_{\cz_G(V)}\ind_K^G V \stackrel{u}{\lgr} \Hh_M(V_{N^0}) \otimes_{\Hh_G(V)} \ind_K^G V\\ 
\qquad\qquad\stackrel{\Theta}{\lgr} \Ind_P^G(\ind_{M^0}^M V_{N^0}),
\end{eqnarray*}
for some canonical map  $\Theta$.
Since $\Hh_M(V_{N^0})$ is the localization of $\Hh_G(V)$ at some central element, and $\cz_M(V_{N^0})$ is the localization of $\cz_G(V)$ at the same element, the map $u$ is an isomorphism. 

\subsection{}\label{III.14}  The main result of \cite{HV2} is, taking into account \ref{III.10} Lemma~(ii):

\vskip2mm

\noindent \textbf{Theorem} \textit{Let   $(\psi_V,\Delta(V))$ be the parameter of $V$. If $\Delta(V)\subset \Delta_M$ then the map 
$$\Hh_M(V_{N^0}) \otimes_{\Hh_G(V)} \ind_K^G V \xrightarrow{\Theta}  \Ind_P^G(\ind_{M^0}^M V_{N^0})$$ of {\upshape \ref{III.13}} is an isomorphism.}

\vskip2mm

We derive some consequences.

\vskip2mm

\noindent \textbf{Corollary 1} \textit{ Let $\varphi\in \Hom_K(V,\Ind_P^GW)$ be an eigenvector for $\cz_G(V)$.  If $\Delta(V)\subset\Delta_M$ and if $\varphi_M(V_{N^0})$ generates $W$ as a \repr\ of $M$, then $\varphi(V)$ generates $\Ind_P^GW$ as a \repr\ of~$G$.}

\vskip2mm

\noindent\textbf{Proof} By the theorem, $\Theta$ is surjective. By hypothesis $\varphi_M:\ind_{M^0}^M V_{N^0}\rg W$ 
is surjective, 
 so by  \ref{III.13} Lemma the map induced by~$\varphi$
$$
\cz_M(V_{N^0}) \otimes_{\cz_G(V)}\ind_K^G V \lgr \Ind_P^G W
$$
is surjective. But $\cz_G(V)$ acts on $\varphi$ via a character which extends to $\cz_M(V_{N^0})$ (\ref{III.12} Proposition (ii)) so we conclude that $\varphi(\ind_K^G V) = \Ind_P^G W$, hence the result. $\square$

\vskip2mm

\noindent \textbf{Corollary 2} \textit{Assume that $(\tau,W)$ is irreducible and admissible. Then $\Ind_P^G W$ is irreducible if and only if every non-zero sub\repr\ of it contains a weight $V$ with $\Delta(V) \subset \Delta_M$.}

\vskip2mm

\noindent\textbf{Proof} Since $W$ has some weight, by \ref{III.12} Proposition (i) and \ref{III.10} Lemma~(i), $\Ind_P^G W$ has a weight $V$ with $\Delta(V)\subset \Delta_M$. Conversely if a sub\repr\ $X$ of $\Ind_P^G W$ contains a weight $V$ with $\Delta(V)\subset \Delta_M$, there is an eigenvector $\varphi\in \Hom_K(V,X)$ for $\cz_G(V)$. As $\tau$ is irreducible, $\varphi_M(V_{N^0})$ generates $W$ and by the proposition $X=\Ind_P^G W$. $\square$

\paragraph{{\large D) Determination of $P(\sigma)$ for supersingular $\sigma$}}
\addcontentsline{toc}{section}{\ \ \ D) Determination of \texorpdfstring{$P(\sigma)$}{P(sigma)} for supersingular \texorpdfstring{$\sigma$}{sigma}}

\subsection{}\label{III.15} We want to apply the preceding corollary to prove the irreducibility of $I(P,\sigma,Q)$ for a supersingular triple $(P,\sigma,Q)$. That can only be done in stages. First we  determine $P(\sigma)$ in terms of weights and eigenvalues of $\sigma$. In other words, we determine the set $\Delta(\sigma)$ of $\alpha\in \Delta-\Delta_M$ such that $\sigma$ is trivial on $Z\cap M_\alpha'$~(\ref{II.7}).

As the generality will be useful in Chapter \ref{V}, we consider the situation where $P=MN$ is a parabolic subgroup of $G$ containing $B$, and $(\sigma,W)$ is a representation of $M$ satisfying the following hypothesis:

\vskip2mm

(H) There is an irreducible \repr\ $(\rho,V)$ of $M^0$ and  some $\varphi$ in  $\Hom_{M^0}(V,W)$ such that $\sigma$ is generated by $\varphi(V)$ as a \repr\ of~$M$. 
\vskip2mm

Hypothesis (H) is certainly true if $\sigma$ is irreducible and admissible, and then we can take $\varphi$ to be an eigenvector for $\cz_M(V)$, and the corresponding eigenvalue is supersingular if $\sigma$ is. 
As before, write $(\psi_V,\Delta(V))$ for the parameter of~$V$.

\vskip2mm

\noindent \textbf{Lemma}  \textit{Assume Hypothesis (H). Let $\alpha\in \Delta$. If  $\sigma$ is trivial on  $Z\cap M_\alpha'$, then $\psi_V$ is trivial on  $Z^0 \cap M_\alpha'$.} 
\vskip2mm

\noindent\textbf{Proof} 
If $\sigma$ is trivial on $Z\cap M_\alpha'$, then certainly $Z\cap M_\alpha'$ acts trivially on $\varphi(V)$. As $\varphi\in \Hom_{M^0}(V,W)$ is injective, $Z^0\cap M_\alpha'$ acts trivially on $V$ hence on $V_{U^0}$ and $\psi_V$ is trivial on $Z^0\cap M_\alpha'$. $\square$

\subsection{}\label{III.16}
\noindent\textbf{Proposition} \textit{Let }$\alpha\in \Delta$.

\textit{(i) If $\psi_V$ is trivial on $Z^0 \cap M_\alpha'$ then $Z\cap M_\alpha' \subset Z_{\psi_V}
$.}

\textit{(ii) $|\alpha|$ $(Z\cap M_\alpha')$ is isomorphic to $\Z$.}

\textit{(iii) Let $z\in Z\cap M_\alpha'$. Then $|\alpha|(z)=1$ if and only if $z\in Z^0 \cap M_\alpha'$.}

\vskip2mm

\noindent \textbf{Notation} By the proposition the group $( Z\cap M_\alpha')/(Z^0 \cap M_\alpha')$ is isomorphic to $\Z$.  By (ii) there is an element $a_\alpha$ in $Z\cap M_\alpha'$ with $|\alpha|(a_\alpha)>1$, such that $|\alpha|(a_\alpha)$ generates $|\alpha|(Z\cap M_\alpha')$; by (iii) the element $a_\alpha$ is well-defined modulo $Z^0\cap M_\alpha'$. Note that if $\alpha$ is orthogonal to $\Delta_M$ then $a_\alpha \in Z_{\Delta_M}^\bot$ (see proof of \ref{III.7} Corollary) and $\tau_{a_\alpha}$ is a unit of $\cz_M(V)$. If $\psi_V$ is trivial on $Z^0\cap M_\alpha'$ the element $\tau_{a_\alpha}$ of $\cz_Z(V_{U\cap M^0})$ does not depend on the choice of $a_\alpha$, so we write it $\tau_\alpha$.

\vskip2mm

\noindent\textbf{Proof of the proposition}
Assume that $\psi_V$ is trivial on $Z^0\cap M_\alpha'$, and take $z\in Z\cap M_\alpha'$; then, for $z'\in Z$ (in particular for $z'\in Z^0$), $zz'z^{-1}z'^{-1}$ belongs to $Z^0\cap M_\alpha'$ (because $Z/Z^0$ is abelian and $Z\cap M_\alpha'$ is normal in $Z$), so we get $\psi_V(zz'z^{-1}z'^{-1})=1$. That shows that $z$ normalizes $\psi_V$ and belongs to $Z_{\psi_V}$, hence (i). 

Let us introduce the isotropic part $\tilde{\bm}_\alpha=\bm_\alpha^{\is}$ of the simply connected covering of the derived group of $\bm_\alpha$, its minimal Levi subgroup $\tilde{\bz}_\alpha$ lifting $\bz$, and the maximal split torus $\tilde{\gs}_\alpha$ of $\tilde{\bz}_\alpha$. Write $j$ for the canonical map $\tilde{\bm}_\alpha\rg \bm_\alpha$. We have $M_\alpha'=j(\tilde{M}_\alpha)$ and $j^{-1}(Z)=\tilde{Z}_\alpha$, so $Z\cap M_\alpha'=j(\tilde{Z}_\alpha)$.

Let $v_Z:Z \rg X_*(\gs) \otimes \Q$ be the homomorphism such that $\chi(v_Z(z))= \val_F(\chi(z))$ for all $z\in S$ and all $F$-characters $\chi$ of $\mathbf S$, where $\val_F$ is the valuation of $F$ with image $\mathbb Z$; its kernel is the maximal compact subgroup of $Z$. Let $w_G : G \to X^*(Z(\widehat\bg)^{I_F})^{\sigma_F}$ be the Kottwitz  homomorphism of $G$ \cite[\S7.7]{Kot}, where $\widehat\bg$ denotes the dual group, $I_F$ the inertia subgroup of $\Gal(F^{\sep}/F)$, and $\sigma_F$ a Frobenius element of $\Gal(F^{\sep}/F)$. The parahoric subgroup $Z^0$ of $Z$ is equal to $\Ker w_Z$. By \cite[6.2]{HV1},   $\Ker w_Z=\Ker v_Z\cap \Ker w_G$.

We have the analogous map $v_{\tilde{Z}_\alpha}$ and a commutative diagram
\begin{align*}
 \begin{CD}
  \tilde{Z}_\alpha @>v_{\tilde{Z}_\alpha}>> X_*(\tilde{\gs}_\alpha) \otimes\Q \\
  @VVV  @VVV \\
  Z @>v_Z>> X_*({\gs}) \otimes\Q \\
 \end{CD}
\end{align*}
where the vertical maps are induced by $j$.

As $\tilde{\bm}_\alpha$ is semisimple and simply connected, $w_{\tilde{M}_\alpha}$ is trivial and by functoriality of the Kottwitz   homomorphism $w_G$ is trivial on $M_\alpha'=j(\tilde{M}_\alpha)$; in particular $Z^0\cap M_\alpha'=\Ker v_Z \cap M_\alpha'$. The vertical map on the right of the above diagram is injective so $j^{-1}(Z^0\cap M_\alpha')=\Ker v_{\tilde{Z}_\alpha}$. Thus $(Z\cap M_\alpha')/(Z^0 \cap M_\alpha')$ is isomorphic to $\tilde{Z}_\alpha/\Ker v_{\tilde{Z}_\alpha}$, i.e.\ to the image of $v_{\tilde{Z}_\alpha}$. Since $\tilde{\gs}_\alpha$ has dimension 1, that image is isomorphic to $\Z$. Now for $z\in \tilde{Z}_\alpha$ we have $|\alpha|(j(z))=q^{-\langle\alpha,v_Z(j(z))\rangle}=
q^{-\langle\alpha,v_{\tilde{Z}_\alpha}(z)\rangle}$ and (ii), (iii) follow. $\square$

\medskip

\noindent\textbf{Remark} From the above proof it is clear that  $v_Z(a_\alpha)$ is a (negative) rational multiple of $\alpha^\vee$. See also \ref{IV.11}  Example 3.

\subsection{}\label{III.17} Let us  derive consequences of \ref{III.16}.

\vskip2mm

\noindent \textbf{Proposition} \textit{Assume Hypothesis (H) {\upshape (\ref{III.15})}. Let $\alpha\in \Delta$ be orthogonal to $\Delta_M$. Then the following conditions are equivalent}: 

\textit{(i) $\sigma$ is trivial on $Z\cap M_\alpha'$,}

\textit{(ii) $\psi_V$ is trivial on $Z^0\cap M_\alpha'$ and $(\varphi\tau_\alpha)(v)=\varphi(v)$ for $v\in V^{U_{\op}\cap M^0}$.}

\vskip2mm

\noindent\textbf{Proof} Apply first \ref{III.12} Lemma to get 
$$(\varphi\tau_\alpha)(v)=a_\alpha^{-1}\varphi(v) \leqno(*)
$$
 for $v\in V^{U_{\op}\cap M^0}$.  Now assume (i). By \ref{III.15} Lemma,   $\psi_V$ is trivial on $Z^0\cap M_\alpha'$; then, since $\alpha$ is orthogonal to $\Delta_M$, $a_\alpha$
belongs to $Z_{\Delta_M}^\bot$  and $(*)$ implies   (ii).
Conversely assume (ii). Applying \ref{III.16} Proposition and   $(*)$ again we get that $Z\cap M_\alpha'$ acts trivially on the line  $\varphi(V^{{U}_{\op}\cap M^0})$. But as $\alpha$ is orthogonal to $\Delta_M$, $M$ normalizes $M_\alpha'$   and  hence also $Z\cap M_\alpha'$; consequently, the set of fixed points of $Z\cap M_\alpha'$ in $W$ is invariant under $M$. As it contains $\varphi(V^{{U}_{\op}\cap M^0})$ it contains $\varphi(V)$ since $V^{{U}_{\op}\cap M^0}$ generates $V$ over $M^0$,  and by hypothesis (H), $Z\cap M_\alpha'$ acts trivially on~$W$. $\square$

\vskip2mm

\noindent \textbf{Corollary} \textit{Assume Hypothesis (H) and that moreover $\varphi$ is a $\cz_M(V)$-eigenvector with supersingular eigenvalue $\chi$.
  Then $\Delta(\sigma)$ as in {\upshape (\ref{II.7})}
is the set of $\alpha\in \Delta$, orthogonal to $\Delta_M$, such that $\psi_V$ is trivial on $Z^0\cap M_\alpha'$ and $\chi(\tau_\alpha)=1$.}

\vskip2mm

\noindent\textbf{Proof}  Assume $\alpha\in \Delta(\sigma)$ is not orthogonal to $\Delta_M$. By \ref{II.7} Corollary 2 and Remark~5, there is a proper parabolic subgroup $Q=M_QN_Q$ of $M$ (containing $M\cap B$) such that $\sigma$ is trivial on $N_Q$ and is a subrepresentation of $\Ind_Q^M(\sigma_{|M_Q})$. By~\ref{III.12} Proposition~(iii), no eigenvalue of $\sigma$ can be supersingular. Consequently, any  $\alpha$ in 
$\Delta(\sigma)$ is orthogonal to $\Delta_M$ and the result follows from the proposition. $\square$

\vskip2mm

In particular, we have determined $P(\sigma)$ for a supersingular representation $\sigma$ of $M$.

\paragraph{{\large E) Weights and eigenvalues of $I(P,\sigma,Q)$}}
\addcontentsline{toc}{section}{\ \ \ E) Weights and eigenvalues of \texorpdfstring{$I(P,\sigma,Q)$}{I(P,sigma,Q)}}

\subsection{}\label{III.18} In this section, for a supersingular triple $(P,\sigma,Q)$ (\ref{I.5}), we determine the weights and eigenvalues of $I(P,\sigma,Q)$. A slightly more general situation is useful in part G) though.

\vskip2mm

\noindent \textbf{Proposition} \textit{Consider a $B$-triple $(P,\sigma,Q)$ as in {\rm \ref{I.5}} with $P=MN$, and assume that $\Delta(\sigma) $ is orthogonal to $\Delta_M$. 
Let $V$ be an irreducible \repr\ of $K$, with parameter~$(\psi_V,\Delta(V))$}.

$1)$ \textit{The following conditions are equivalent:}

\textit{(i)} $V$ \textit{is a weight of} $I(P,\sigma,Q)$,

\textit{(ii)} $V_{N^0}$ \textit{is a weight of} $\sigma$ \textit{and} $ \Delta(V)\cap \Delta(\sigma)=\Delta_Q\cap \Delta(\sigma)$.

 $2)$ \textit{If $V$ is a weight of $I(P,\sigma,Q)$, then the eigenvalues of $\cz_G(V)$ in $I(P,\sigma,Q)$ are in bijection with those of $\cz_M(V_{N^0})$ in $\sigma$ via~$\cs_M^G$.}

\vskip2mm
The proof of 1)  is in \ref{III.19}--\ref{III.21} below, that of 2) in \ref{III.22}, which actually gives more precise information.

\vskip2mm

\noindent \textbf{Remark 1}  Consider the case where $P=B$ and $\sigma$ is the  trivial representation of $B$. Then $P(\sigma)=G$ and $I(B,\sigma,Q)=\St_Q^G$. From \cite[\S 8]{Ly1} we get that $\St_Q^G$ has a unique weight $V_Q^G$, with multiplicity one, and parameter $(1,\Delta_Q)$. That weight also occurs with multiplicity one in $\Ind_Q^G1$ and the natural map $\Hom_K(V_Q^G,\Ind_Q^G1) \rg \Hom_K(V_Q^G,\St_Q^G)$ is an isomorphism; similarly $V_Q^G$ occurs with multiplicity one in $\Ind_B^G1$ and the natural map $\Hom_K(V_Q^G,\Ind_Q^G1) \rg \Hom_K(V_Q^G,\Ind_B^G1)$ is an isomorphism. Those isomorphisms are $\Hh_G(V_Q^G)$-equivariant, and  the algebra $\Hh_G(V_Q^G)$, isomorphic to the monoid algebra $C[Z^+/Z^0]$, acts via the augmentation character sending $\tau_z$ to 1 for $z\in Z^+$. That special case will be used in the proof of part 2) of the proposition.

\vskip2mm

The proposition may be applied to a supersingular triple, by \ref{III.17} Corollary.

\vskip2mm

\noindent \textbf{Corollary} \textit{
Assume $(P,\sigma,Q)$ is a supersingular  triple;
if $V$ is a weight of $I(P,\sigma,Q)$ then for any eigenvalue $\chi$ of $\cz_G(V)$ in $I(P,\sigma,Q)$, we have $\Delta_0(\chi)=\Delta_M$.}

\vskip2mm
\noindent \textbf{Proof}  By  part 2) of the proposition, $\chi$ extends to a character of $\cz_M(V_{N^0})$ so $\Delta_0(\chi)\subset \Delta_M$. On the other hand the extended character is an eigenvalue of $\sigma$ which is supersingular so $\Delta_M\subset \Delta_0(\chi)$. $\square$

\vspace{2mm}

\noindent \textbf{Remark 2} In the context of the corollary, if $P\not= G$, then no eigenvalue of $I(P,\sigma,Q)$ is supersingular.

\subsection{}\label{III.19}

By \ref{III.12} Proposition, we immediately reduce the proof of part 1) of  the proposition to the case where $P(\sigma)=G$. In the course  of the proof we shall glean more information on the weights and eigenvalues.

  We put  $\Delta_1=\Delta_M$ and $\Delta_2=\Delta(\sigma)$, so that  $\Delta$ is the  union of two orthogonal subsets $\Delta_1 $ and $ \Delta_2$.  As in  \ref{II.4} we introduce the group $\tilde{\bg}=\bg^{\is}$. It appears as the product of two factors $\tilde{\bg}_1$ and $\tilde{\bg}_2$ attached to $\Delta_1$, $\Delta_2$. Note that $\tilde G$ and $G$ have the same semisimple building and their actions on it are compatible. Let $\tilde{K}$ be the parahoric subgroup of $\tilde{G}$ attached to the point $\mathbf{x}_0$. It decomposes as $\tilde{K}_1\times \tilde{K}_2$ where for $i=1,2$, $\tilde{K}_i=\tilde{K}\cap \tilde{G}_i$ is a parahoric subgroup of $\tilde{G}_i$. Write $\iota$ for the natural map $\tilde{G}\rg G$. For $i=1,2$, let $M_i$ be the Levi subgroup $M_{\Delta_i} $ of $G$. Then  $M_i'=\iota(\tilde G_i)$ and $M_i=ZM_i'$. By \ref{II.7} Remark~4, $M_1'$ and $M_2'$ commute with each other, $Z$ normalizes each of them and~$G=ZM_1'M_2'$.

\vskip2mm

\noindent \textbf{Proposition} \textit{(i) $\tilde{K}= \iota^{-1}(K)$, $\tilde{Z}^0=\iota^{-1}(Z^0)$ and $\iota(\tilde{K_i})= K\cap M_i'$ for $i=1,2$.}

\textit{(ii) Let $\alpha\in \Phi$; then $\iota$ induces a group isomorphism of} $\tilde{U}_\alpha^0=\tilde{U}_\alpha \cap \tilde{K}$ onto $U_\alpha^0=U_\alpha\cap K$.

\vskip2mm

Here, $\tilde{U}_\alpha$ denotes the root subgroup of $\tilde{G}$ attached to~$\alpha\in \Phi$.

\vskip2mm

\noindent\textbf{Proof} By functoriality of the Kottwitz   homomorphism,  since $\tilde G$ is semisimple simply connected, $ w_G\circ \iota$ is trivial; on the other hand an element $x\in \tilde{G}$ fixes the point $\bx_0$ if and only if $\iota(x)$ fixes~$\bx_0$. 
 So we have $\tilde{K}=\iota^{-1}(K)$ and intersecting with $\tilde{Z}=\iota^{-1}(Z)$ we get $\tilde{Z}^0=\iota^{-1}(Z^0)$. If $x\in \tilde{K}_i$ then $\iota(x)\in K \cap \iota(\tilde{G}_i)=K\cap M_i'$. Conversely if $x\in \tilde{G}_i$ and $\iota(x)\in K$ then $x\in \tilde{K}\cap \tilde{G}_i=\tilde{K}_i$. This proves~(i).

(ii) Let $\alpha\in \Phi$.  As $\iota(\tilde{K})\subset K$ we have $\iota(\tilde{U}_\alpha^0)\subset U_\alpha^0$. Conversely for $x\in \tilde{U}_\alpha$, $\iota(x)\in U_\alpha^0$ implies $x\in \tilde{U}_\alpha \cap \iota^{-1}(K)= \tilde{U}_\alpha^0$ by~(i). $\square$

\vskip2mm

\noindent \textbf{Corollary} \textit{We have $K=Z^0 \iota(\tilde{K})$. For $i=1,2$, $M_i^0=Z^0 \iota(\tilde{K}_i)$.
}

\vskip2mm

\noindent\textbf{Proof} This comes from (ii) of the proposition, given \ref{III.7} Lemma. $\square$

\medskip

\noindent \textbf{Remark} By \ref{II.7} Remark 4, $M_1'\cap M_2'$ is finite and central in $G$. As it is contained in $\Ker w_G$, it follows that $Z^0$ contains $M_1'\cap M_2'$, which is equal to $\iota(\tilde{K}_1)\cap \iota(\tilde{K}_2)$.

\subsection{}\label{III.20} Let now $(\rho,V)$ be an irreducible \repr\ of $K$.  We want to write $V$ as a tensor product adapted to the orthogonal decomposition $\Delta=\Delta_1\sqcup \Delta_2$.

Write $(\tilde{\rho},\tilde{V})$ for the representation of $\tilde{K}$ obtained from $\rho$ via $\iota:\tilde{K}\rg K$. By \ref{III.19} Proposition (ii) $\overline {\iota(\tilde K)}$ contains $\overline G'$, so by \ref{III.11} Corollary $\tilde \rho$  is irreducible. Since $\tilde{K}=\tilde{K}_1\times \tilde{K}_2$, $\tilde{V}$ decomposes as a tensor product $\tilde{V}_1 \otimes \tilde{V}_2$ where  for $i=1,2$, $\tilde{V}_i$ is an irreducible \repr\ of $\tilde{K}_i$ which is trivial on $\tilde K_{3-i}$.

 To decompose $V$ as a tensor product $V_1\otimes V_2$ of irreducible representations of $K$, where $V_1$ restricts to $\tilde{V}_1$ via $\iota$, and $V_2$ to $\tilde{V}_2$, we have to take some care, as $K$ is not the direct product $M_1^0\times M_2^0$. 

\medskip

\noindent \textbf{Proposition} \textit{(i) For $i=1,2$, let $V_i$ be an irreducible representation of $K$  trivial on $K\cap M'_{3-i}$. Then $V_1\otimes V_2$ is irreducible with parameter $(\psi_{V_1}\psi_{V_2},\Delta(V_1) \cap \Delta(V_2))$. Moreover, $\Delta(V_i) $ contains $\Delta_{3-i}$.}

\textit{(ii) Let  $V$ be an  irreducible \repr\ of $K$.  If $V_2$ is  an irreducible  representation of $K$  trivial on $K\cap M'_{1}$ with $\Hom_{K\cap M'_2}(V_2,V)\neq 0$, then $V_1=\Hom_{K\cap M'_2}(V_2,V)$ is an irreducible  representation of $K$  trivial on $K\cap M'_{2}$  and $V\simeq V_1\otimes V_2$.}

\textit{(iii) Let  $V$ be an  irreducible \repr\ of $K$.  Then $V\simeq V_1\otimes V_2$ with $V_i$ as in (i) if and only if $V$ is trivial on $M'_1\cap M'_2$.}

\medskip

We will not need part (iii), we only included it for completeness.

\medskip

\noindent \textbf{Proof} (i)   Let $\tilde V_i$ be the pullback of $V_i$ to $\tilde K$ via $\iota$. Then $\tilde V_i$ is trivial on $\tilde K_{3-i}$, so $\tilde V_1 \otimes \tilde V_2$ is an irreducible  representation of $\tilde K$. Hence $V:=V_1\otimes V_2$ is an irreducible  representation of $ K$. If $Q=M_QN_Q$ is a parabolic subgroup containing $B$, then
$$V_{N_Q^0}\simeq (V_1)_{N_Q^0}\otimes (V_2)_{N_Q^0}, \ \text{as} \ N_Q^0=(N_Q^0\cap M_1')\times (N_Q^0\cap M_2').$$
Hence  by \ref{III.10},  $\Delta_Q \subset \Delta(V)$ if and only if $\Delta_Q \subset \Delta(V_i)$ for $i=1,2$, so $ \Delta(V)=\Delta(V_1) \cap \Delta(V_2)$. Taking $Q=B$, we deduce $\psi_V=\psi_{V_1}\psi_{V_2}$. As  $K\cap M'_{3-i}$ is trivial on $V_i$, we get $\Delta_{3-i}\subset \Delta(V_i)$.

(ii) This follows from Clifford theory \cite[Lemma 5.3]{Abe}.

(iii) The ``if'' direction is obvious. Assume that $V$ is trivial on $M'_1\cap M'_2$.  Let $W$ be an irreducible representation of $K\cap M'_{2}$ such that $\Hom_{K\cap M'_{2}}(W,V)\neq 0$. Via $\iota$, $W$ is an irreducible representation of $\tilde K_2$, which we consider as a representation $\tilde W$ of $\tilde K$ trivial on $\tilde K_1$. As $V$, hence $W$, is trivial $\iota(\tilde K_1)\cap \iota(\tilde K_2)$ by assumption, it follows that $\tilde W$ is trivial on $\Ker \iota$, so we have extended $W$ to an irreducible representation of $K\cap G'$, which is trivial on $K\cap M'_1$.  We may view $W$ as an irreducible representation of $\overline {G'}$ and we choose an irreducible representation $V_2$ of  $\overline {G}$ such that $W$ occurs in $V_2|_{\overline {G'}}$. By \ref{III.11} Corollary $W \simeq V_2|_{\overline {G'}}$ and hence $\Hom_{K\cap M'_2}(V_2,V)\neq 0$. By part (ii), $V\simeq  V_1\otimes V_2$ with $V_i$ as in (i). 
$\square$

\subsection{}\label{III.21} Let $(P,\sigma,Q)$ be a  $B$-triple with $P(\sigma)=G$.
We are now finally ready to determine the weights of $^e\sigma \otimes \St_Q^G$. We keep the notation of \ref{III.19}. Recall that 
 by construction $^e\sigma$ is trivial on $M_2'$ and $\St_Q^G$ is trivial on $M_1'$.
 
 Let us fix  a weight $V$ of $I(P,\sigma,Q)={}^e\sigma \otimes \St_Q^G$. We decompose the pullback $\tilde V$ of $V$ to   a representation of   $\tilde{K}=\tilde{K}_1\times \tilde{K}_2$, via $\iota$, as $\tilde V \simeq \tilde V_1\otimes \tilde V_2$. Therefore $\Hom_K(V, {}^e\sigma \otimes \St_Q^G)$ injects into
 $$
 \Hom_{\tilde K}(\tilde V, {}^e\sigma \otimes \St_{\tilde Q}^{\tilde G})\simeq  \Hom_{\tilde K}(\tilde V_1, {}^e\sigma )\otimes  \Hom_{\tilde K}(\tilde V_2,  \St_{\tilde Q}^{\tilde G}),
 $$
where we used that $\tilde{K}_1$ acts trivially on $\tilde V_2,  \St_{\tilde Q}^{\tilde G} $ and $\tilde{K}_1$ acts trivially on $\tilde V_1, {}^e\sigma$.  As $ \St_{\tilde Q}^{\tilde G} $ has a unique weight (\ref{III.18}), $\tilde V_2$ is the pullback via $\iota$ of the unique weight  $V_2$ of $ \St_{Q}^{G} $. By lifting via $\iota:\tilde K_2\to \iota(\tilde K_2)= K\cap M'_2$, we deduce $\Hom_{K\cap M'_2}(V_2,V)=\Hom_{\tilde K_2}(\tilde V_2, \tilde V) \neq 0$. By \ref{III.20} Proposition (ii), $V\simeq V_1\otimes V_2$ for some irreducible representation $V_1$ of $K$ trivial on $K\cap M'_2$. We also see by \ref{III.20} Proposition (i) and \ref{III.18}  Remark 1 that $\Delta(V)\cap \Delta_2=\Delta(V_2)\cap \Delta_2=\Delta_Q\cap \Delta_2$.  The natural injection 
$\Hom_{ K}( V_2, \St_Q^G) \hookrightarrow \Hom_{ K\cap M'_2}( V_2, \St_Q^G)$ is an isomorphism of $1$-dimensional vector spaces, because the right-hand side is isomorphic to $ \Hom_{\tilde K}(\tilde V_2,  \St_{\tilde Q}^{\tilde G})$ via $\iota$.
 Thus the following lemma, in our situation, implies that  $V_1$ is a weight of $^e\sigma$, so $V_{N^0}$ is a weight of $\sigma$. This proves that (i) implies (ii) in \ref{III.18} Proposition 1).

\medskip

\noindent\textbf{Lemma} \textit{Let $\sigma_1$ be a representation of $G$ trivial on $M_2'$, $\sigma_2$ a representation of $G$ trivial on $M_1'$. Let $V_1$ be an irreducible representation of $K$ trivial on $ K\cap M'_2$, $V_2$ an irreducible representation of $K$ trivial on $ K\cap M'_1$. Assume that the inclusion $\Hom_{ K}(V_2,\sigma_2)\rg \Hom_{K\cap M'_2}(V_2,\sigma_2)$ is an isomorphism. Then the natural inclusion of $\Hom_K(V_1,\sigma_1) \otimes \Hom_K(V_2,\sigma_2)$ into $\Hom_K(V_1\otimes V_2,\sigma_1\otimes\sigma_2)$ is an isomorphism.}

\medskip

\noindent\textbf{Proof} Look first at points fixed by $ K\cap M'_2$ in $\Hom(V_1\otimes V_2,\sigma_1\otimes \sigma_2)$. As $ K\cap M'_2$ acts trivially in $V_1$ and $\sigma_1$, it is simply $\Hom(V_1,\sigma_1) \otimes \Hom_{ K\cap M'_2}(V_2,\sigma_2)$, so by the assumption it is also $\Hom(V_1,\sigma_1) \otimes \Hom_K(V_2,\sigma_2)$. Now $K$ acts trivially on $\Hom_K(V_2,\sigma_2)$, so taking fixed points under $K$ indeed gives $\Hom_K(V_1,\sigma_1) \otimes \Hom_K(V_2,\sigma_2)$. $\square$

\medskip

We now prove that (ii) implies (i) in \ref{III.18} Proposition 1). Let $V$ be an irreducible representation of $K$ satisfying (ii). From \ref{III.12} Proposition (i), $V$ is a weight of $\Ind_P^G\sigma \simeq {}^e \sigma \otimes \Ind_P^G 1$. Therefore, $V$ is a weight of $I(P,\sigma, Q')$ for some parabolic $Q'\supset P$. As we have already proved that (i) implies (ii) in \ref{III.18} Proposition 1), we deduce that $\Delta_{Q'}\cap \Delta_2=\Delta_Q\cap \Delta_2 $, so $Q'=Q$.
  $\square$

\subsection{}\label{III.22} It remains to prove part 2) of \ref{III.18} Proposition. We in fact establish something more precise, which gives what we need by \ref{III.12} Proposition. Also, by that proposition we may assume $P(\sigma)=G$.

\medskip

\noindent\textbf{Lemma 1} \textit{Let $(\rho,V)$ be a weight of $I(P,\sigma,Q)$ where $P(\sigma)=G$.}

\textit{(i) The quotient map $\Ind_Q^{G}1\rg \St_Q^{G}$ induces an $\Hh_G(V)$-isomorphism}
$$
\Hom_K(V,\Ind_Q^G {}^e\sigma) \lgr \Hom_K(V,I(P,\sigma,Q)).
$$

\textit{(ii) The inclusion $\Ind_Q^{G}1 \rg \Ind_P^{G}1$ induces an $\Hh_G(V)$-isomorphism}
$$
\Hom_K(V,\Ind_Q^G {}^e\sigma) \lgr \Hom_K(V,\Ind_P^G\sigma).
$$

\medskip

\noindent\textbf{Proof} It is clear that the maps in (i), (ii) are $\Hh_G(V)$-equivariant.  As in \ref{III.21} write $V$ as $V_1\otimes V_2$ where $V_2$ is the   unique weight of  $\St_Q^G$ (it has parameter $(1, \Delta_Q)$). By \ref{III.21} Lemma (the hypothesis is verified by pulling back via $\iota$, as in \ref{III.21}), we get isomorphisms
\begin{equation*}
\begin{matrix}
\Hom_K(V_1\otimes V_2, {}^e\sigma \otimes \St_Q^G) &\simeq &\Hom_{K}(V_1,{}^e\sigma) \otimes \Hom_K(V_2,\St_Q^G),\\
\Hom_K(V_1\otimes V_2, {}^e\sigma \otimes \Ind_Q^G1) &\simeq &\Hom_{K}(V_1,{}^e\sigma) \otimes \Hom_K(V_2,\Ind_Q^G1),\\
\Hom_K(V_1\otimes V_2, {}^e\sigma \otimes \Ind_P^G1) &\simeq & \Hom_{K}(V_1,{}^e\sigma) \otimes \Hom_K(V_2, \Ind_P^G1).\\
\end{matrix}
\end{equation*}

The maps $\Ind_Q^G1\rg \St_Q^G$ and $\Ind_Q^G1\rg \Ind_P^G1$ induce on each side vertical maps which give commutative diagrams. As the vertical maps on the right-hand side are isomorphisms by \ref{III.18} Remark 1, so are the vertical maps on the left-hand side, and (i), (ii) are implied by the following well-known lemma. $\square$

\vskip2mm

\noindent \textbf{Lemma 2} \textit{Let  $H'$ be a closed subgroup of a locally profinite group $H$ and $\ind_{H'}^H$  the smooth compact induction functor. Let  $V$ be a smooth \repr\ of  $H'$ and $W$  a smooth \repr\ of   $H$. Then there is an isomorphism $\Phi$ of \repr s of $H$, $W\otimes \ind_{H'}^{H} V \stackrel{\sim}{\longrightarrow} \ind_{H'}^H(W\otimes V)$, given by the formula
}
$$
\Phi(w\otimes f): h \longmapsto hw\otimes f(h) \qquad \mathrm{for}\ w\in W,\ f\in \ind_{H'}^HV.
$$

\vskip2mm

\paragraph{{\large F) Irreducibility of $I(P,\sigma,Q)$}}
\addcontentsline{toc}{section}{\ \ \ F) Irreducibility of \texorpdfstring{$I(P,\sigma,Q)$}{I(P,sigma,Q)}}

\subsection{}\label{III.23} 
 \textbf{Proposition} \textit{Let $(P,\sigma,Q)$ be a supersingular triple. Then $I(P,\sigma,Q)$ is irreducible.  }
\vskip2mm
 
\noindent \textbf{Proof} 
 It is enough to prove that if $V$ is an  irreducible \repr\ of $K$ and $\varphi\in \Hom_K(V,I(P,\sigma,Q))$ is  a $\cz_G(V)$-eigenvector with eigenvalue $\chi$, then the sub\repr\ $X$ of $I(P,\sigma,Q)$ generated by $\varphi(V)$ is $I(P,\sigma,Q)$. So we fix such a situation and write $(\psi_V,\Delta(V))$ for the parameter of $V$. We prove the result by induction on the cardinality of~$\Delta(V)$.

By \ref{III.14} Corollary 1  we have $X= I(P,\sigma,Q)$ if $\Delta(V)\subset \Delta_{P(\sigma)}$, so we assume that  this is not the case. We pick $\alpha$ in $\Delta(V)$ but not in $\Delta_{P(\sigma)}$, and let $V'$ be an irreducible \repr\ of $K$ with parameters $(\psi_V,\Delta(V)-\{\alpha\})$. Note that $V_{U^0}'$ and $V_{U^0}$ are isomorphic, so that $\chi$  defines a character of $\cz_G(V')$ via the Satake isomorphism, which we also denote by $\chi$.

Via $\varphi$, $X$ is a quotient of $\chi\otimes_{\cz_G(V)}\ind_K^GV$. By \ref{III.18} Corollary  $\Delta_0(\chi)=\Delta_M$, hence $\alpha \notin \Delta_0(\chi) $.
By the change of weight theorem (\ref{IV.2} Corollary),
 $\chi\otimes_{\cz_G(V)}\ind_K^GV$ and $\chi\otimes_{\cz_G(V')}\ind_K^GV'$ are isomorphic unless $\alpha$ is orthogonal to $\Delta_0(\chi)$, $\psi_V$ is trivial on $Z^0 \cap M_\alpha'$ and $\chi(\tau_\alpha)=1$ 
 (see \ref{III.16} for the notation $\tau_\alpha$). By induction then, we are reduced to the case where $\alpha$ is orthogonal to $\Delta_0(\chi)$, $\psi_V$ is trivial on $Z^0\cap M_\alpha'$ and $\chi(\tau_\alpha)=1$. As $\Delta_0(\chi)=\Delta_M$,  the conditions imply  (\ref{III.17} Corollary) that $\alpha$ belongs to $\Delta(\sigma)\subset \Delta _{P(\sigma)}$ contrary to assumption. $\square$

\paragraph{{\large G) Injectivity of the parametrization}}
\addcontentsline{toc}{section}{\ \ \ G) Injectivity of the parametrization}

\subsection{}\label{III.24} Let $(P_1,\sigma_1,Q_1)$ and $(P_2,\sigma_2,Q_2)$ be supersingular triples such that
$$I(P_1,\sigma_1,Q_1) \simeq I(P_2,\sigma_2,Q_2).$$
Let $V$ be a weight of $I(P_1,\sigma_1,Q_1)$, with parameter $(\psi_V,\Delta(V))$, and $\chi$ an eigenvalue of $\cz_G(V)$ in $I(P_1,\sigma_1,Q_1)$. We have seen $\Delta_0(\chi)=\Delta_{P_1}$ (\ref{III.19} Corollary) so we deduce $\Delta_{P_1}=\Delta_{P_2}$ and $P_1=P_2$. Write $P_i=M_iN_i$ as usual. By \ref{III.18} Proposition, $V_{N_i^0}$ is a weight of 
$\sigma_i$ with supersingular eigenvalue $\chi$ (via $\cs_{M_i}^G)$. Then \ref{III.17} Corollary implies that  $P(\sigma_1)=P(\sigma_2)$. Taking the ordinary part functor \cite{Eme,Vig3} with respect to $P(\sigma_1)$, we deduce that $^e{\sigma_1} \otimes \St_{Q_1}^{P(\sigma_1)}$ and $^e{\sigma_2} \otimes \St_{Q_2}^{P(\sigma_2)}$ are isomorphic as \repr s of $P(\sigma_1) = P(\sigma_2)$. From \ref{II.8} Remark, we get $Q_1=Q_2$  and $\sigma_1 \simeq \sigma_2$. This completes the proof of the uniqueness in  \ref{I.5} Theorem 4.

\medskip

We insert here a consequence of the irreducibility of $I(P,\sigma,Q)$ and of the injectivity of the parametrization, which we shall use in part H) and generalize in Chapter~\ref{VI}.

\medskip

\noindent \textbf{Proposition} \textit{Let $P=MN$ be a parabolic subgroup of $G$ containing $B$, and $\sigma$ a supersingular \repr\ of $M$, inflated to $P$.   Then the irreducible components of $\Ind_P^G\sigma$ are the $I(P,\sigma,Q)$, $Q$ a parabolic subgroup of $G$ with $P\subset Q\subset P(\sigma)$; each occurs with multiplicity $1$. In particular $\Ind_P^G\sigma$ has finite length.}

\vskip2mm

\noindent\textbf{Proof}  The representation  $\Ind_P^{P(\sigma)}\sigma$ is isomorphic to $^e\sigma\otimes \Ind_P^{P(\sigma)}1$ (\ref{III.22} Lemma 2), which has a filtration with subquotients $^e\sigma\otimes \St_Q^{P(\sigma)}$, one for each parabolic subgroup $Q$ with $P\subset Q\subset P(\sigma)$. The proposition then follows from 
\ref{III.23} Proposition  by parabolic induction from~$P(\sigma)$~to~$G$.  $\square$

\paragraph{{\large H) Surjectivity of the parametrization}}
\addcontentsline{toc}{section}{\ \ \ H) Surjectivity of the parametrization}

\subsection{}\label{III.25} Let $(\pi,W)$ be an irreducible admissible \repr\ of $G$. To prove that $\pi$ has the form $I(P,\sigma,Q)$ for a supersingular triple $(P,\sigma,Q)$, we use induction on the semisimple rank~of~$G$.

If $\Delta_0(\chi)=\Delta$ for all weights $V$ of $\pi$ and corresponding eigenvalues $\chi$, then $\pi$ is supersingular and $\pi\simeq I(G,\pi,G)$. So we fix a weight $V$ for $\pi$   with $\cz_G(V)$-eigenvalue $\chi$ such that $\Delta_0(\chi)\not= \Delta$. By construction $\pi$ is a quotient of $\chi\otimes_{\cz_G(V)}\ind_K^G V$.

Let $P = MN$ be the parabolic subgroup such that $\Delta_P = \Delta_0(\chi)$.
Consider $\sigma = \chi\otimes\ind_{M^0}^MV_{N^0}$.
By the filtration theorem (\ref{I.6} Theorem 6, proved in Chapter~\ref{V}), $\chi\otimes\ind_K^GV$ has a filtration with subquotients $I_e(P,\sigma,Q)=\Ind_{P_e}^G({}^e\sigma\otimes\St_Q^{P_e})$ where $P\subset Q\subset P_e$.
So $\pi$ is a quotient of some $I_e(P,\sigma,Q) $.
If $P_e \ne G$, then by \cite[Proposition~7.9]{HV2} (note that $\sigma$ has a central character by \ref{III.12} Lemma) there is an irreducible admissible \repr\ $\rho$ of the Levi quotient of $P_e$ such that $\pi$ is a quotient of $\Ind_{P_e}^G \rho$. By the induction hypothesis  and \ref{III.24} Proposition, $\rho$ is an irreducible constituent of $\Ind_{P_1}^{P_e}\rho_1$ where $P_1 $ is a  parabolic subgroup of $P_e$ containing $B$, and $\rho_1$ is a supersingular \repr\ of the Levi quotient of $P_1$. Then $\pi$ is an irreducible constituent of $\Ind_{P_1}^G\rho_1$, so by \ref{III.24} Proposition it is isomorphic to $I(P_1,\rho_1,Q')$ for some $Q'$.

If $P_e = G$, $\pi$ is a quotient of some $^e\sigma \otimes \St_Q^G$.
By \ref{II.8} Proposition and Remark, $\pi $ is isomorphic to $ ^e\sigma_\pi \otimes \St_Q^G$ for some irreducible admissible representation $\sigma_\pi$ of $M$.
The eigenvalues of $\sigma_\pi$ are those of $\pi$ by \ref{III.18} Proposition, and since $\Delta_M=\Delta_0(\chi)$, $\sigma_\pi$ has a supersingular eigenvalue. As $\Delta_M\neq \Delta$, the induction hypothesis implies that $\sigma_\pi$ is supersingular, cf.~\ref{III.18} Remark 2, and $\pi\simeq I(P,\sigma_\pi,Q)$. $\square$

\subsection{}\label{III.26} It is worth commenting on the admissibility assumptions in our results. The reader may notice that, since admissibility plays no r\^ole in Chapters \ref{IV} and \ref{V}, our results would still be true if instead of irreducible admissible representations, we considered irreducible representations $(\sigma,W)$ such that for some weight $(\rho,V)$ of $\sigma$, $\Hom_K(V,W)$ contains an eigenvector for $\mathcal{Z}_G(V)$. But the classification thus obtained would depend on the choice of $K$, $\bs$, $\bb$, whereas we shall see in Chapter \ref{VI} that with the admissibility assumption it does not depend on those choices. Of course one may hope that the condition above actually implies admissibility or even, as is the case for  complex representations, that any irreducible representation  of $G$ is admissible. 
Note that because of our admissibility condition we do not assert that $G$ has any supersingular
representation. When $G=\GL_n(F)$ and $F$ has characteristic $0$, we will show in forthcoming work that 
supersingular representations of $G$ exist.

\section{Change of weight}\label{IV}

\subsection{}\label{IV.1} The main goal of this chapter is to establish our change of weight theorem (\ref{IV.2} Corollary)  used in \ref{III.23}. Before commenting on the method of proof, let us state precisely what we prove here. We fix an irreducible representation $\rho$ of $K$ on a space $V$, with parameter $(\psi_V,\Delta(V))$  as defined in \ref{III.9}. We consider the ``universal''  representation $\ind_K^GV$, which we see as a sub-representation of $\Ind_B^G(\ind _{Z^0}^Z(V_{U^0}))$ via the injective canonical intertwiner (III.13.1). 

We assume that $ \Delta(V)$ is non-empty, and we  choose $\alpha\in  \Delta(V)$  and let $(\rho',V')$ be the irreducible representation of $K$ with parameter $(\psi_V, \Delta(V)-\{\alpha\})$. Similarly we consider the universal representation $\ind_K^GV'$ as a subrepresentation of $\Ind_B^G(\ind_{Z^0}^Z V'_{U^0})$.

To compare the two universal representations, we fix non-zero vectors $v$ in $V$ and $v'$ in $V'$ which are invariant under $U^0_{\op}$. The image of $v$ in $V_{U^0}$ is then a basis of $V_{U^0}$, and similarly for $v'$. Using those images as basis vectors, we obtain embeddings of $\ind_K^GV$ and $\ind_K^GV'$ into the same representation $\Ind_B^G(\ind _{Z^0}^Z\psi_V)$. Moreover the Satake isomorphism induces an algebra homomorphism $\Hh_G(V) \rg \Hh_Z(\psi_V)$; the algebra $\Hh_Z(\psi_V)$ acts on $\ind _{Z^0}^Z\psi_V$, hence on $\Ind_B^G(\ind _{Z^0}^Z\psi_V)$, and the embedding $\ind_K^G(V) \rg  \Ind_B^G(\ind _{Z^0}^Z\psi_V)$ is $\Hh_G(V)$-equivariant. We have similar properties for $V'$. Note that $\Hh_G(V)$ and $\Hh_G(V')$ have the same image in $\Hh_Z(\psi_V)$, so we identify them with that common image, which we write $\Hh_G$, and similarly we write $\cz_G$ for their common centre.

For $z$ in $Z$ normalizing $\psi_V$, we have the function $\tau_z$ in $\Hh_Z(\psi_V)$ with support $Z^0z$ and value $1_C$~at~$z$. Recall from \ref{III.16} the notation $a_\alpha\in  Z\cap M_\alpha'$ and $\tau_\alpha=\tau_{a_\alpha}  \in \cz_Z(\psi_V)$, when  $\psi_V$ is trivial on $Z^0\cap M_\alpha'$.

\medskip

\noindent\textbf{Theorem} \textit{Let $z\in Z^+$. Assume that $z$ normalizes $\psi_V$ and that $|\alpha|(z)<1$. We have: }

\noindent$(i)$ $ 
\tau_z(\ind_K^GV) \subset \ind_K^GV'. 
$ 

\noindent$(ii)$ \textit{If $\psi_V$ is not trivial on $Z^0\cap M_\alpha'$, then $\tau_z(\ind_K^GV') \subset \ind_K^GV$.}

\noindent$(iii)$ \textit{If $\psi_V$ is trivial on $Z^0\cap M_\alpha'$, then $
\tau_z(1-\tau_\alpha)(\ind_K^GV') \subset \ind_K^GV$.}

\medskip

\noindent\textbf{Remark} In (iii) $\tau_z(1-\tau_\alpha)=\tau_z-\tau_{za_\alpha}$ belongs to $\cz_G(V)$ if $z\in Z_{\psi_V}$ and $za_\alpha$ belongs to $Z^+$; moreover, if $|\alpha|(z)$ is small enough, $za_\alpha$ belongs to $Z^+$.

\subsection{}\label{IV.2} 

We obtain our change of weight theorem:
\medskip

\noindent \textbf{Corollary} \textit{ Let $\chi$ be a character of $\cz_G$ and assume that $\alpha\notin \Delta_0(\chi)$. Then $\chi\otimes_{\cz_G}\ind_K^GV$ and $\chi\otimes_{\cz_G}\ind_K^GV'$ are isomorphic unless $\alpha$ is orthogonal to $\Delta_0(\chi)$, $\psi_V$ is trivial on $Z^0 \cap M_\alpha'$ and $\chi(\tau_\alpha)=1$. 
}

\vskip2mm  We remark that  $\chi(\tau_\alpha)$ is well defined if  $\alpha$ is orthogonal to $\Delta_0(\chi)$ (\ref{III.4}, \ref{III.16} Notation).

\vskip2mm 

\noindent\textbf{Proof} Choose $z$ as in the theorem, with $\chi(\tau_z)\not=0$. For example, we can take for $z$ the element $z_\alpha$ of \ref{III.4}, since $\alpha\notin \Delta_0(\chi)$: then $\chi(\tau_{z_\alpha})\not=0$. Multiplying by $\tau_z$ in $\Ind_B^G(\ind_{Z^0}^Z\psi_V)$ is $\cz_G$-linear, so, when $\psi_V$ is not trivial on $Z^0 \cap M_\alpha'$, by (i) and (ii) of the theorem, $\tau_z$ induces $G$-equivariant maps from $\ind_K^G V$ to  $\ind_K^GV'$ and back. The composites in both directions are given by the action of $\tau_z^2$. Tensoring with $\chi$, we see that the representations $\chi\otimes_{\cz_G}\ind_K^G V$ and $\chi\otimes_{\cz_G}\ind_K^G V'$ are isomorphic, because $\chi(\tau_z^2)\not=0$. That gives the desired result when $\psi_V$ is non-trivial on $Z^0\cap M_\alpha'$.

Assume then that $\psi_V$ is trivial on $Z^0\cap M_\alpha'$. Replacing $z$ by a positive power, we may assume $za_\alpha\in Z^+$. If $\alpha$ is not orthogonal to $\Delta_0(\chi)$ then there is $\beta$ in $\Delta_0(\chi)$ with $|\beta| (za_\alpha)<1$ and then $\chi(\tau_{za_\alpha})=0$, so the same reasoning applies, using (iii) instead of (ii). It similarly applies if $\alpha$ is orthogonal to $\Delta_0(\chi)$ and $\chi(\tau_\alpha)\not=1$. $\square$

\subsection{}\label{IV.3} Let us now comment on the proof of \ref{IV.1} Theorem. 
We abbreviate $\psi=\psi_V, J=\Delta(V), J'=J-\{\alpha\}$, and $ \cx= \Ind_B^G(\ind_{Z^0}^Z\psi)$. Let $I$ be the pro-$p$ Iwahori subgroup of $G$ which is the inverse image in $K$ of $U_k^{\op}$.\footnote{Beware of the notation: here, for convenience, we write $I$ for a pro-$p$ ``lower'' Iwahori subgroup.}

We first remark that $\ind_K^GV$ is generated, as a representation of $G$, by a single element, the function with support $K$ and value $v$ at $1_G$. We write $f$ for its image in $\cx$;   it is described explicitly in \ref{IV.4} below.
 Similarly we have a function $f'$ in $\cx$, corresponding to $v'$, which generates the subrepresentation~$\ind_K^GV'$.
 We use work of the fourth-named author \cite{Vig4} which determines the structure of the Hecke algebra $\Hh=\Hh(G,I)$, the intertwining algebra in $G$ of the trivial character of $I$. The space $\cx^I$ is a right module over $\Hh$, and for $x\in \cx^I$ and $T$ in $\Hh$, $xT$ belongs to the $G$-subspace generated by $x$. By construction, $f$ and $f'$ belong to $\cx^I$ and to prove the theorem we  show that: for (i) $\tau_zf\in f'\Hh$; for (ii) $\tau_zf'\in \cz_G f+f \Hh$; for (iii) $\tau_z(1-\tau_\alpha)f'\in \cz_G f+f\Hh$. That is not an easy matter and takes up the rest of this chapter.

\subsection{}\label{IV.4} Let us first identify  the  element  $f\in \cx^I$; the obvious analogue will hold for~$f'$.

As $G=BK$ it is enough to specify $f$ at $g\in K$. Going through the construction of the embedding $\ind_K^GV\rg \Ind_B^G(\ind _{Z^0}^Z\psi)$ we get that for $g$ in $K$, $f(g)$ is the function in $\ind_{Z^0}^Z\psi$ with support $Z^0$ and value $\varepsilon(g)$ at $1$, where $\overline{gv}=\varepsilon(g)\overline{v}$ in $V_{U^0}$, bars indicating the images under~$V\rg V_{U^0}$.

The value $\varepsilon(g)$ depends only on the image $\overline{g}$ of $g$ in $K/K(1)$, we write accordingly $\varepsilon(\overline{g})$. By \cite[Corollary 3.19]{HV2} we have $\varepsilon(\overline{g})\not=0$ if and only if $\overline{g}$ belongs to $B_k P_{J,k}B_k^{\op}$ (recall from \ref{III.9} Definition that $P_{J,k}$ is the stabilizer in $G_k$ of the kernel of the quotient map $V\rg V_{U^0})$; that last set is also $P_{J,k}U_k^{\op}$. We can be more precise; we obviously have $\varepsilon(\overline{g}x)=\varepsilon(\overline{g})$ for $x\in U_k^{\op}$, so it is enough to describe $\varepsilon_{|P_{J,k}}$. Since $P_{J,k}$ is the stabilizer in $G_k$ of the kernel of $V\rg V_{U^0}$,  the restriction $\varepsilon_{|P_{J,k}}$ is a character $P_{J,k}\rg C^\times$; as such it is trivial on unipotent elements. On $Z_k$ it is given by the action of $Z_k$ on $V^{U_k^{\op}}$ or $V_{U_k}$, so it is equal to $\psi$ there. In other words, on $P_{J,k}$ the character $\varepsilon$ is simply the (unique) extension of $\psi$  to $P_{J,k}$.

\subsection{}\label{IV.5} To relate $f$ and $f'$ we shall express both of them in terms of Hecke operators in the subalgebra $\Hh(K,I)$ of $\Hh(G,I)$ acting on a single function $f_0$ in~$\cx^I$.

We first describe the double coset spaces $I\ba G/I$ and $B\ba G/I$. Recall that the Weyl group $W_0$ of $G$ can be seen as $\cn^0/Z^0$ or $\cn_k/Z_k$. As $G=BK$ the inclusion of $K$ in $G$ induces a bijection $B^0\ba K/I\simeq B\ba G/I$; as moreover $I$ contains the normal subgroup $K(1)$ of $K$, reduction mod $K(1)$ induces a bijection $B^0\ba K/I  \simeq B_k\ba G_k/U_k^{\op}$ and the Bruhat decomposition in $G_k$ gives a bijection $\cn_k/Z_k \simeq B_k\ba G_k/U_k^{\op}$. All in all, we see that the map $\cn^0 \rg B\ba G/I$ $g\mapsto BgI$ induces a bijection $W_0=\cn^0/Z^0 \simeq B\ba G/I$.

On the other hand, the map $\cn\rg I\ba G/I$ induces a bijection $\cn/(Z\cap K(1) )\simeq I\ba G/I$ and, by restriction, a bijection $\cn^0/(Z\cap K(1)) \simeq I\ba K/I$. Under reduction modulo $K(1)$ we get the bijection $\cn_k \simeq U_k^{\op}\ba G_k /U_k^{\op}$ given by the Bruhat decomposition.

\medskip

\noindent\textbf{Notation} Recall that $Z(1)=Z\cap K(1)$ is  the  unique  pro-$p$ Sylow subgroup of $Z^0$ and that it is normal in $Z$. We write $_1{W}$ for the group $\cn/Z(1)$ and $_1W_0$ for the group $\cn^0/Z(1)$ (naturally isomorphic to $\cn_k$), ${W}$ for the group $\cn/Z^0$. We have obvious exact sequences of groups
\begin{align*}
1 &\lgr Z_k \lgr {}_1W_0 \lgr W_0 \lgr 1,\\
1 &\lgr Z_k \lgr {}_1{W} \lgr {W} \lgr 1.
\end{align*}
Moreover ${W}$ is the semi-direct product of $\Lambda=Z/Z^0$ with $W_0$ viewed  as $\cn^0/Z^0$. We also put $_1\Lambda=Z/Z(1)$ and $_1\Lambda^+=Z^+/Z(1)$.

For $g$ in $G$ we write $T(g)$ for the double coset $IgI$ viewed as an element of $\Hh(G,I)$. On an element $\varphi$ in $\cx^I$ it acts via
$$
(\varphi T(g))(h) = \sum_{x\in I/(I\cap g^{-1}Ig)} \varphi(hxg^{-1}) \ \mathrm{for}\ h\in G. \leqno(\mathrm{IV}.5.1)
$$
When $g\in \cn$, $T(g)$ depends only on the class $w$ of $g$ modulo $Z(1)$, and we write $T(w)$ for $T(g)$. In a similar manner, reduction modulo $K(1)$ gives an isomorphism of $\Hh(K,I)$ onto $\Hh(G_k,U_k^{\op})$; accordingly for $g\in K$, $T(g)$ depends only on the reduction $\overline{g}$ of $g$ in $G_k$ and we write also $T({\overline{g}})$. In fact we shall also have use of the Hecke algebras with integer coefficients $\Hh_\Z(G,I)$ and $\Hh_\Z(K,I)$ (isomorphic to $\Hh_\Z(G_k,U_k^{\op}))$ and we use the same notations~$T(g)$,~$T(w), T({\overline{g}})$.

\subsection{}\label{IV.6} Basic generators and relations for $\Hh_\Z(G,I)$ and  $\Hh_\Z(K,I)$ are given in \cite{Vig4}. By tensoring with $C$ they give generators and relations for $\Hh(G,I)$ and $\Hh(K,I)$. We now state  the results we use, referring to \cite{Vig4} for details. We need a bit more notation, though.

For $\beta\in \Delta$, we let $s_\beta$ be the corresponding reflection in $W_0$. We put $\Sigma_0=\{s_\beta\mid \beta\in \Delta\}$. The pro-$p$ Iwahori subgroup $I$ is attached to an alcove $\mathfrak{a}$ in the (semisimple) Bruhat-Tits building of $G$, with vertex the special point $\mathbf{x}_0$, and 
we let $\Sigma$ be the set of reflections across the  walls of $\mathfrak{a}$, so that $\Sigma_0$ appears as the subset of reflections across walls passing through $\bx_0$. Then ${\Sigma}$ generates an affine Weyl group   ${W}^a$ canonically identified with the subgroup $( \cn\cap \Ker w_G)/Z^0 $ of ${W}$; also ${W}$ is the semi-direct product of its normal subgroup ${W}^a$ and the subgroup $\Omega$ stabilizing the alcove $\mathfrak{a}$. We let $\ell$ be  the length function of  the Coxeter system $({W}^a, \Sigma)$  and we extend it to ${W}$, trivially on $\Omega$, i.e.\ so that $\ell(w\omega)=\ell(w)$ for $w\in W^a$, $\omega\in \Omega$; on $W_0$ it restricts to the length function of the Coxeter system $(W_0,\Sigma_0)$. Inflating through $_1{W}\rg {W}$ we get a length function on $_1{W}$ and $_1W_0$, still written $\ell$. The operators $T(w)$ in $\Hh_\Z(G,I)$ for $w\in {} _1 {W}$ satisfy the ``braid relations''
$$
T(w)T({w'})=T({ww'})\ \mathrm{when}\ \ell(ww')=\ell(w)+\ell(w'). \leqno(\mathrm{IV}.6.1)
$$
There are other relations, the ``quadratic relations'' \cite[Proposition~4.3]{Vig4}. Essentially there is one such relation for each $s\in \Sigma$.  It comes directly from the finite field case, treated in \cite[6.8]{CE}. For $s\in \Sigma_0$, $s=s_\beta$ for some $\beta\in \Delta$, we may describe the relation as follows: let $n_s$ be a lift of $s_\beta$ in $  \cn_k \cap M_{\beta,k}'$ and define $Z_{k,s} = Z_k\cap M_{\beta,k}'$ (so that $n_s^2$ belongs to $Z_{k,s}$); then the quadratic relation for $T(n_s)$ is
$$
T({n_s})(T({n_s})-c_{n_s}) = q_sT({n_s^2}), \leqno(\mathrm{IV}.6.2)
$$
where $q_s > 1$ is a power of $p$ and 
\[
c_{n_s}=\sum_{t\in Z_{k,s}} c_{n_s}(t)T(t)
\]
for positive integers $c_{n_s}(t)= c_{n_s}(-t)$, constant on each coset of $\{x s(x)^{-1} \ | \ x\in Z_k\}$, of sum $q_s-1$. Moreover, we have $c_{n_s}\equiv c_{s} \mod p$, where
$$
c_{s}:= (q_s-1) |Z_{k,s}|^{-1}\sum_{t\in Z_{k,s}} T(t). \leqno(\mathrm{IV}.6.3)
$$
We have  $T({n_s})c_{n_s}=c_{n_s}T({n_s})$.   \medskip

\noindent \textbf{Remark} In the $C$-algebra $\Hh(G,I)$, $q_s$ equals $0$  and $c_{n_s}$ equals $ -  |Z_{k,s}|^{-1} \sum_{t\in Z_{k,s}} T(t)$, so the relations simplify somewhat. We always embed the group algebra of $Z_k$ over $C$ into $\Hh(G,I)$ by sending
$t$ to $T(t)$; for $s=s_\beta$ as above we have $\psi(c_{n_s})=-1$ if $\psi$ is trivial on $Z_{k,s}$ (i.e.\ $\beta$ belongs to the set $\Delta(\psi)$ of \ref{III.8}, which contains $J$), and 
$\psi(c_{n_s})=0$ otherwise.

\medskip

\noindent \textbf{Proposition} \textit{
There is a unique extension of the map $s\mapsto n_s$ from $\Sigma_0$ to $\cn_k$ to a map $w \mapsto n_w$ from $W_0$ to $\mathcal{N}_k$ such that  $n_{ww'}=n_wn_{w'}$ for $w$, $w'$ in $W_0$ such that $\ell (ww')=\ell(w)+\ell(w')$.
}

\medskip

\noindent\textbf{Proof} (Another proof is in \cite[Proposition~3.4]{Vig4}.) Uniqueness is obvious, as we must have $n_w=n_{s_1}\cdots n_{s_r}$ for each reduced decomposition $w=s_s\cdots s_r$ of $w$ in $W_0$ with the $s_i$ in $\Sigma_0$. Existence will be consequence of \cite[\S1, n$^o$~5, Proposition~5]{Bk} once we prove:

($*$) For $s$, $s'$ distinct in $\Sigma_0$, and $m$ the order of $ss'$,
then $(n_sn_{s'})^\ell=(n_{s'}n_s)^\ell$ if $m=2\ell$
and $(n_sn_{s'})^\ell n_s = (n_{s'}n_s)^\ell n_{s'}$  if $m=2\ell+1$.

To prove ($*$) we may assume that $\bg_k$ is semisimple simply connected of relative rank 2, with $W_0$ generated by $s$ and $s'$, corresponding to the two simple roots $\beta$ and $\beta'$. But then the result follows from \cite[6.1.8]{BT1} applied to the valued root datum associated to $(\bg_k,\gs_k,\bb_k)$: indeed, we can always put reduced roots of $\Phi$ in a ``circular order'' as required by \cite[6.1.8]{BT1}, with $\beta$ first and $\beta'$ in the $m$-th position, in which case formula (9) of \cite[6.1.8]{BT1} gives exactly the required equality ($*$) above. $\square$

\medskip
Henceforward we use the extension $w\mapsto n_w$, and we put $\nu_w=n_{w^{-1}}^{-1}$ for $w\in W_0$; in particular if $w$, $w'$ in $W_0$ satisfy $\ell(ww')=\ell(w)+\ell(w')$, then $\nu_{ww'}=\nu_w\nu_{w'}$.

\subsection{}\label{IV.7} We are now ready to define $f_0$ (as promised in \ref{IV.5}) and study the action of $\Hh(K,I)$ on it. We let $w_0$ be the longest element in~$W_0$.

\medskip

\noindent\textbf{Definition} \textit{The function $f_0$ in $\cx^I$ has support $B \nu_{w_0}I$ and its value at $\nu_{w_0}$ is the function $e_\psi$ in $\ind _{Z^0}^Z\psi$ with support $Z^0$ and equal to $\psi$ on~$Z^0$.}

\medskip

Note the abuse of notation: we should choose a representative $\tilde{\nu}_{w_0}$ of $\nu_{w_0}$ in $\cn^0$ but neither the coset $B\tilde{\nu}_{w_0}I$ nor the value at $\tilde{\nu}_{w_0}$ depend on that choice. We shall allow similar abuse of notation below. Note also that $f_0$ depends on the choice of $\nu_{w_0}$ (but the support of $f_0$ is independent of this choice).

\medskip

\noindent \textbf{Notation} For $z\in Z_k$ and $w\in W_0$ we put $w\cdot z=n_wzn_w^{-1}$ (it is simply the natural action of $w\in W_0=\cn_k/Z_k$ on $Z_k$); more generally we shall use a dot to denote a conjugation action, which will be clear from the context.
 
\medskip

\noindent\textbf{Lemma} \textit{For $z\in Z_k$ we have $z^{-1}f_0=\psi({w_0}\cdot z^{-1})f_0=f_0 T(z)=\tau_{w_0\cdot z}f_0$.}

\medskip

The last equality in the lemma will be generalized below (\ref{IV.10}). 
\medskip

\noindent\textbf{Proof} Since $Z^0$ normalizes $I$, the first equality in the lemma comes from an immediate computation, whereas the equality $f_0T(z)= z^{-1}f_0$ comes from (IV.5.1). The equality $\tau_zf_0=\psi(z^{-1})f_0$ is equally easy. $\square$

\medskip

\noindent\textbf{Proposition} \textit{
Let $w\in W_0$. Then $f_0T({n_w})$ has support $B\nu_{w_0w}I$ and value $e_\psi$ at~$\nu_{w_0w}$.
}

\medskip

\noindent\textbf{Proof} As $f_0T(n_w)$ is $I$-invariant, it is enough to compute its value at $\nu_{w'}$ for $w'$ in $W_0$. By definition $(f_0T(n_w))(g) = \sum f_0(ghn_w^{-1})$ for $g\in G$, where the sum runs over $h$ in $I/(n_w^{-1}In_w \cap I)$. Assume that for such an $h$, $f_0$ is not $0$ at $\nu_{w'}hn_w^{-1}$. Then looking modulo $K(1)$, we get that $\nu_{w'}U_k^{\op} n_w^{-1} \cap B_k \nu_{w_0}U_k^{\op}$ is non-empty, and, multiplying on the right by $\nu_{w_0}^{-1}$, that $\nu_{w'} U_k^{\op} n_w^{-1}\nu_{w_0}^{-1} \cap B_k\not=\emptyset$ and hence $B_k \nu_{w'} U_k^{\op} \cap B_k \nu_{w_0}n_w U_k^{\op}\not= \emptyset$; by the Bruhat decomposition in $G_k$, that implies $w'=w_0w$. Assume that $w'=w_0w$; then $h$ belongs to $\nu_{w'}^{-1} B^0 \nu_{w_0} In_w$. However note that $\ell(w_0w)+\ell(w^{-1})=\ell(w_0)$ (because $w_0$ is the longest element in $W_0$), so that $\nu_{w'}\nu_{w^{-1}}=\nu_{w_0}$; we deduce that the image of $h$ in $G_k$ belongs to $n_w^{-1} B_k^{\op}n_w \cap U_k^{\op}=n_w^{-1}U_k^{\op}n_w\cap U_k^{\op}$. But that shows that $h$ belongs to $n_w^{-1} In_w \cap I$ and consequently $(f_0T(n_w))(\nu_{w_0w})=f_0(\nu_{w_0w}n_w^{-1}) = f_0(\nu_{w_0})= e_\psi$. $\square$

\medskip

\noindent\textbf{Corollary} $f=\sum\limits_{w\in w_0W_J}f_0 T({n_w})$.

\medskip

\noindent\textbf{Proof}  By the description in \ref{IV.4}, for $w$ in $W_0$, $f(\nu_w)$ is equal to $e_\psi$ if $w$ belongs to  $W_J$ and is $0$ otherwise: we only have to remark that $P_{J,k}U_k^{\op}= B_k W_J U_k^{\op}$, and 
since $\psi(Z_k\cap M_{\beta,k}')=1$ for $\beta\in J$, the character $\varepsilon$ of \ref{IV.4} is trivial on $\nu_w$ for $w\in W_J$. $\square$

\subsection{}\label{IV.8} We need to determine the action of $c_{n_s}$ on  $f_0T({n_w})$  for $s=s_\beta$, $\beta\in J$. We recall that $J\subset \Delta(\psi)$.

\medskip

\noindent\textbf{Proposition} \textit{Let $\beta\in \Delta(\psi)$, $s=s_\beta$ and $z\in Z_k\cap M_{\Delta(\psi),k}'$. For $w\in w_0W_{\Delta(\psi)}$, we have }
\begin{align*}
f_0\, T({n_w})T(z) &=f_0\, T({n_w})\ \mathit{and}\\
f_0\, T({n_w})c_{n_s} &=-f_0\, T({n_w}).
\end{align*}
 \textit{In particular $fc_{n_s}=-f$}.

\medskip

\noindent\textbf{Proof} By \ref{III.10} Example 2, $\psi$ is trivial on $Z_k\cap M_{\Delta(\psi),k}'$. By \ref{IV.7} Lemma then, we get $f_0T(z)=f_0$ for $z\in Z_k\cap {w_0}M_{\Delta(\psi),k}'w_0^{-1}$. The braid relation gives $T({n_w})T(t)=T(w\cdot t)T({n_w})$ for $t\in Z_k$, $w\in W_0$. For $z\in Z_k\cap M_{\Delta(\psi),k}'$ we have ${w\cdot z}\in Z_k\cap M_{\Delta(\psi),k}'$ for $w\in W_{\Delta(\psi)}$, hence $(w_0w)\cdot z\in Z_k\cap {w_0}M_{\Delta(\psi),k}'w_0^{-1}$, and consequently $f_0T({n_{w_0w}})T(z)=f_0T(n_{w_0w})$. That gives the first assertion.

The second one comes from the expression of $c_{n_s}$ in (IV.6.3), noting that $q_s$ gives $0$ in $C$; the last assertion follows from \ref{IV.7} Corollary. $\square$

\subsection{}\label{IV.9} 

\noindent\textbf{Notation } Let $w_J$ be the longest element in $W_J\subset W_0$ and put $w^J=w_0w_J$ (note that $w_J$ and $w_0$ have order $2$). We put $f_J=f_0T(n_{w^J})$.

\medskip
\noindent\textbf{Lemma 1} \textit{For $w\in W_J$ we have} 
(i) $\ell(w^Jw) =\ell(w^J)+\ell(w)$, 
(ii) $T(n_{w^Jw})=T(n_{w^J})T({n_w})$, and 
(iii) $f_0T(n_{w^Jw})=f_JT({n_w})$.

\medskip
\noindent\textbf{Proof} We have $\ell(w^Jw)=\ell(w_0w_Jw)=\ell(w_0)-\ell(w_Jw)$; if $w\in W_J$ we also have $\ell(w_Jw) =\ell(w_J)-\ell(w)$ so we get $\ell(w^Jw)=\ell(w^J)+\ell(w)$; by the braid relation 
 $T(n_{w^Jw})=T(n_{w^J})T({n_w})$, and the last assertion follows. $\square$
 \medskip
 
By   Lemma 1, and  \ref{IV.7} Corollary, \ref{IV.8} Proposition, we have 
 $f=\sum\limits_{w\in W_J}f_JT({n_w})$   and   for  $w\in W_J$ $$f_JT(n_w)c_{n_s}=-f_JT(n_w) . \leqno(\mathrm{IV}.9.1)
 $$
  For $s\in\Sigma_0$ we put $T^*({n_s})=T({n_s})-c_{n_s}$, so that in $\Hh_\Z(K,I)$ we get 
$$
T({n_s})T^*({n_s}) =T^*({n_s})T({n_s})=q_s T({n^2_s}) \quad (=0\ \mathrm{in}\ \Hh(K,I)).
$$
That definition can be extended to defining $T^*({n_w})$ for $w\in W_0$, so that $T^*({n_{ww'}})=T^*({n_w}) T^*({n_{w'}})$ if $\ell(ww')=\ell(w)+\ell(w')$ \cite[Proposition 4.13]{Vig4}. We now use the Bruhat order $\le$ on the Coxeter group~$W_J$ (see for example \cite{Deo}).

\medskip

\noindent\textbf{Proposition} \textit{ For $w\in W_J$ we have}
\begin{align*}
f_J \big(\sum\limits_{v\le w} T({n_v})\big) &= f_J T^*({n_w})\ \mathit{and\ in\ particular}\\
f=f_J T^*({n_{w_J}}) &=f_0 T({n_{w^J}})T^*({n_{w_J}}).
\end{align*}

\medskip A  similar proposition can be found in   \cite[Lemma 5.1]{Oll2}.

\medskip 
\noindent\textbf{Proof} We use induction on $\ell(w)$. The result is true for $w=1$. If $\ell(w)=\ell\ge1$, we write $w=w's$ with $\ell(w')=\ell-1, \ell(s)=1$.  As $\ell(w)=\ell(w')+\ell(s)$ we have $T^*(n_w)=T^*(n_{w'})T^*(n_s)=T^*(n_{w'})(T(n_s)-c_{n_s})$. By induction $f_JT^*(n_{w'})= \sum\limits_{v\le w'} f_JT(n_v)$. Remembering that for $v$ in $W_J$ we have $T(n_{w^J})T(n_v)=T(n_{w^Jv})$ and by (IV.9.1)  $f_JT(n_v)c_{n_s}=-f_JT(n_v)$. So finally we obtain
$$
f_JT^*({n_w})=f_J T^*({n_{w'}})(T(n_s)+1).
$$
By  induction $f_JT^*(n_{w'})=\sum\limits_{v\le w'}f_J T(n_v)$, so we want to compute  $A=\sum\limits_{v\le w'}f_J T(n_v)T(n_s)$.

Divide the set of $v\le w'$ in the disjoint union $X\sqcup Y\sqcup Ys$ where
$$
\begin{array}{lll}
Y&=&\{ v\in W_J,\ v<vs\le w'\},\\
Ys&=&\{ v\in W_J,\ vs<v\le w'\},\\
X&=&\{ v\in W_J,\ v\le w'\ \mathrm{and}\ vs \nleqslant w'\}.
\end{array}
$$
 In $A$, the subsum over $Y\sqcup Ys$ is
$$
\sum_{v\in Y} f_J(T(n_{vs}) + T(n_v))T(n_s).
$$
But for $v\in Y$, we have $v<v s$ so $T(n_{vs})=T(n_v)T(n_s)$ and 
$f_J(T(n_{vs})+T(n_v))T(n_s) = f_J T(n_v)( T(n_s)+1)T(n_s)$. By (IV.9.1) that equals  $f_JT({n_v})T^*(n_s)T(n_s)$  which is $0$ because $T^*(n_s)T(n_s)=0$ in $\Hh$.
 So $A=\sum\limits_{v\in X} f_J T(n_v)T(n_s)$. Since for $v\in X$, we have $v<vs$ we get $A=\sum\limits_{v\in X}f_J T(n_{vs}
)$.

The proof will be complete once we get:

\medskip

\noindent\textbf{Lemma 2} $Xs =\{v\in W_J$, $v\le w$ and $v\nleqslant w'\}$.

\medskip

\noindent\textbf{Proof} We use properties of the Bruhat order \cite[Theorem 1.1 (II) (ii)]{Deo}. Let $a$, $b$ in $W_J$ with $a\le b$. Then:

\medskip
(1) If $a<as$ then $a\le bs$;\ \ (2) if $b>bs$ then $as\le b$.
\medskip

Let $v\in X$, i.e.\ $v\le w'$, $vs \nleqslant w'$. Then by (2) applied to $a=v, b=w$, we get $vs\le w$. Conversely let $v\in W_J$ verify $v\le w$ and $v\nleqslant w'$; if $v<vs$ then $v\le w'$ by (1) applied to $a=v$, $b=w$, which is a contradiction; so $vs<v\le w$, which gives $vs\le w'$ by (1) applied to $a=vs$ and $b=w$. That proves the lemma. $\square$

\subsection{}\label{IV.10}

We now turn to the promised generalization of \ref{IV.7} Lemma which will be used in \ref{IV.15}.
\medskip

\noindent \textbf{Proposition} \textit{
Let $z\in Z$ with $z^{-1}\in Z^+$. Assume that $\nu_{w_0}\cdot z$ normalizes $\psi$. Then $f_0T(z)=\tau_{\nu_{w_0}\cdot z}f_0$.}

\medskip

\noindent\textbf{Remark} If $z^{-1}$ belongs to $Z^+$, $\nu_{w_0}\cdot z$ also belongs to $Z^+$, and conversely.

\medskip

\noindent \textbf{Proof} As both terms are $I$-invariant, we only need to check that they are equal at $\nu_w$ for $w\in W_0$. Now $(f_0T(z))(g)=\sum f_0 (ghz^{-1})$ for $g\in G$, where the sum runs over $h\in I/(z^{-1}Iz\cap I)$. But $I$ has an Iwahori decomposition and the assumption that $z^{-1}$ belongs to $Z^+$ gives $z^{-1}(I\cap U)z \subset I\cap U$, $z^{-1}(I\cap U_{\op})z \supset I\cap U_{\op}$, thus the inclusion of $I\cap U$ into $I$ induces of bijection of $(I\cap U)/(z^{-1} Iz\cap U)$ onto $I/(z^{-1}Iz\cap I)$, and it is enough to let $h$ run through $(I\cap U)/(z^{-1} Iz\cap U)$. For such an $h$, $\nu_w h z^{-1}$ belongs to $B\nu_{w_0}I$ only if $w=w_0$: indeed, $\nu_w hz^{-1} \in B n_wU$ and $B \nu_{w_0}I \subset B n_{w_0}U$, so the Bruhat decomposition in $G$ implies $w=w_0$. Consequently, both terms of the desired equality vanish at $\nu_w$ for $w\not=w_0$.

Consider now $(f_0T(z))(\nu_{w_0})$. Let $h\in I\cap U$ with $\nu_{w_0}hz^{-1}=b\nu_{w_0}j$ for some $b$ in $B$, $j$ in $I$; again by the Iwahori decomposition of $I$, we may assume that $j$ belongs to $I\cap U$ and then the equality $h=(\nu_{w_0}^{-1} b\nu_{w_0})z(z^{-1}jz)$, where $\nu_{w_0}^{-1} b \nu_{w_0}z \in B^{\op}$ and $z^{-1}jz \in  U$, shows that $h$ is equal to $z^{-1}jz$ and belongs to $z^{-1}Iz\cap U$; consequently, $(f_0T(z))(\nu_{w_0})=f_0(\nu_{w_0}z^{-1}) 
=f_0((\nu_{w_0}\cdot z^{-1})\nu_{w_0})=(\nu_{w_0}\cdot z^{-1})f_0(\nu_{w_0})$. That is equal to $(\nu_{w_0}\cdot z^{-1})e_\psi$, which sends $z'$ to $e_\psi(z'(\nu_{w_0}\cdot z^{-1}))$.
On the other hand if $\nu_{w_0}\cdot z$ normalizes $\psi$, we have $(\tau_{\nu_{w_0}\cdot z}f)(\nu_{w_0})= \tau_{\nu_{w_0}\cdot z}e_\psi$, sending $z'$ to $e_\psi((\nu_{w_0}\cdot z^{-1})z')$. That gives the result since $\nu_{w_0}\cdot z$ normalizes~$\psi$. $\square$

\subsection{}\label{IV.11}    To go further, we need more notation. We have the vector space $V_{\ad}=X_*(\gs_{\ad})\otimes \R$, where $\gs_{\ad}$ is the torus image of $\gs$ in the adjoint group $\bg_{\ad}$ of $\bg$,   the dominant Weyl chamber $\mathcal{D}^+=\{v\in V_{\ad},\ \beta(v)>0$ for $\beta\in \Delta\}$, and the antidominant Weyl chamber $\mathcal{D}^-=-\mathcal{D}^+ =w_0 \mathcal{D}^+$. 
We recall the natural map $\nu:Z \rg V_{\ad}$ used in \cite[3.3]{Vig4}: the action of $z\in Z$ on $V_{\ad}$ is via translation by $\nu(z)$. We remark that $\nu$ is the composite of $-v_Z:Z\to X_*(\gs )\otimes \R$ with $X_*({\bf S})\otimes \R\to V_{\ad}$.
 By \cite{Vig4}, $Z^+$ is the set
of $z\in Z$ such that $\nu(z)$ belongs to the closure of $\mathcal{D}^-$ (i.e.\ $\beta\circ \nu(z) \le 0$ for $\beta\in \Delta$). The map $\nu$ factors through $_1\Lambda$ and $\Lambda$, and we still write $\nu$ for the corresponding maps. 
\medskip

Note however that in citing  \cite[Ch.\ 5]{Vig4}, some care is needed:

Firstly, the roots in   \cite[Ch.\ 5]{Vig4}  are in the reduced root system $\Phi_a$ on $V_{\ad}$ attached to the collection of affine root hyperplanes in $V_{\ad}$  (it is denoted by $\Sigma$ in  \cite[Ch.\ 5]{Vig4}).  It is \textbf{not} in general the root system $\Phi$ attached to $(\bg_{\ad},\gs_{\ad})$.
Let us describe what is happening. The space $V_{\ad}=X_*(\gs_{\ad})\otimes\R$ is naturally a quotient of $X_*(\gs)\otimes \R$, and its dual $X^*(\gs_{\ad})\otimes\R$ appears as the subspace of $X^*(\gs)\otimes\R$ generated by the roots in $\Phi$, which are then the same for $(\bg,\gs)$ and $(\bg_{\ad},\gs_{\ad})$. The coroot in $V_{\ad}$ attached to a given root $\beta$ in $\Phi$ is the image of $\beta^\vee\in X_*(\gs)$, we also write it $\beta^\vee$. The root system $\Phi_a$ on $V_{\ad}$ can be described from $\Phi$ as follows. For each $\beta\in\Phi$, there is a positive integer $e_\beta$ such that $\Phi_a$ is the set of $\beta_a:=e_\beta\beta$ for $\beta\in \Phi$; in particular, $e_{2\beta}=e_\beta/2$ if $2\beta\in\Phi$. The root systems $\Phi_a$ and $\Phi$ share the same Weyl group $W_0$, and consequently the same Weyl chambers. The choice of Weyl chamber defining $\Phi^+$ also defines $\Phi_a^+$ and $\beta\mapsto  \beta_a$ gives a bijection of $\Delta$ onto the set $\Delta_a$ of simple roots in~$\Phi_a$. Note also that $(\beta_a \circ \nu)(\Lambda)\subset \mathbb Z$  and that the coroot in $V_{\ad}$ associated to $\beta_a\in \Phi_a$ is $\beta_a^\vee=e_\beta^{-1}\beta^\vee$.

\medskip

\noindent \textbf{Examples} 1) If $\bg$ is split, then $\Phi_a=\Phi$, $e_\beta=1$ for $\beta\in\Phi$.

2) For $G= \GL_r(D)$, where $D$ is a central division algebra over $F$, of finite degree $d^2$, then $e_\beta=d$ for all $\beta\in \Phi$.

3) Assume that $\bg$ is semisimple simply connected of relative rank 1. Then there is only one positive root $\beta$ and $\beta_a\circ \nu(\Lambda)=2\Z$ \cite[5.14]{Vig4}. Going back to the situation of \ref{III.16} with no condition on the reductive group $\bg$ we deduce that $\nu(a_\beta)=\beta_a^\vee$, since $\langle\beta_a, \beta_a^\vee\rangle=2$. In particular $v_Z(a_\beta)= -   e_\beta^{-1}\beta^\vee$.
\medskip

Secondly, the choice of Iwahori subgroup corresponds to a choice of alcove with vertex $\bx_0$, and positivity conditions are with respect to that choice. As we work with the ``lower'' pro-$p$ Iwahori subgroup $I$, the alcove with vertex $\bx_0$ which corresponds to $I$ is the one contained in $\mathcal{D}^-$,
so positive roots in  \cite[Ch.\ 5]{Vig4} correspond to negative roots here. In citing  \cite[Ch.\ 5]{Vig4} therefore we either have to exchange positive and negative roots or replace $\nu$ with $-\nu$; we choose the first solution. For example $\Sigma^+, \mathcal D^+$ in  \cite[Ch.\ 5]{Vig4} correspond to $\Phi_a^-, \mathcal D^-$ here.

\subsection{}\label{IV.12} Other bases  of $ \Hh_\Z(G,I)$  are constructed in \cite[Ch.\ 5]{Vig4} using (spherical) orientations. They generalize the Bernstein-Lusztig basis of an affine Hecke algebra. We need not know what such an object is, only that it is determined by a Weyl chamber in $V_{\ad}$; the action of $W_0$ on Weyl chambers determines an action on  orientations; but as in    \cite[Ch.\ 5]{Vig4}, we let $W_0$ (and hence $_1W$ via $_1W\to W_0$) act \textbf{on the right} on orientations by $(o,w) \mapsto o\cdot w$, so that if an orientation $o$ corresponds to the Weyl chamber $\mathcal{D}$, then $o\cdot w$ corresponds to $w^{-1}(\mathcal{D})$.

Let $o$ be an orientation. By \cite[Corollary 5.26]{Vig4} it gives a basis $(E_{o}(w))_{w\in {}_1\!  W}$ for $ \Hh_\Z(G,I)$.
In $ \Hh_\Z(G,I)$ some computations are easier because it is a ``characteristic zero'' algebra.  The above basis  of $ \Hh_\Z(G,I)$  specializes to a basis $(E_{o}(w))_{ w\in {}_1\!W}$ of  $\Hh$ over $C$: we use the same notation, making the context precise when necessary.

To $w\in {}_1W$ is attached an element $q_w$ in $\Z$, such that $q_{n_s}=q_s$ for $s\in \Sigma_0$ and $q_w=1$ if $\ell(w)=0$. The main relations in $ \Hh_\Z(G,I)$ satisfied by the  $E_{o}(w)$ are the following relations: for $w$, $w'$ in $_1W$,
$$
E_{o}(w) E_{o\cdot w}(w') =q_{w,w'} E_{o}(ww')\ \mathrm{with}\ q_{w,w'}=(q_w q_{w'}q_{ww'}^{-1})^{1/2}.
\leqno(\mathrm{IV}.12.1)
$$
\textbf{Beware} that in general $o\cdot w\not=o$, although it is the case when $w\in{}_1\Lambda$. Note  that $q_{w,w'}=1$   if and only  if $\ell(ww')=\ell(w)+\ell(w')$,  and   $q_{w,w'}$ gives $0$ in $C$ otherwise \cite[Remark 4.18 and Lemma~4.19]{Vig4}.

In particular, if $\ca _{o}$ is the subspace of $\Hh$ with basis $(E_{o}(\lambda))$ for $\lambda\in {}_1\Lambda$, the multiplication in $\ca_{o}$ is straightforward:
\begin{equation*} 
E_{o}(\lambda)E_{o}(\lambda') = \begin{cases}E_{o}(\lambda\lambda')\ &\mathrm{if}\ \ell(\lambda\lambda')=\ell(\lambda)+\ell(\lambda'),\\
 0 \ &\mathrm{otherwise}.
 \end{cases}
 \leqno(\mathrm{IV}.12.2)
\end{equation*}
 Thus $\ca_{o}$ is a subalgebra of $\Hh$. In fact the condition $\ell(\lambda\lambda')=\ell(\lambda)+\ell(\lambda')$ means that $\nu(\lambda)$ and $\nu(\lambda')$ belong to the same closed Weyl chamber in $V_{\ad}$ \cite[5.12]{Vig4}.

If $o$ is an orientation, we let $\Lambda_{o}$ be the set of $\lambda \in \Lambda $ such that $\nu(\lambda)$ belongs to the closure of the corresponding Weyl chamber; we similarly define $_1\Lambda_{o}$. For $\lambda$ in $_1\Lambda_{o}$, we have $E_{o}(\lambda)=T(\lambda)$ \cite[Example 5.30]{Vig4}.

We shall need the orientation $o_I$ attached to a subset $I$ of $\Delta$: by definition it is the orientation corresponding to the Weyl chamber $w_I(\mathcal{D}^-)$. Hence  $o_\Delta$ corresponds to $\mathcal{D}^+$,  $o_\emptyset$ corresponds to $\mathcal{D}^-$,  $o_I=o_\Delta \cdot w^I$, $\Lambda_{o_I}=  w_I \cdot  \Lambda^+$ (hence $_1\Lambda_{o_I}= \nu_{w_I}\cdot  {}_1\Lambda^+$).
For $w\in  W_I$ we then have $E_{o_I}(n_w)=T(n_w)$ \cite[Example 5.32]{Vig4}. (Note that $w_I(\mathcal{D}^-)$ here  equals  $w_I(\mathcal{D}^+)$ in  \cite{Vig4}, which corresponds to $o_{w_I (\Delta)}$ in \cite{Vig4}.)

\subsection{}\label{IV.13}

We  need some length  formulas  (\cite[Corollaries 5.10 and 5.11]{Vig4}). We have to be careful to remember that $\Sigma^+$ in \cite{Vig4}  corresponds to $\Phi_a^-$. For $\lambda\in \Lambda$, $w\in W_0$, we have

\medskip
\begin{align} 
 &\ell(w\cdot\lambda)= \sum_{\beta\in \Phi_a^+} |\beta\circ \nu(\lambda)|=\ell(\lambda) ,\tag{IV.13.1}\\
&\ell(w\lambda)= \sum\limits_{\beta\in\Phi_a^+\cap w^{-1}(\Phi_a^+)}|\beta\circ \nu(\lambda)| +\sum\limits_{\beta\in\Phi_a^+\cap w^{-1}(\Phi_a^-)}|\beta\circ \nu(\lambda)-1|, \tag{IV.13.2}\\
  &\ell(\lambda w)=
\sum\limits_{\beta\in\Phi_a^+\cap w(\Phi_a^+)}
|\beta\circ \nu(\lambda)| +
\sum\limits_{\beta\in\Phi_a^+\cap w(\Phi_a^-)}
|\beta\circ \nu(\lambda)+1|.\tag{IV.13.3}
\end{align}
Note that for $\beta\in\Delta$ and $w=s_\beta=w^{-1}$, $s_\beta$ permutes $\Phi_a^+-\{\beta_a\}$ and sends $\beta_a$ to $-\beta_a$ so $\Phi_a^+\cap s_\beta(\Phi_a^-)=\{\beta_a\}$.

\medskip

\noindent\textbf{Lemma} \textit{ Let $I\subset\Delta$. Then, for $\lambda\in \Lambda_{o_I}$, $\ell(w^I\lambda)=\ell(w^I)+\ell(\lambda)$.
}

\medskip

\noindent\textbf{Proof} By (IV.13.2) we need to check that $\beta\circ \nu(\lambda)\le 0$ for $\beta\in\Phi^+\cap (w^I)^{-1}(\Phi^-)$; but $\lambda\in \Lambda_{o_I}$ means that $\beta\circ \nu(\lambda)\ge0$ for $\beta\in w_I(\Phi^-)=(w^I)^{-1}(\Phi^+)$.  $\square$

\subsection{}\label{IV.14} 
An important result in this chapter is the following.\medskip

\noindent\textbf{Theorem} {\itshape\textit{
Let $w\in W_J$. Then for $\lambda\in{}_1\Lambda$, 
}
$$
f_J T^*(n_w) E_{o_J}(\lambda)\ \mathrm=\begin{cases} 
\tau((\nu_{w_J}n_w)\cdot \lambda)f_JT^*(n_w) &\text{if $(\nu_{w_J}n_w)\cdot\lambda \in {} _1 \Lambda^+$ and normalizes $\psi$},  \\
0 & \text{if $(\nu_{w_J}n_w)\cdot \lambda\notin {}_1\Lambda^+$.}
\end{cases}
$$}

 The proof of the theorem is in \ref{IV.15}--\ref{IV.18}.
 Taking $w=w_J$ we get by \ref{IV.9} Proposition:

\medskip
\noindent \textbf{Corollary} {\itshape\textit{For} $\lambda\in {}_1\Lambda$,
$$
f E_{o_J}(\lambda) = \begin{cases} \tau(\lambda)f &\text{if $\lambda\in {}_1\Lambda^+$ and normalizes $\psi$},\\
 0 & \text{if $\lambda\notin {}_1\Lambda^+$.}
\end{cases}$$
}

\noindent \textbf{Remarks}
\begin{itemize}
\item[1)]  We have used the notation $\tau(\mu)$ for $\mu\in {} _1 \Lambda^+$ to mean $\tau_z$ for $z\in Z^+$ with image $\mu\in {} _1 \Lambda^+$. The shift of indices is only for typographical convenience.
\item[2)] As $\psi$ extends to a character of $M_{J,k}$ by \ref{IV.4}, each $n_w$ for $w\in W_J$ normalizes $\psi$, and it follows that $\lambda$ normalizes $\psi$ if and only if so does $(\nu_{w_J}n_w)\cdot\lambda$.
\item[3)] The subspace of $\ca _{o_J}$ generated by the $E_{o_J}(\lambda)$ for $\lambda$ in $_1\Lambda$ normalizing $\psi$ is a subalgebra $\ca _{o_J}(\psi)$ of $\ca _{o_J}$ The map $\ca _{o_J}(\psi)\rg \Hh_Z(\psi)$ sending $E_{o_J}(\lambda)$ to $\tau((\nu_{w_J}n_w)\cdot\lambda)$ if $(\nu_{w_J}n_w)\cdot\lambda \in {} _1 \Lambda^+$ and to $0$ otherwise is an algebra homomorphism $\theta_{\nu_{w_J}n_w}$, and for $T\in \ca _{o_J}(\psi)$ we have
$$
f_JT^*(n_w)T = \theta_{\nu_{w_J}n_w}(T) f_J T^*(n_w).
$$
\item[4)] The theorem says nothing when $(\nu_{w_J}n_{w})\cdot \lambda\in {}_1\Lambda^+$ and does not normalize $\psi$.
We do not use this case.
\end{itemize}

\subsection{}\label{IV.15} 

 We treat first the case where $w=1$. Recalling that $f_J=f_0T(n_{w^J})$, we want to compute $f_0T(n_{w^J})E_{o_J}(\lambda) $.
By \ref{IV.12}, we have $T(n_{w^J})=E_{o_\Delta}(n_{w^J})$, so we look at $E_{o_\Delta}(n_{w^J})E_{o_J}(\lambda) $.

Assume first that $\nu_{w_J}\cdot \lambda$ belongs to $_1\Lambda^+$, i.e.\ that $\lambda$ belongs to $_1\Lambda_{o_J}$,
and that $\nu_{w_J}\cdot \lambda$ normalizes $\psi$. Then \ref{IV.13} Lemma gives $\ell(n_{w^J})+\ell(\lambda)=\ell(n_{w^J}\lambda)$, hence $E_{o_\Delta}(n_{w^J}\lambda)=E_{o_\Delta}(n_{w^J})E_{o_J}(\lambda)$. Since $\ell(n_{w^J}\cdot\lambda)=\ell(\lambda)$ by (IV.13.1), we also obtain $\ell(n_{w^J}\cdot\lambda)+\ell(n_{w^J})=\ell(n_{w^J}\lambda)$ hence $E_{o_\Delta}(n_{w^J}\lambda) = E_{o_\Delta}(n_{w^J}\cdot\lambda)E_{o_\Delta}(n_{w^J})$, and finally 
$$
E_{o_\Delta}(n_{w^J}\cdot\lambda)E_{o_\Delta}(n_{w^J})=E_{o_\Delta}(n_{w^J})E_{o_J}(\lambda) = T(n_{w^J})E_{o_J}(\lambda).
$$
We can apply \ref{IV.10} Proposition to $n_{w^J}\cdot\lambda$. Indeed, $\nu_{w_0}\cdot (n_{w^J}\cdot\lambda)=(\nu_{w_0}n_{w^J})\cdot \lambda$ and $\nu_{w_0}n_{w^J}=\nu_{w_J}$. Since by \ref{IV.12} $E_{o_\Delta}(n_{w^J}\cdot\lambda)=T(n_{w^J}\cdot\lambda)$, that gives  $f_0E_{o_\Delta}(n_{w^J}\cdot\lambda)=\tau(\nu_{w_J}\cdot\lambda)f_0$, so $\tau(\nu_{w_J}\cdot\lambda)f_J=f_0 E_{o_\Delta}(n_{w^J}\cdot\lambda)T(n_{w^J})=f_J E_{o_J}(\lambda)$, which is the desired formula when $\nu_{w_J}\cdot\lambda$ belongs to $_1\Lambda^+$.

Fix a regular such $\lambda$ and let $\lambda'\in {}_1\Lambda- {}_1\Lambda_{o_J}$. Then $E_{o_J}(\lambda) E_{o_J}(\lambda')=0$ by (IV.12.2), and $f_JE_{o_J}(\lambda)E_{o_J}(\lambda')=0$, implying $\tau(\nu_{w_J}\cdot\lambda)f_JE_{o_J}(\lambda')=0$. Since $\tau(\nu_{w_J}\cdot\lambda)$ is invertible in $\Hh_Z(\psi)$, we get $f_JE_{o_J}(\lambda')=0$, which is the formula we want for $\lambda'$. The theorem is proved for $w=1$.

\subsection{}\label{IV.16} We prove the theorem by induction on $\ell(w)$ (see  \cite[Section 5]{Oll} for  $\GL_n$). Let $\ell(w)=\ell\ge 1$, and write $w=w's$ with $\ell(w')=\ell-1$ and $s=s_\beta$ for some $\beta\in J$ -- note that $w'(\beta)\in \Phi^+$ since $\ell(w's)=\ell(w')+1$. In particular $n_w=n_{w'}n_s$ and $T^*(n_w)= T^*(n_{w'})T^*(n_s)$.

We need to investigate $T^*(n_s)E_{o_J}(\lambda)$ for $\lambda\in {}_1\Lambda$. Suppose we can prove
$$
f_J T^*(n_{w'})T^*(n_s)E_{o_J}(\lambda) = f_J T^*(n_{w'})E_{o_J}(n_s\cdot \lambda)T^*(n_s); \leqno(*)
$$
then the desired formula follows from the induction hypothesis. So we need to compare $E_{o_J}(n_s \cdot \lambda)T^*(n_s)$ and $T^*(n_s)E_{o_J}(\lambda)$.
By  \cite[Corollary 5.53]{Vig4} we have, for any orientation $o$ such that $\Ker \beta$ is a wall of the Weyl chamber corresponding to $o$:
\medskip

\begin{tabular}{llll}
\hskip-0.8cm(IV.16.1) &{If} &$\beta\circ \nu(\lambda)=0$, &$E_{o}(n_s\cdot \lambda)E_o(n_s)=E_o(n_s)E_o(\lambda)$;\\
&{if} &$\beta\circ \nu(\lambda)>0$, &$E_o(n_s\cdot \lambda)E_o(n_s)=E_{o\cdot s}(n_s)E_o(\lambda)$;\\
&{if} &$\beta\circ \nu(\lambda)<0$, &$E_o(n_s\cdot \lambda)E_{o\cdot s}(n_s)=E_o(n_s)E_o(\lambda)$.
\end{tabular}

\medskip

We now apply the results in  \cite[\S 5.4]{Vig4} to our case, where $o=o_J$. (We need to point out that since $\beta\in J$, $\Ker(\beta)$ is a wall of the Weyl chamber corresponding to $o_J$; also  \cite{Vig4} uses the notation $s$ for an element of $_1W_0$, where we use $n_s$, but we do have $n_s^2\in Z_k$ as required by \cite[5.35 and~5.36]{Vig4}.) Since $\beta\in J$, we have $E_{o_J}(n_s)= T (n_s)$ (\ref{IV.12}) and $E_{{o_J}\cdot s}(n_s)= T^*(n_s)$ by  \cite[Example 5.32]{Vig4}. So we get:

\medskip

\begin{tabular}{llll}
\hskip-0.8cm(IV.16.2) &{If} &$\beta\circ \nu(\lambda)=0$, &$E_{o_J}(n_s\cdot \lambda)T(n_s)=T(n_s)E_{o_J}(\lambda)$;\\
&{if} &$\beta\circ \nu(\lambda)>0$, &$E_{o_J}(n_s\cdot \lambda)T(n_s)=T^*(n_s)E_{o_J}(\lambda)$;\\
&{if} &$\beta\circ \nu(\lambda)<0$, &$E_{o_J}(n_s\cdot \lambda)T^*(n_s)=T(n_s)E_{o_J}(\lambda)$.
\end{tabular}

\subsection{}\label{IV.17} Accordingly we distinguish the three cases.

Assume first $\beta\circ \nu(\lambda)=0$; then formula ($*$) of \ref{IV.16} follows from (IV.16.2) and the following lemma.

\medskip

\noindent\textbf{Lemma} \textit{
Assume $\beta\circ \nu(\lambda)=0$. Then $E_{o_J}(n_s\cdot \lambda)c_{n_s}=c_{n_s}E_{o_J}(\lambda)$.}

\medskip 

\noindent\textbf{Proof} We work within the Levi subgroup $M_\beta$ of $G$. As $\beta\circ \nu(\lambda)=0$, $\lambda$ normalizes $K\cap M_\beta$ (\ref{III.7} Corollary). (Note that $K\cap M_\beta$ is the parahoric subgroup of $M_\beta$ attached to our special point $\bx_0$; $\lambda$ also normalizes the pro-$p$ radical $K(1)\cap M_\beta$ of $K\cap M_\beta$.) Consequently, $\lambda$ acts via conjugation on $M_{\beta,k}$; that action stabilizes $U_{\beta,k}$ and $U_{\beta,k}^{\op}$, so it also stabilizes the subgroup $M_{\beta,k}'$ they generate. Consequently, $\lambda$ acts via conjugation on $Z_{k,s}=Z_k\cap M_{\beta,k}'$. On the other hand, an element $t$ in $Z_{k,s}$ has length $0$, implying $E_{o_J}(n_s\cdot \lambda)T(t)=E_{o_J}((n_s\cdot \lambda)t)$ and $T(t)E_{o_J}(\lambda)=E_{o_J}(t\lambda)$. Now, computing in $_1W$, $(n_s\cdot \lambda)t\lambda^{-1}=(n_s\lambda n_s^{-1}\lambda^{-1})(\lambda t \lambda^{-1})$. As $t$ runs through $Z_{k,s}$, so does $\lambda t\lambda^{-1}$; on the other hand, by construction $n_s$ belongs to $M_{\beta,k}'$ so $n_s\lambda n_s^{-1} \lambda^{-1}$ belongs to $Z_{k,s}$. The result follows. $\square$

\subsection{}\label{IV.18} Assume now that $\beta\circ \nu(\lambda)<0$. Since $w'(\beta)$ is positive, $(w_J w's)(\beta)=-w_Jw'(\beta)$ is positive too. But $((w_Jw's)(\beta))\circ \nu$, evaluated on $\nu_{w_J}n_{w'}(\lambda)$ gives $(s(\beta)\circ \nu)(\lambda)=-\beta\circ \nu(\lambda)>0$ so $\nu_{w_J}n_{w'}(\lambda)$ is not in $_1\Lambda^+$, and consequently $f_JT^*(n_{w'})E_{o_J}(\lambda)=0$ by the induction hypothesis. But by (IV.16.2)
$$
f_JT^*(n_{w'})[T^*(n_s)E_{o_J}(\lambda) -E_{o_J}(n_s\cdot \lambda)T^*(n_s)]=
-f_J T^*(n_{w'})c_{n_s}E_{o_J}(\lambda).
$$
Since $f_JT^*(n_{w'})c_{n_s}=-f_JT^*(n_{w'})$ by \ref{IV.8} Proposition, 
$-f_JT^*(n_{w'})c_{n_s}E_{o_J}(\lambda)$ is equal to $f_JT^*(n_{w'})E_{o_J}(\lambda)$,  which is 0 by the above, and ($*$) is true in that case too. 

The case where $\beta\circ \nu(\lambda)>0$ is dealt with similarly: in that case we find 
$$
f_JT^*(n_{w'})[T^*(n_s)E_{o_J}(\lambda) - E_{o_J}(n_s\cdot\lambda)T^*(n_s)] = f_JT^*(n_{w'})E_{o_J}(n_s\cdot\lambda)c_{n_s}
$$
by (IV.16.2) and that is $0$ by induction because $\nu_{w_J}n_w(\lambda)$ is not in $_1\Lambda^+$ (as $(w_Jw(\beta))\circ \nu$ is positive on it). This completes the proof of \ref{IV.14} Theorem. $\square$

\subsection{}\label{IV.19} We now reach the easier part of our change of weight  (\ref{IV.19} Theorem (i)), which is a consequence of the following theorem.

\medskip

\noindent\textbf{Theorem} \textit{
Assume that $\lambda\in {}_1\Lambda$ normalizes $\psi$. Then}  
\begin{equation*} f'E_{o_{J'}}(\lambda n_{w_Jw_{J'}}^{-1})T^*(n_{w_Jw_{J'}})= 
\begin{cases}\tau(\lambda)f&\lambda\in {} _1 \Lambda^+\ \text{\textit{and}}\
\alpha\circ\nu(\lambda) < 0, \\ 0 &  \text{\textit{otherwise.}}
\end{cases}
\end{equation*}

\medskip

Taking $z\in Z^+$, normalizing $\psi$, and with $|\alpha|(z)<1$, we get  $f'T=\tau_z f$ for some $T$ in $\Hh$, which gives \ref{IV.1} Theorem (i). To prove the theorem, we first prove:

\medskip

\noindent\textbf{Lemma} $f'=f_J T^*(n_{w_Jw_{J'}w_J})T(n_{w_J w_{J'}})$.

\medskip

\noindent\textbf{Proof} By \ref{IV.7} Corollary, $f'=f_0 \sum\limits_{w\in w_0 W_{J'}}T(n_w)$, which can also be written as $f'=f_0 \sum\limits_{v\in W_{J'}}T(n_{w_0vw_{J'}})$. For $v$ in $W_{J'}$, write $w_0 vw_{J'}= w^J(w_J vw_J)(w_Jw_{J'})$. We have
$$
\ell(w_0 vw_{J'})=\ell(w_0)-\ell(vw_{J'})=\ell(w_0)-\ell(w_{J'})+\ell(v)\ \mathrm{(since}\ v\in W_{J'}),
$$
and
$$
\ell(w^J) = \ell(w_0)-\ell(w_J),\ \ell(w_Jvw_J)=\ell(v),\ \ell(w_Jw_{J'})=\ell(w_J)-\ell(w_{J'}),
$$
so $\ell(w_0 vw_{J'})= \ell(w^J) + \ell(w_J v w_J) +\ell(w_Jw_{J'})$. Consequently,
$$
\sum_{v\in W_{J'}} T(n_{w_0 v w_{J'}}) =T(n_{w^J})\big(\sum_{v\in W_{J'}} T(n_{w_J v w_J})\big) T(n_{w_Jw_{J'}})
$$
and $f'= f_J(\sum\limits_{v\in W_{J'}}T(n_{w_Jvw_{J}}))T(n_{w_Jw_{J'}})$.

Now $J''=-w_J(J')$ is a subset of $J$ and $w_JW_{J'}w_J=W_{J''}$; the element $w_Jw_{J'}w_J$ is the longest element of that group, hence
$$
\sum_{v\in W_{J'}} T(n_{w_J v w_J})=\sum_{v\le w_J w_{J'}w_J} T(n_v) .
$$
By \ref{IV.9} Proposition 
$$
f_J\big(\sum_{v\le w_Jw_{J'}w_J}T(n_v)\big)= f_J T^*(n_{w_J w_{J'}w_J})
$$
so $f'=f_J T^*(n_{w_Jw_{J'}w_J})T(n_{w_Jw_{J'}})$. $\square$

\medskip
\noindent\textbf{Proof of the theorem}
Put $v=w_Jw_{J'}$.
Note that since $v\in W_J$, $n_v\cdot \lambda$ normalizes $\psi$, see \ref{IV.14} Remark 2).

By the relations (IV.12.1) we get
\[
q_{n_v,\lambda n_v^{-1}}E_{o_J}(n_v\cdot\lambda)=E_{o_J}(n_v) E_{o_J\cdot v}(\lambda n_v^{-1}).
\]
On the other hand $E_{o_J}(n_w)=T(n_w)$ for $w\in W_J$, so we get
\[
T(n_v) E_{o_J\cdot v}(\lambda n_v^{-1})=q_{n_v,\lambda n_v^{-1}}E_{o_J}(n_v\cdot \lambda).
\]
We now compute
\begin{align*}
f' E_{o_J\cdot v}(\lambda n_v^{-1}) &=f_J T^*(n_{w_J w_{J'}w_J})T(n_v) E_{o_J\cdot v}(\lambda n_v^{-1})\\
&=q_{n_v,\lambda n_v^{-1}}f_J T^*(n_{w_J w_{J'}w_J}) E_{o_J}(n_v\cdot \lambda).
\end{align*}

By \ref{IV.14} Theorem we see that $f_JT^*(n_{w_Jw_{J'}w_J})E_{o_J}(n_v\cdot \lambda)$ is $0$ if $\lambda\notin {}_1\Lambda^+$.
If $\alpha\circ\nu(\lambda) = 0$, since $v(\alpha)\in\Phi^-$, $\ell(\lambda v^{-1}) > \ell(\lambda)-\ell(v^{-1})=\ell(v\cdot \lambda)-\ell(v)$ by \ref{IV.13}, so $q_{n_v,\lambda n_v^{-1}} = 0$.

Assume $\alpha\circ\nu(\lambda) < 0$ and $\lambda\in{}_1\Lambda^+$.
Let $\beta\in \Phi^+$ with $v(\beta)\in \Phi^-$. Since $v\in W_J$, $\beta$ is a linear combination of roots in $J$ (with non-negative integer coefficients). Moreover $w_{J'}(\beta)\in \Phi^+$, so the coefficient of $\alpha$ in $\beta$ is positive. Then for $\beta\in \Phi^+\cap v^{-1}(\Phi^-)$ we have $\beta\circ \nu(\lambda)\le \alpha\circ \nu(\lambda)$ by the above, so $\beta\circ \nu(\lambda)<0$, which implies by \ref{IV.13} that  $\ell(\lambda v^{-1}) = \ell(v\cdot \lambda)-\ell(v^{-1})$ and $f_JT^*(n_{w_Jw_{J'}w_J})E_{o_J}(n_v\cdot \lambda)=\tau(\lambda)f_J T^*(n_{w_Jw_{J'}w_J})$ by \ref{IV.14} Theorem (indeed,
$\ell(w_Jw_{J'}w_J)+\ell(v)=\ell(w_J)$ implies $n_{w_Jw_{J'}w_J}n_v = n_{w_J}$, so   $\nu_{w_J}n_{w_Jw_{J'}w_J}n_v=1$).
 The theorem follows on multiplying by $T^*(n_v)$, noting $T^*(n_{w_Jw_{J'}w_J})T^*(n_v)=T^*(n_{w_J})$ and $o_{J'}=o_J\cdot v$. $\square$

\subsection{}\label{IV.20} We now turn to the other part of the change of weight theorem (\ref{IV.19} Theorem (ii), (iii)), which is harder.
From now on, we put $s=s_\alpha$.

\medskip

\noindent\textbf{Lemma} $f'=f-f_JT^*(n_{w_Js})=f_J T^*(n_{w_Js})T(n_s)$.

\medskip

\noindent \textbf{Proof} By \ref{IV.9} we have $f=f_J\big(\sum\limits_{w\le w_J}T(n_w)\big)$  and  $f_J T^*(n_{w_Js})=f_J\Bigl(\sum\limits_{w\le w_J s}T_w\Bigr)$ so $f-f_J T^*(n_{w_Js})=f_J(\sum T(n_w))$ where the sum runs over $w$ in $W_J$ with $w\nleq w_Js$; but for $w\in W_J$, $w\le w_J s$ is equivalent to $s\le w_Jw$, so $w\nleq w_Js$ means that $w_Jw$ belongs to $W_{J'}$. Consequently, $f-f_JT^*(n_{w_Js})=f_J\big(\sum\limits_{w\in W_{J'}} T(n_{w_Jw_{J'}w})\big)=f_0 T(n_{w^J})\big(\sum\limits_{w\in W_{J'}} T(n_{w_Jw_{J'}w})\big)$. For $w$ in $W_{J'}$, 
$\ell(w_J w_{J'}w) = \ell(w_J)-\ell(w_{J'}w)=\ell(w_J)-\ell(w_{J'})+\ell(w)=\ell(w_J w_{J'})+\ell(w)$ so
 $T(n_{w_J w_{J'}w})=T(n_{w_Jw_{J'}})T(n_w)$. On the other hand, $\ell(w^J)+\ell(w_Jw_{J'})=\ell(w_0)-\ell(w_J)+\ell(w_J)-\ell(w_{J'})=\ell(w^{J'})$ so $T(n_{w^{J'}})=T(n_{w^J})T(n_{w_Jw_{J'}})$. It follows that $f-f_JT^*(n_{w_Js}) = f_0T(n_{w^{J'}})\big(\sum\limits_{w\in W_{J'}} T(n_w)\big)=f'$ by \ref{IV.7} Corollary applied~to~$J'$.
Moreover, as $\ell(w_Js)+\ell(s) = \ell(w_J)$ we have $T^*(n_{w_J})=T^*(n_{w_Js})T^*(n_s)$ and $f=f_JT^*(n_{w_J})=f_J T^*(n_{w_Js})(T(n_s)+1)$, as seen in \ref{IV.9} above, so $f'=f_J T^*(n_{w_Js})T(n_s)$. $\square$

\subsection{}\label{IV.21} Let now $\lambda\in {}_1\Lambda^+$ and put $\lambda'=n_s\cdot \lambda$. It is the element $f E_{o_J\cdot s}(n_s\lambda')$ that we want to relate to $f'$. To get an expression for it, we again need to distinguish cases, according to the integer $r=-\alpha_a \circ \nu(\lambda)\geq 0$ (recall that $\alpha_a$ is the simple root in $\Phi_a$ corresponding to $\alpha$). We first deal with the ``easy'' relations in~$\Hh$.

\medskip

\noindent\textbf{Lemma}\textit{
(i)  $\lambda'(n_s\cdot\lambda')=n_s\cdot(\lambda\lambda')\in {}_1\Lambda^+$.} 

\textit{(ii) If $r>0$, $\ell(n_s \lambda') = \ell(\lambda')-1$ and $T(n_s)E_{o_J\cdot s}(n_s\lambda')=T(n_s^2)E_{o_J} (\lambda')$.}

\textit{(iii) If $r\ge 2$, then $E_{o_J}(\lambda') E_{o_J}(n_s\lambda')=0$.}

\textit{(iv) If $r=1$, then $E_{o_J\cdot s}(n_s\lambda')=E_{o_J}(n_s\lambda')$ and}
$$
E_{o_J}(\lambda')E_{o_J}(n_s\lambda') =E_{o_J}(\lambda'(n_s\cdot \lambda')) T(n_s).
$$

\noindent\textbf{Proof} (i) The first equality is clear. Let us prove that $\lambda'(n_s\cdot \lambda')$ is in $_1\Lambda^+$. We have $\alpha_a \circ \nu(\lambda'(n_s\cdot\lambda'))=0$. For $\beta\in \Delta$, $\beta\not=\alpha$ we compute 
$$
\begin{array}{rl}
\beta_a \circ \nu(\lambda'(n_s\cdot \lambda')) &=\beta_a\circ \nu((n_s\cdot \lambda)(n_s^2\cdot\lambda))\\
&=(\beta_a+s(\beta_a))(\nu(\lambda)).
\end{array}
$$
It is $\le 0$ since $\beta_a,s(\beta_a) > 0$ and $\lambda$ is in $_1\Lambda^+$.
So we get (i).

(ii) Assume $r>0$. We need to work   in $\Hh_\Z$, and then specialize to $\Hh$. By (IV.13.2), we get $\ell(n_s\lambda')=\ell(\lambda')-1$ because $\alpha\circ \nu(\lambda')>0$. So the relation (IV.12.1) gives $E_{o_J\cdot s}(n_s) E_{o_J}(\lambda')=q_s E_{o_J\cdot s}(n_s\lambda')$. We  also have $E_{o_J}(n_s) E_{o_J\cdot s}(n_s) =q_sE_{o_J}(n_s^2)$, which gives
$$
q_s E_{o_J}(n_s)E_{o_J\cdot s}(n_s\lambda') = q_s E_{o_J}(n_s^2) E_{o_J}(\lambda').
$$
 Cancelling $q_s$, using  $E_{o_J}(n_s^2)=T(n_s^2)$, and specializing to $\Hh$ we get~(ii).

(iii) We proved $\ell(n_s\lambda') =\ell(\lambda')-1$ in (ii), so $\ell(\lambda')+\ell(n_s\lambda')=2\ell(\lambda')-1$. On the other hand $\lambda'n_s\lambda' =\lambda'(n_s\cdot \lambda') n_s$ and $\alpha_a \circ \nu(\lambda'(n_s\cdot \lambda'))=0$ so $\ell(\lambda' n_s \lambda')=\ell(\lambda'(n_s\cdot \lambda'))+1$   by (IV.13.3). But $\ell(\lambda'(n_s\cdot \lambda'))=\ell(\lambda')+\ell(n_s\cdot\lambda')-2r$   by (IV.13.1)  so we get $\ell(\lambda')+\ell(n_s\lambda')-\ell(\lambda' n_s\lambda')=2r-2$. This is $>0$ if $r\ge2$, so in that case $E_{o_J}(\lambda') E_{o_J}(n_s\lambda')=0$ by the relations (IV.12.1).

(iv) Assume now that $r=1$. The first formula is given by  \cite[Lemma 5.34]{Vig4}. In the proof of (ii) we have seen that
$ 
\ell(\lambda') + \ell (n_s\lambda')= \ell(\lambda' n_s \lambda')  $
  so we get 
$ 
E_{o_J}(\lambda') E_{o_J}(n_s\lambda') = E_{o_J}(\lambda' n_s \lambda')$. On the other hand
$\lambda'n_s \lambda'=\lambda'(n_s\cdot \lambda')n_s$ and we have seen $\ell(\lambda'n_s \lambda') =\ell(\lambda'(n_s\cdot\lambda'))+1$, so $E_{o_J}(\lambda' n_s\lambda') =E_{o_J}(\lambda'(n_s\cdot \lambda')) E_{o_J}(n_s) =E_{o_J}(\lambda'(n_s\cdot \lambda'))T(n_ s)$. $\square$

\subsection{}\label{IV.22}  In the sequel it is convenient to put $\varphi=f_JT^*(n_{w_Js})$ so that $f'=\varphi T(n_s)$, $f=\varphi +f'$. From \ref{IV.14} Theorem, we get the following: for $\mu\in {} _1  \Lambda$,
\begin{equation}\varphi E_{o_J}(n_s\cdot\mu)=\begin{cases} 
 \tau(\mu)\varphi  &\text{if $\mu\in {} _1  \Lambda^+$ and normalizes $\psi$},\\
 0 \ & \text{if $\mu\notin{}_1\Lambda^+$.}
 \end{cases}\tag{IV.22.1}
\end{equation}
 Put $E=E_{o_J\cdot s}(n_s\lambda')$ with $\lambda'=n_s\cdot \lambda$ as in \ref{IV.21} -- note that $\lambda'$ also normalizes~$\psi$.

By (ii) of \ref{IV.21} Lemma, $T(n_s) E=T(n_s^2) E_{o_J}(\lambda')$, so $\varphi T(n_s)E=\tau(n_s^2)\varphi E_{o_J}(\lambda')$ by (IV.22.1). But $\tau(n_s^2)\varphi=\varphi$ because $n_s^2$, which belongs to $Z_k\cap M_{\alpha,k}'$, acts trivially on $\varphi$ by \ref{IV.7} Lemma. We deduce $\varphi T(n_s)E=\varphi E_{o_J}(\lambda')=\tau(\lambda)\varphi$, again by (IV.22.1).  

We are now ready to prove a change of weight formula, in the special case where $\lambda\in {}_1\Lambda^+$ normalizes $\psi$ and $\alpha_a \circ \nu(\lambda)=-1$. Indeed, by (IV.22.1)
and (iv) of \ref{IV.21} Lemma we get $\tau(\lambda)\varphi E= \varphi E_{o_J}(\lambda')E=\varphi E_{o_J}(\lambda'(n_s\cdot \lambda'))T(n_s)$, hence $\tau(\lambda) \varphi E =\tau(\lambda\lambda')\varphi T(n_s)$, using again (IV.22.1).
We deduce that $\varphi E =\tau(\lambda')\varphi T(n_s)$, as $\tau(\lambda)$ is invertible~in~$\Hh_Z(\psi)$.

Consequently, $f E=\varphi E +\varphi T(n_s) E=\tau(\lambda) \varphi +\tau(\lambda') \varphi T (n_s) =\tau(\lambda)(f-f')+\tau(\lambda')f'$. We have proved:

\medskip

\noindent\textbf{Proposition} \textit{
Let $\lambda\in {}_1 \Lambda^+$ normalize $\psi$, and assume $\alpha_a \circ \nu(\lambda)=-1$. Then}
$$
\tau(\lambda) f-f E= (\tau(\lambda)-\tau(\lambda'))f'.
$$

\noindent\textbf{Remark} Note that $\tau(\lambda)f$  belongs to $\ind_K^G V$ because $\lambda\in {}_1 \Lambda^+$, so we see that $\ind_K^G V$ contains $(\tau(\lambda)-\tau(\lambda'))(\ind_K^G V')$. Note also that $\tau(\lambda)f'$ belongs to $\ind_K^G V'$ for the same reason;  but $\tau(\lambda')f'$ does not necessarily belong to $\ind_K^GV'$ because $\lambda'$ is not in $_1\Lambda^+$.

\subsection{}\label{IV.23}We now seek a similar formula in the case where $\lambda\in {}_1\Lambda^+$ normalizes $\psi$, $r=-\alpha_a \circ \nu(\lambda)\ge2$, still with $\lambda'=n_s\cdot\lambda$ and $E=E_{o_J\cdot s}(n_s \lambda')$. By  \cite[Proposition 5.48]{Vig4} we have, in    $\Hh_\Z$, an identity
$$
E_{o_J\cdot s}(n_s\lambda')-E_{o_J}(n_s\lambda') = \sum_{k=1}^{r-1} q(k,\lambda')q_s^{-1} c(k,\lambda') E_{o_J}(\mu(k,\lambda')) \leqno(*)
$$
and by \cite[Proposition 5.49]{Vig4}, in $\Hh$ only the terms for $k=1$ and $k=r-1$ may be non-zero, so we get, in~$\Hh$,
$$
E=E_{o_J}(n_s\lambda')+c_1 E_{o_J}(\mu_1\lambda') +c_{r-1} E_{o_J}(\mu_{r-1}\lambda'),
$$
where the last term disappears if $r=2$. For the moment we need not know what $c_1$, $c_{r-1}$ are in $C[Z_k]$, nor what $\mu_1$ and $\mu_{r-1}$ are in ${}_1\Lambda$ except that they do not depend on $\lambda$ and $\nu(\mu_k)=-k \alpha_a^\vee$ by \cite[formula (87)]{Vig4}, so $\nu(n_s^{-1}\cdot \mu_k)=k \alpha_a^\vee$. From that it follows that $(n_s^{-1}\cdot \mu_1)\lambda$ is in $_1\Lambda^+$, but not $(n_s^{-1}\cdot \mu_{r-1})\lambda$ if $r>2$.  Also by  \cite[5.49]{Vig4}, the $q$-terms in the  identity ($*$) above give 1 in $C$ for $k=1$ or $k=r-1$. 
Indeed, we have to show that $\ell(\lambda')-\ell(\mu_{-\alpha_a}^{-1}\lambda')=2$: remarking that $\nu(\mu_{-\alpha_a}^{-1}\lambda')=\nu(\lambda')-\alpha_a^\vee$, that comes from the length formula in \ref{IV.13}.
As in \ref{IV.22} we have $\varphi T(n_s)E=\tau(\lambda)\varphi$. On the other hand
$$
\varphi E=\varphi E_{o_J}(n_s \lambda') +\varphi c_1 E_{o_J}(\mu_1\lambda') +\varphi c_{r-1} E_{o_J}(\mu_{r-1}\lambda')
$$
where the last term disappears if $r=2$.

But $\tau(\lambda)\varphi=\varphi E_{o_J} (\lambda')$ by (IV.22.1), so $\tau(\lambda) \varphi E_{o_J} (n_s \lambda')= \varphi E_{o_J} (\lambda') E_{o_J}(n_s \lambda')$ which is $0$ by \ref{IV.21} Lemma (iii), and hence $\varphi E_{o_J} (n_s \lambda')=0$.
For $z\in Z_k$ we have $\varphi E_{o_J}(n_s\cdot z)=\tau_z\varphi=\psi(z^{-1})\varphi$ so we get $\varphi c_1 E_{o_J}(\mu_1\lambda')=\psi^{-1}(n_s^{-1}\cdot c_1) \varphi E_{o_J}(\mu_1\lambda')$, with the obvious notation for the conjugation action on $C[Z_k]$, and the obvious extension of $\psi^{-1}$ from $Z_k$ to $C[Z_k]$.
Similarly, if $r\ge 3$, $\varphi c_{r-1} E_{o_J} (\mu_{r-1}\lambda')=\psi^{-1}(n_s^{-1} \cdot c_{r-1}) \varphi E_{o_J}(\mu_{r-1}\lambda')$, which is  $0$ by (IV.22.1) because $(n_s^{-1}\cdot \mu_{r-1})\lambda$ is not in ${}_1\Lambda^+$. Thus   for $r\geq 2$,
$$
\varphi E = \psi^{-1}(n_s^{-1} \cdot c_1) \varphi E_{o_J}(\mu_1 \lambda').
$$
As $\varphi T(n_s)E=\tau(\lambda)\varphi$ we obtain:

\medskip

\noindent \textbf{Proposition} \textit{
Let $\lambda\in {}_1\Lambda^+$ normalize $\psi$, and assume $-\alpha_a \circ \nu(\lambda)\ge 2$. Then $$fE = \tau(\lambda)\varphi + \psi^{-1}(n_s^{-1}\cdot c_1)\varphi E_{o_J}(\mu_1\lambda').$$
}

\subsection{}\label{IV.24} We now apply the formulas given  by \ref{IV.22} Proposition and \ref{IV.23} Proposition to the case where $\lambda\in {}_1\Lambda^+$ normalizes $\psi$, and deduce \ref{IV.1} Theorem (ii) and (iii). We first assume $\alpha_a \circ \nu(\lambda)=-1$. As we have seen in \ref{IV.22} Remark, $\lambda'$ normalizes $\psi$ and $(\tau(\lambda)-\tau(\lambda'))(\ind_K^G V') \subset \ind_K^G V$.

\medskip

\noindent\textbf{Proposition} \textit{
Let $\lambda\in {}_1\Lambda^+$ normalize $\psi$, and assume $\alpha_a\circ \nu(\lambda)=-1$. Then $\psi$ is trivial on $Z^0 \cap M_\alpha'$ and $\tau(\lambda')=\tau(\lambda)\tau_\alpha$.
}

\medskip

\noindent\textbf{Proof} We work within $M_\alpha$. The semisimple Bruhat-Tits building of $M_\alpha$ is a tree, the apartment corresponding to $\gs$ is the line in $V_{\ad}$ generated by $\alpha_a^\vee$; the group $Z$ acts on that line via its quotient $\Lambda$, and $\lambda\in \Lambda$ acts via translation by $v$ with $\alpha_a\circ \nu(\lambda)=\alpha_a(v)$ and as $\alpha_a \circ \nu(\lambda)=-1$, $\lambda$ sends the (special) vertex $\bx_0$ to the adjacent (special) vertex $\bx_1=\bx_0-\frac{1}{2} \alpha_a^\vee$ in the apartment. We shall later prove the following claim.

For the claim the situation is the following: 

\medskip 
\noindent\textbf{Assumption} Assume that $\bg$ has relative semisimple rank 1, and let $\bx_1$ be a vertex in  $V_{\ad}$ (a line) adjacent to $\bx_0$, and  $K_1$ the corresponding (special) parahoric subgroup of $G$. Let $\bg_{1,k}$ be the group over $k$ attached to the parahoric subgroup $K_1$. (Note that both $K=K_0$ and $K_1$ contain $Z^0$ and $G_k$, $G_{1,k}$ contain $Z_k$.)

\medskip

\noindent\textbf{Claim}\textit{ The subgroup of $Z_k$ generated by $Z_k\cap G_k'$ and $Z_k\cap G_{1,k}'$ is  the image of $Z^0 \cap G'$ in~$Z_k$.}

\medskip

We apply the claim to $\bm_\alpha$. Since $\lambda$ sends $\bx_0$ to $\bx_1$, it conjugates $K_0$ to $K_1$, and conjugation by $\lambda$ induces an isomorphism of $M_{\alpha,k}$ onto $M_{\alpha,1,k}$ and of $M_{\alpha,k}'$ onto $M_{\alpha,1,k}'$. As $\psi$ is trivial on $Z_k\cap M_{\alpha,k}'$ by hypothesis, and $\lambda$ stabilizes $\psi$, $\psi$ is also trivial on $Z_k\cap M_{\alpha,1,k}'$ and by the claim $\psi$ is trivial on  $Z^0\cap M_\alpha'$. By  the second line after formula (87) in \cite{Vig4},  from $\alpha_a \circ \nu(\lambda)=-1$ we get $\nu(\lambda^{-1}\lambda')=\alpha_a^\vee$; but $\lambda'=n_s\cdot \lambda$ by definition, so $\lambda^{-1}\lambda'=\lambda^{-1}n_s\lambda n_s^{-1}$. Take $z\in Z$ with image $\lambda$ in $_1\Lambda$ and $\tilde{n}_s$ in $K\cap M_\alpha'\cap \cn$ with image $n_s$ in $M_{\alpha,k}$ (the existence follows from \ref{III.7} Lemma, for instance). Since $M_\alpha'$ is normal in $M_\alpha$, $z^{-1}\tilde{n}_s z$ is in $M_\alpha'$ so $\lambda^{-1}n_s \lambda n_s^{-1}$ is the image in $_1\Lambda$ of an element of $Z\cap M_\alpha'$. It follows that we can take $\lambda^{-1}\lambda'$  as the image in $_1\Lambda$  of  $a_\alpha$ of \ref{III.16} Notation (which verifies $\nu(a_\alpha)=\alpha_a^\vee$, cf.~\ref{IV.11} Example 3), and then $\tau(\lambda')=\tau(\lambda)\tau_\alpha$. $\square$

\medskip From the above proposition and \ref{IV.22} Proposition, we get case (iii) of \ref{IV.1} Theorem when $\alpha_a \circ \nu(\lambda)=-1$.

\medskip

\noindent\textbf{Corollary} \textit{
Let $\lambda\in {}_1\Lambda^+$ normalize $\psi$, and assume $\alpha_a \circ \nu(\lambda)=-1$. Then $\psi$ is trivial on $Z^0 \cap M_\alpha'$ and $\tau(\lambda)(1-\tau_\alpha)\ind_K^G V'\subset \ind_K^{G}V$.}

\medskip
We note that $\lambda a_\alpha\notin Z^+$ so in particular $\tau(\lambda)(1-\tau_\alpha)\notin \cz_G$.

\subsection{}\label{IV.25}  We  investigate the  term $\psi^{-1}(n_s^{-1}\cdot c_1)\varphi E_{o_J}(\mu_1\lambda')$ in \ref{IV.23} Proposition.

\medskip

\noindent\textbf{Proposition} {\itshape Let $\lambda\in {}_1\Lambda^+$ normalize $\psi$, and assume $-\alpha_a \circ \nu(\lambda)\geq 2$. 

(i) The element $n_s^{-1}\cdot\mu_1 \in {}_1\Lambda$ is in the image of $Z\cap M_\alpha'$.

(ii) If $\psi$ is not trivial on $Z^0\cap M_\alpha'$, then $\psi^{-1}(n_s^{-1}\cdot c_1)=0$.

(iii) If $\psi$ is trivial on $Z^0\cap M_\alpha'$, then $\psi^{-1}(n_s^{-1}\cdot c_1)=-1$ and $\tau((n_s^{-1}\cdot \mu_1)\lambda)=\tau(\lambda) \tau_\alpha $.

}
\medskip

Note that from (i) and \ref{III.16} Proposition (i), $n_s^{-1}\cdot\mu_1$ normalizes $\psi$ if $\psi$ is trivial on $Z^0\cap M'_\alpha$. In particular, in (iii) the element $\tau((n_s^{-1}\cdot \mu_1)\lambda)$ is defined. Using \ref{IV.23} Proposition and (IV.22.1) we get
$$
	fE = 
	\begin{cases}
	\tau(\lambda)(f - f') & \text{if $\psi$ is not trivial on $Z^0\cap M_\alpha'$},\\
	\tau(\lambda)(1 - \tau_\alpha)(f - f') & \text{if $\psi$ is trivial on $Z^0\cap M_\alpha'$}.
	\end{cases}
$$
This formula immediately yields \ref{IV.1} Theorem (ii), (iii) when  $-\alpha_a \circ \nu(\lambda)\ge 2$  (note that this implies $\lambda a_\alpha \in {}_1\Lambda^+$):

\medskip

\noindent\textbf{Corollary}  \textit{
Let $\lambda\in {}_1\Lambda^+$ normalize $\psi$, and assume $- \alpha_a \circ \nu(\lambda)\geq 2$. }

\textit{
(i) If $\psi$ is not trivial on $Z^0\cap M_\alpha'$ then $\tau(\lambda) \ind_K^G V' \subset \ind_K^G V$.
}

\textit{(ii) If $\psi$ is trivial on $Z^0\cap M_\alpha'$ then }
$$
\tau(\lambda)(1-\tau_\alpha)\ind_K^G V' \subset \ind_K^G V.
$$

To prove the proposition we need to know precisely what $c_1$ and $\mu_1$ are.
We have to distinguish cases: $\alpha_a \circ \nu(\Lambda)=\delta\Z$ for $\delta=1$ or $2$  \cite[Remark 5.3]{Vig4}. The generic case is $\delta=1$, which we tackle first. In that case choose $\lambda_s\in \Lambda$ with $\alpha_a\circ \nu(\lambda_s)=1$; then $\mu_1=(n_s\cdot \lambda_s)\lambda_s^{-1}$ and $c_1=(n_s\cdot\lambda_s)\cdot c_{n_s}$. Recall that $c_{n_s}=\frac{-1}{|Z_{k,s}|} \sum\limits_{z\in Z_{k,s}}z$ in $C[Z_k]$. In particular,
$\psi^{-1}(n_s^{-1}\cdot c_1)= \frac{-1}{|Z_{k,s}|} \sum\limits_{z\in Z_{k,s}} \psi^{-1}(\lambda_s\cdot z)$.
So we see that $\psi^{-1}(n_s^{-1}\cdot c_1)$ is non-zero if and only if $\psi$ is trivial on $\lambda_s Z_{k,s}\lambda_s^{-1}$, in which case it is equal to $-1$. Reasoning as in \ref{IV.24} with $\lambda_s$ instead of $\lambda$ we see that $\psi^{-1}(n_s^{-1}\cdot c_1)\not=0$ if and only if $\psi$ is trivial on $Z^0 \cap M_\alpha'$ and the other assertions of the proposition are obtained as in \ref{IV.24} as well (when $\delta=1$), noting that $\tau_\alpha$ is in the centre of $\Hh_Z(\psi)$.

\subsection{}\label{IV.26} We continue the proof of \ref{IV.25} Proposition. Now assume that $\delta=2$. One situation where this may happen is when $\bg$ has relative semisimple rank 1, or more generally when the connected component of the relative Dynkin diagram of $\bg$ containing $\alpha$ has rank 1. In that case, let $\tilde{s}$ be the reflection in the affine Weyl group of $\bm_\alpha$ corresponding to the affine root $\alpha_a+1$; it corresponds to a vertex $\bx_1$ in the semisimple Bruhat-Tits building of $\bm_\alpha$ (a tree) adjacent to the vertex $\bx_0$. As in \ref{IV.24} we let $K_1$ be the parahoric subgroup of $M_\alpha$ corresponding to the vertex $\bx_1$ (which is special), and $K_1(1)$ its pro-$p$ radical. Then $Z\cap K_1=Z^0$, $Z\cap K_1(1)=Z(1)$. The image of $\cn \cap K_1$ in $K_1/K_1(1)=M_{\alpha,1,k}$ is the group $\cn_{1,k}$ of $k$-points of the normalizer of $\bz_k$ in $\bm_{\alpha,1,k}$
 and we can choose in $\cn_{1,k}$ a lift $n_{\tilde{s}}$ of $\tilde{s}$ which actually belongs to $M_{\alpha,1,k}'$ -- note that $\tilde{s}$ generates $(\cn \cap K_1)/Z^0$ which we identify, via reduction with $\cn_{1,k}/Z_k$. Then, inside $_1W=\cn/Z(1)$, we can take (cf.\ \cite[Notation 5.37]{Vig4}) $\lambda_s=n_sn_{\tilde{s}}$, $\mu_1=\lambda_s^{-1}$, $c_1=c_{\tilde{s}}n_s^2$, where $c_{\tilde{s}}=\frac{-1}{|Z_{k,\tilde{s}}|} \sum\limits_{z\in Z_{k,\tilde{s}}} z$, with $Z_{k,\tilde{s}}= Z_k \cap M_{\alpha,1,k}'$\footnote{In principle those elements are defined in \cite{Vig4} with respect to $\bg$, not $\bm_\alpha$, but the above choices in $M_\alpha$ also work in $G$. The same remark applies in~\ref{IV.27}.}. We see that $\psi^{-1}(n_s^{-1}\cdot c_1)\not=0$ if and only if $\psi$ is trivial on $Z_{k,\tilde{s}}$. As $\psi$ is already trivial on $Z_{k,s}$, we get by  \ref{IV.24} Claim that $\psi^{-1}(n_s^{-1}\cdot c_1)\not=0$ if and only if $\psi$ is trivial on $Z^0 \cap M_\alpha'$, in which case $\psi^{-1}(n_s^{-1}\cdot c_1)=-1$. On the other hand,  $n_s^{-1}\cdot \mu_1$ is in the image of $Z \cap M_\alpha'$ (by lifting $n_s$ and $n_{\tilde{s}}$ to $\cn \cap M_\alpha'$ as in \ref{IV.24}). Moreover, by construction $\nu(\mu_1)=-\alpha_a^\vee$ and as in \ref{IV.24} we deduce that we can take the image of $a_\alpha$ in $_1\Lambda$ to be  $n_s^{-1}\cdot \mu_1$   and that $\tau((n_s^{-1}\cdot \mu_1)\lambda)= \tau(\lambda)\tau_\alpha$ if $\psi$ is trivial on $Z^0\cap M_\alpha'$.

\subsection{}\label{IV.27}  The only other case when $\delta=2$ may happen  is when the connected component of the Dynkin diagram of $\Phi_a$ containing $\alpha$ has type $C_n$, $n\ge 2$, and $\alpha$ is a long root \cite[Proposition 5.14]{Vig4}. Let then $\tilde{\alpha}_a$ be the highest root in $\Phi_a^+$ lying in the same component as $\alpha$, and $\tilde{s}$ be the reflection associated with $\tilde{\alpha}_a+1$. Then (cf.\ \cite[Lemma 5.15 and Notation 5.37]{Vig4}) $\mu_{-\alpha_a}=sw\tilde{s} w^{-1}$ for some $w\in W^a$ such that $\ell(\mu_{-\alpha_a})=2\ell(w)+2$ and $w\tilde{s} w^{-1}$ is the reflection $s'$ associated with the affine root $\alpha_a+1$ (whereas $s$ is associated with $\alpha_a$). Moreover $\mu_{-\alpha_a}=ss'$ satisfies $\nu(\mu_{-\alpha_a})=\alpha_a^\vee$. In that case (cf.\ \cite{Vig4}) $c_1=(w\cdot c_{\tilde{s}})n_s^2$ and $\lambda_s=n_s(w\cdot n_{\tilde{s}})$, $\mu_1=\lambda_s^{-1}$ with $n_{\tilde{s}}$, $c_{\tilde{s}}$ defined similarly as before \cite[\S4]{Vig4}; but conjugating by $w$ yields $w\cdot c_{\tilde{s}}=c_{s'}$ and $w\cdot n_{\tilde{s}}=n_{s'}$ where now $c_{s'},n_{s'}$ have a similar meaning, but in the relative semisimple rank 1 group $M_\alpha$. The same reasoning as in \ref{IV.26} then gives the desired result. 

\subsection{}\label{IV.28} To finish the proof of  \ref{IV.25} Proposition 
  we need only prove   \ref{IV.24}  Claim. It is convenient to deal first with the case where  $\bg=\bg^{\is}$. Then $W=W^a$ is generated by the involutions $s_0$ (generating $\cn^0/Z^0$) and $s_1$ (generating $(\cn \cap K_1)/Z^0)$. As $s_0s_1$ acts as a non-trivial translation on the apartment, $s_0s_1$ has infinite order.

Identify $\cn^0/Z(1)$ with $\cn_k$ and similarly $(\cn\cap K_1)/Z(1)$ with the group $\cn_{1,k}$ of $k$-points of the normalizer of $\bz_k$ in $\bg_{1,k}$. Choose a lifting $n_0$ of $s_0$ in $\cn_k\cap G_k'\subset {}_1W$ and a lifting $n_1$ of $s_1$ in $\cn_{1,k}\cap G_{1,k}'\subset {}_1W$. An element $w$ of $W$ has a unique reduced expression $w=\sigma_1\cdots \sigma_h$ with $\sigma_i=s_0$ or $s_1$ and we put $n_w=x_1\cdots x_h$ with $x_i=n_0$ if $\sigma_i=s_0$, $x_i=n_1$ if $\sigma_i=s_1$. We let $X$ be the subgroup of $Z_k$ generated by $Z_k\cap G_k'$ and $Z_k \cap G_{1,k}'$, and put $Y=\{n_wx\mid w\in W,x\in X\}$.

\medskip

\noindent \textbf{Lemma 1} \textit{$X$ and $Y$ are normal subgroups of $_1W$.}

\medskip

\noindent \textbf{Proof} Let $x\in Z_k$; then $n_0^{-1} x n_0 x^{-1}$ belongs to $Z_k$; but $Z_k$ normalizes  $G_k'$ so $n_0^{-1}x n_0 x^{-1}$ belongs to $Z_k \cap G_k'$. Similarly $n_1^{-1}xn_1x^{-1}$ belongs to $Z_k\cap G_{1,k}'$. 
In particular, $n_0$ and $n_1$ normalize $X$. Since $Z_k$ also normalizes $X$, so $_1W$ itself normalizes $X$. As $n_0^2$ and $n_1^2$ belong to $X$, we deduce that for $w$, $w'\in  W$ $n_wn_{w'}\in n_{ww'}X$ and $n_w^{-1}\in n_{w^{-1}}X$, so $Y$ is indeed a normal subgroup of $_1W$. $\square$

\medskip

Now let $H=IYI$ with the usual abuse of notation.

\medskip

\noindent\textbf{Lemma 2} \textit{
$H$ is a normal subgroup of $G$ and $(H\cap Z^0)/Z(1)=X$
}

\medskip

\noindent \textbf{Proof}
We first prove that $H$ is a subgroup of $G$. By Lemma 1, $H$ is closed under inverses. Working in $\Hh_\Z$, it is enough to show that for $y$, $y'$ in $Y$, the product $T(y)T(y')$ in $\Hh_\Z$ is a linear combination of $T(y'')$ for $y''$ in $Y$. But that is given by the relations in $\Hh_\Z$: the braid relations and the two quadratic relations $T(n_i)^2=q_i T(n_i^2)+ c_iT(n_i)$ where $q_i\in \Z$ and $c_i\in \Z [Z_k \cap G_{i,k}']$ for $i=0,1$.

As $Z^0$ normalizes $I$, and $Z_k$ normalizes $Y$, $Z^0$ normalizes $H$. The normalizer of $H$ contains $n_0$, $n_1$ (which belong to $H$), $Z^0$ and $I$, so it is $G$ itself. If an element $x$ of $H$ in a class $IyI$, $y\in Y$, is in $Z^0$ then $y$ has to belong to $Z_k$ so by the very definition of $Y$, $y$ belongs to $X$ and $x$ itself has image $y$ in $Z^0/Z(1)=Z_k$. That gives the last assertion of the lemma.
$\square$

\medskip

Clearly $H$ is not central in $G$, so $H=G$ because the only non-central normal subgroup of $G$ is $G$ itself (\ref{II.3} Proposition). But then $H\cap Z^0=G\cap Z^0=Z^0$ so $X=Z_k$, which gives the claim for $\bg=\bg^{\is}$.

 Let us now prove   \ref{IV.24}  Claim in the general case. We show first that the claim is equivalent to
$$
 Z(1)(Z^0\cap G')=Z(1)\big\langle{Z^0\cap \langle{U^0, U_{\op}^0}\rangle ,   Z^0\cap \langle{U\cap K_1, U_{\op}  \cap K_1}\rangle }\big\rangle. \leqno{(*)}
 $$
 It suffices to show that the image of $Z^0\cap \langle{U^0, U_{\op}^0}\rangle$ in $Z_k$ equals $Z_k\cap G'_k$ (and similarly  for the other term). It is clear that an arbitrary element of $Z_k\cap G'_k$ lifts to an element of  $\langle{U^0, U_{\op}^0}\rangle \cap Z^0 K(1)$. Using the Iwahori decomposition of $K(1)$ (\ref{III.7}) we can modify the lift so that it is contained in $Z^0\cap \langle{U^0, U_{\op}^0}\rangle$.   

 The only non-trivial part of the equality  $(*)$ is the inclusion $\subset$. The inclusion is true for $G^{\is}$, and we deduce it for $G$ by  applying the natural homomorphism
 $\iota:G^{\is}\rg G$,    using that $(Z^{\is})^0=\iota^{-1}(Z^0)$  (\ref{III.19} Proposition) and  that $Z(1)$ is the pro-$p$ Sylow of $Z^0$.
 This completes the proof of \ref{IV.24} Claim and hence of \ref{IV.1} Theorem. $\square$

\section{Universal modules}\label{V}

\subsection{}\label{V.1} In this chapter our goal is, for an irreducible representation $V$ of $K$, to study the ``universal'' representation $\ind_K^G V$ as a module over the centre $\cz_G(V)$ of the Hecke algebra $\Hh_G(V)$. In fact that structure is difficult to elucidate, so we consider various algebra homomorphisms $\chi:\cz_G(V) \rg A$ and the corresponding $A$-module $A\otimes_\chi \ind_K^G V$.  As an application, for a character $\chi:\cz_G(V)\rg C$, we prove Theorem~6 of the introduction -- used in Chapter~\ref{III} at the end of our classification -- which gives a nice filtration of $C\otimes_\chi \ind_K^GV$ as a representation of $G$. In this chapter we fix an irreducible representation $V$ of $K$ and let $(\psi,\Delta(V))$ be its parameter   as defined in \ref{III.9}.

\paragraph{{\large A) Freeness of the supersingular quotient of $\ind_K^GV$}}
\addcontentsline{toc}{section}{\ \ A) Freeness of the supersingular quotient of \texorpdfstring{$\ind_K^GV$}{ind\_K\^{}G V}}

\subsection{}\label{V.2} Until \ref{V.11} we fix a parabolic subgroup $P=MN$ of $G$ containing $B$.  Recall  from \ref{III.4} the subgroup $Z_{\Delta_M}^\bot $ of $Z$ consisting of those $z\in Z$ with $|\beta|(z)=1$ for all $\beta\in \Delta_M$. We write $Z^{+M}$ for the set of $z\in Z$ with $|\beta|(z)\le 1$ for $\beta\in \Delta_M$. Recall  from \ref{III.4} that $\cz_Z(V_{U^0})$ is spanned by the $\tau_z$ for $z\in Z_\psi$,  and that the natural image of $\cz_M(V_{N^0})$ in $\cz_Z(V_{U^0})$ (via $\cs_Z^M$) is spanned by the $\tau_z$ for $z\in Z^{+M}\cap Z_\psi$ -- we identify $\cz_M(V_{N^0})$ with that image.

\medskip

\noindent \textbf{Notation} We let $R_M$ be the quotient of $\cz_M(V_{N^0})$ by the ideal of elements  supported     on $(Z^{+M}\cap Z_\psi)-  Z_{\Delta_M}^\bot $.

\medskip
As $\cz_M(V_{N^0})$  is viewed as  a subset of  $\cz_Z(V_{U^0})$, we emphasize that the  supports above  are subsets of $Z$. 
Note that the elements of $\cz_M(V_{N^0})$ supported on $Z_{\Delta_M}^\bot $ form a subalgebra which maps isomorphically onto $R_M$. 
\medskip

Our first main result in this chapter is:

\medskip

\noindent\textbf{Theorem} \textit{Let $P=MN$ be a parabolic subgroup of $G$ containing $B$. Then $R_M\otimes_{\cz_G(V)} \ind_K^GV$ is free over $R_M$, where the tensor product is via the composite map $\cz_G(V) \rg \cz_M(V_{N^0})\rg R_M$.}

\medskip
We call $R_M\otimes_{\cz_G(V)} \ind_K^GV$ the supersingular quotient of $\ind_K^G V$
(cf.\ \ref{III.4} Corollary).

The proof of that theorem is rather long (\ref{V.3} to \ref{V.11}). We first treat the case where $P=G$ (\ref{V.3} Proposition). The proof then proceeds by comparing  with situations with a more regular weight (i.e.\ smaller $\Delta(V)$). Using the change of weight results of Chapter~\ref{IV}, we reduce the proof in general to a special case where, in particular, $\Delta_M$ is orthogonal to $\Delta-\Delta_M$ (\ref{V.7}). Finally, we use a filtration argument (\ref{V.8} to \ref{V.11}).

\subsection{}\label{V.3}

\noindent \textbf{Proposition} $R_G\otimes_{\cz_G(V)} \ind_K^GV$ \textit{is free over $R_G$.}

\medskip

The proof in \ref{V.4} requires several lemmas. We  use again the Kottwitz   homomorphism $w_G$  and the map  $v_G$ (\ref{III.16}).

\medskip

\noindent \textbf{Lemma 1}  \textit{Let $z$, $z_1$, $z_2$ in $Z$. If $zz_1z_2\in Kz_1Kz_2K$, then $w_G(z)=0$.}

\medskip

\noindent \textbf{Proof} The Kottwitz   homomorphism $w_G$ is a homomorphism of $G$ into a commutative group; the result follows from $w_G(K)=0$. $\square$

\medskip

\noindent\textbf{Lemma 2}  \textit{Let $z_1\in Z^+$ normalizing $\psi$, and $f\in \Hh_G(V)$ with support in $Kz_1K$. Then $\cs_Z^G(f)\in \Hh_Z(V_{U^0})$ has support in $(Z\cap \Ker w_G)z_1$.}

\medskip

\noindent \textbf{Proof} That is immediate from (III.3.2), once we note that $w_G$ is trivial on~$U$. $\square$

\medskip

 \noindent\textbf{Lemma 3} \textit{Let $z_1\in Z^+$ normalizing $\psi$, and $z_2\in Z$. If $f\in \ikg$ has support in $Kz_2K$, then $\tau_{z_1}* f$ has support in $K(Z\cap \Ker w_G) z_1z_2K$.}
 
 \medskip
 
 \noindent \textbf{Proof} By definition $\tau_{z_1}$, as an element of $\Hh_Z(V_{U^0})$, has support $Z^0z_1$. From Lemma~2, $\tau_{z_1}$, as an element of $\Hh_G(V)$, has support in $K(Z\cap \Ker w_G)z_1K$.
 The result then follows from the convolution formula in $\Hh_G(V)$ and Lemma~1. $\square$

\medskip

\noindent \textbf{Lemma 4}  $Z_{\Delta}^\bot \cap \Ker w_G =Z^0$.

\medskip

\noindent \textbf{Proof} Let $z\in \Ker w_G$. Then $v_G(z)=0$. If moreover $z\in Z_{\Delta}^\bot$, then $v_Z(z)=0$ for the analogous map $v_Z$, cf.\ \cite[6.3 Remark 1]{HV1}; from \cite[6.2 Lemma]{HV1}, (ii) and (iii), it follows that $z\in Z^0$. Conversely $Z^0 \subset Z_{\Delta}^\bot \cap \Ker w_G$ is clear. $\square$

\subsection{}\label{V.4} We prove \ref{V.3} Proposition. We decompose $\ikg $ as $\oplus I(x)$, $x\in Z/(Z\cap\Ker w_G)$, where $I(x)$ consists of the functions in $\ikg$ with support in $Kx(Z\cap\Ker w_G)K$. For $z$ in $Z^+$ normalizing $\psi$, we have $\tau_z* I(x) \subset I(zx)$ by \ref{V.3} Lemma~3, with equality if $z\in Z_{\Delta}^\bot$ since then $\tau_z$ has inverse $\tau_{z^{-1}}$. For $x\in Z/(Z\cap\Ker w_G)$, let $I^+(x)$ be the sum of the subspaces $\tau_z * I(y)$ of $I(x)$, where $z\in Z^+ \cap Z_\psi$, $z\notin  Z_{\Delta}^\bot $, $y\in Z/(Z\cap\Ker w_G)$ and $zy=x$ in $Z/(Z\cap\Ker w_G)$. By definition $R_G\otimes_{\cz_G(V)} \ikg$ is the quotient of $\ikg$ obtained by killing all the subspaces $I^+(x)$; thus it appears as $\oplus_{x\in Z/(Z\cap\Ker w_G)}(I(x)/I^+(x))$. Let $z\in Z_{\Delta}^\bot \cap Z_\psi$; then $\tau_z*I(x)=I(zx)$, $\tau_z* I^+(x)=I^+(zx)$ for $x\in Z/(Z\cap\Ker w_G)$, hence the corresponding element in $R_G$, still written $\tau_z$, sends $I(x)/I^+(x)$ isomorphically onto $I(zx)/I^+(zx)$. As $Z_{\Delta}^\bot  \cap \Ker w_G=Z^0$ by \ref{V.3} Lemma~4, the image of $Z_{\Delta}^\bot  \cap Z_\psi$ in $Z/(Z\cap\Ker w_G)$ acts  by multiplication without fixed points on $Z/(Z\cap\Ker w_G)$; choosing a set of representatives $\Omega$ for the orbits, we deduce that $R_G\otimes_{\cz_G(V)} \ikg$ is isomorphic to the free $R_G$-module $R_G \otimes_C (\bigoplus\limits_{x\in\Omega} I(x)/I^+(x))$. $\square$

\medskip

For further use, we state a result proved in a similar manner.
\medskip

\noindent \textbf{Lemma}  \textit{Let} $z\in Z^+\cap Z_\psi$.

\textit{(i) If $v_G(z)\not=0$, $\tau_z-1$ acts injectively on $\ikg$; if moreover $z\in Z_{\Delta}^\bot$ then $\tau_z-1$ is not a divisor of $0$ in $R_G$.}

\textit{(ii) Let $T\in \cz_G(V)$; if $v_G(z)$ is linearly independent from $v_G(\Supp (T))$, then $(\tau_z-1)\ikg\cap T\ikg=(\tau_z-1)T\ikg$.}

\medskip

\noindent\textbf{Remark}       
The condition $v_G(z)=0$ is equivalent to $v_Z(z)\!\in\! \R \Delta^\vee \!\subset\! X_*(\mathbf S) \otimes \R$. 

\medskip

\noindent\textbf{Proof} (i) Let $f\in \ikg$, and write as above $f=\sum f_x$, $x\in Z/(Z \cap \Ker w_G)$, $f_x\in I(x)$. Then for $z\in Z^+\cap Z_\psi$, $\tau_z*f=\sum_x\,\tau_z*f_x$ with $\tau_z*f_x\in I(zx)$. The equality $\tau_z*f=f$ amounts to $\tau_z*f_x=f_{zx}$ for all $x\in Z/(Z\cap\Ker w_G)$. If $v_G(z)\not=0$ then the image of $z$ in $Z/(Z \cap \Ker w_G)$ has infinite order; since $f_x=0$ for all but a finite number of $x$'s, $\tau_z*f=f$ implies $f=0$, and $\tau_z-1$ acts injectively on $\ikg$; in particular, as $\cz_G(V)$ acts faithfully on $\ikg$, $\tau_z-1$ is not a divisor of $0$ in $\cz_G(V)$. If moreover $z\in Z_{\Delta}^\bot$ then $\tau_z-1$ is not a divisor of $0$ in the subalgebra of $\cz_G(V)$ which maps isomorphically onto $R_G$.  

(ii) Let $\Gamma$ be the subgroup of $Z$ generated by the elements $\xi$ with $v_G(\xi)$ in $v_G(\Supp T)$. For $y\in Z/\Gamma$, let $J(y)$ be the space of functions in $\ikg$ with support in $Ky\Gamma K$; then $TJ(y)\subset J(y)$ and for $z\in Z^+ \cap Z_\psi$, $\tau_z* J(y) \subset J(zy)$. Let $f$, $f'$ in $\ikg$ with $(\tau_z-1)f=f'$. We have $\ikg=\oplus_{y\in Z/\Gamma}J(y)$ and decomposing accordingly $f=\sum f_y$ and $f'=\sum f_y'$ we get $\tau_z* f_y=f_{zy}+f_{zy}'$ for $y\in Z/\Gamma$. Let $f'\in T \ikg$; then $f_y'\in T \ikg$ for all $y\in Z/\Gamma$ so if $f_{y}$ belongs to $T\ikg$, then so do $\tau_z* f_y$ and $f_{zy}$. The hypothesis on $z$ in (ii) implies that its image in $Z/\Gamma$ has infinite order, so $f_{z^{-r}y}$ is $0$ for large $r$. So we get, using descending induction on $r$, that $f_y$ does indeed belong to $T\ikg$. $\square$

\subsection{}\label{V.5} We now turn to the general case of \ref{V.2} Theorem. For each parabolic subgroup $P_1=M_1N_1$ of $G$ containing $P$, we let $V_{P_1}$ be the irreducible representation of $K$ with parameter $(\psi,\Delta_{P_1}\cap \Delta(V))$ -- for $P_1=G$ we have $V_G=V$; we choose a basis vector for $(V_{P_1})_{U^0}$.

For such a $P_1$ consider the sequence of canonical (injective) intertwiners:
$$
\ikg_{P_1}\! \rg\! \Ind_{P_1}^G \ind_{M_1^0}^{M_1} (V_{P_1})_{N_1^0}\! \rg\! \Ind_P^G \ind _{M^0}^M (V_{P_1})_{N^0}\! \rg \! \Ind_B^G \ind _{Z^0}^Z (V_{P_1})_{U^0}.
\leqno(\mathrm V.5.1)
$$
As $(V_{P_1})_{N_1^0}$ has the same parameter as $V_{N_1^0}$, there is a unique isomorphism  between them that is compatible with the choice of basis vectors in $(V_{P_1})_{U^0}$ and $V_{U^0}$; it induces an isomorphism of $(V_{P_1})_{N^0}$ onto $V_{N^0}$. Using those isomorphisms we identify the sequence (V.5.1) with
$$
\ikg_{P_1} \rg \Ind_{P_1}^G \ind_{M_1^0}^{M_1} V_{N_1^0} \rg \Ind_P^G \ind _{M^0}^M V_{N^0} \rg \Ind_B^G \ind _{Z^0}^Z \psi.
\leqno(\mathrm V.5.2)
$$
The sequence (V.5.1) is equivariant for the sequence of Hecke algebras
$$
\Hh_G(V_{P_1}) \rg \Hh_{M_1}((V_{P_1})_{N_1^0})\rg \Hh_M ((V_{P_1})_{N^0}) \rg \Hh_Z((V_{P_1})_{U^0})
\leqno(\mathrm V.5.3)
$$
given by the (injective) Satake homomorphisms. The choice of basis vectors gives an isomorphism $\Hh_Z((V_{P_1})_{U^0}) \simeq \Hh_Z(\psi)$, actually independent of that choice, and inside $\Hh_Z(\psi)$ the Hecke algebras in (V.5.3) do not depend on $P_1$; accordingly we write $\Hh_G$, $\Hh_{M_1}$, $\Hh_M$, $\Hh_Z$, and similarly for the centres. The sequence (V.5.2) is then equivariant for the sequence of algebras $\Hh_G \rg \Hh_{M_1}\rg \Hh_M \rg \Hh_Z$.

As in Chapter \ref{IV} we identify the spaces in (V.5.2) with their images in $\Ind_B^G\, \ind _{Z^0}^Z\psi$, and similarly $\Hh_G$, $\Hh_{M_1}$, $\Hh_M$ with their images in $\Hh_Z$.

\medskip

\noindent\textbf{Notation} For $P_1$ as above containing $P$, we let $\pi_{P_1}$ be the  $\cz_{M}[G]$-submodule  $\cz_M\otimes_{\cz_G} \ikg_{P_1}$ of $\Ind_P^G \ind _{M^0}^M V_{N^0}$ (which is then $\pi_P$).

\medskip

\noindent \textbf{Remark} By \ref{III.14} Theorem, $\pi_{P_1}$ is also $\cz_M \otimes_{\cz_{M_1}} \Ind_{P_1}^G \ind _{M_1^0}^{M_1} V_{N_1^0}$, which we also see as $\Ind_{P_1}^G (\cz_M\otimes_{\cz_{M_1}} \ind _{M_1^0}^{M_1}V_{N_1^0})$, cf.~\cite{HV2} Corollary 1.3.

For further use, let us recall some useful facts. Let $X$ be a locally profinite space with a countable basis. Then the functor $X \mapsto C_c^\infty(X,A)$ is exact on $\Z$-modules $A$, $C_c^\infty(X,\Z)$ is free and $C_c^\infty(X,\Z) \otimes A \rg C_c^\infty(X,A)$ is an isomorphism; if $A$ is a free module over some ring $R$, then so is $\Cci(X,A)$ and if $R\rg R'$ is a ring homomorphism, then $R'\otimes_R \Cci(X,A) \rg \Cci(X,R'\otimes_R A)$ is an isomorphism of $R'$-modules. If $Y$ is an open subset of $X$, we have an exact sequence $0\rg \Cci(Y,\Z) \rg \Cci(X,\Z) \rg \Cci(X-Y,\Z)\rg 0$ of free $\Z$-modules. We are particularly interested in the case $X=J\ba H$ where $H$ is a locally profinite second countable group, and $J$ a closed subgroup of $H$. If $A$ is a smooth $R[J]$-module for some ring $R$,  choosing a continuous section  of $H \rg J\ba H$ gives an isomorphism  of $R$-modules $\Cci(J\ba H,A) \simeq \ind_J^HA$, so we deduce that $\ind_J^H$ is an exact functor on smooth $R[J]$-modules, that $\ind_J^H A$ is free over $R$ if $A$ is, and that $R'\otimes_R \ind_J^H A\rg \ind_J^H(R'\otimes_R A)$ is an isomorphism of $R'[H]$-modules for any ring homomorphism $R\rg R'$.

\subsection{}\label{V.6} We gather consequences of the change of weight results of Chapter \ref{IV}.

\medskip

\noindent\textbf{Proposition} \textit{Let $P_1,P_2$ be  parabolic subgroups of $G$ containing $P$, with $\Delta_{P_2}=\Delta_{P_1} \sqcup \{\alpha\}$.}

\textit{(i) $\pi _{P_2} \subset \pi_{P_1}$ with equality if $\alpha\notin \Delta(V)$ or if $\psi$ is not trivial on $Z^0 \cap M_\alpha'$.}

\textit{(ii) If $\alpha\in \Delta(V)$ and $\psi$ is trivial on $Z^0\cap M_\alpha'$, then $(\tau_\alpha-1)\pi_{P_1}\subset \pi_{P_2}$ (with $\tau_\alpha$ as in} \ref{III.16} \textit{Notation). If moreover $\alpha$ is not orthogonal to $\Delta_M$, the inclusion $\pi_{P_2}\subset \pi_{P_1}$ induces an isomorphism $R_M\otimes_{\cz_M}\pi_{P_2}\stackrel{\sim}{\longrightarrow} R_M\otimes_{\cz_M}\pi_{P_1}$.}

\medskip

\noindent\textbf{Proof} First note that if $\alpha\notin \Delta(V)$ then $V_{P_1}$ and $V_{P_2}$ are isomorphic, so $\pi_{P_2}=\pi_{P_1}$ is immediate. Assume $\alpha\in \Delta(V)$. We apply \ref{IV.1} Theorem to $V_{P_2}$ (in lieu of $V$) and $V_{P_1}$ (in lieu of $V'$). Choose $z\in Z_\psi$ with $|\alpha|(z)<1$ and $|\beta|(z)=0$ for $\beta\in \Delta$, $\beta\not=\alpha$; thus $\tau_z$ is an invertible element of $\cz_M$. By \ref{IV.1} Theorem (i), we have the inclusion $\tau_z \ikg _{P_2}\subset \ikg _{P_1}$ of subspaces of $\Ind_B^G(\ind _{Z^0}^Z\psi)$. As $\tau_z$ is invertible in $\cz_M$, we get $\pi_{P_2}\subset\pi_{P_1}$. If $\psi$ is not trivial on $Z^0\cap M_\alpha'$ then \ref{IV.1} Theorem (ii) gives $\tau_z \ikg _{P_1} \subset \ikg _{P_2}$ hence $\pi_{P_2}=\pi_{P_1}$. If $\psi$ is trivial on $Z^0\cap M_\alpha'$, \ref{IV.1} Theorem (ii) gives $\tau_z(1-\tau_\alpha)\ikg _{P_1} \subset \ikg _{P_2}$ so $(\tau_\alpha-1) \pi_{P_1} \subset \pi_{P_2}$. Now $\tau_\alpha=\tau_{a_\alpha}$ for $a_\alpha \in Z_\psi$ with $\nu(a_\alpha)=r \alpha^\vee$ with some positive rational number $r$ (\ref{III.16} Proposition (i), \ref{IV.12} Example). If $\alpha$ is not orthogonal to $\Delta_M$, we have $|\beta|(a_\alpha)<1$ for some $\beta\in \Delta_M$;  but $\tau_\alpha$ is sent to $0$ in $R_M$. This implies the last assertion. $\square$

\medskip

\subsection{}\label{V.7} We deduce a reduction for the proof of \ref{V.2} Theorem. Let $\Delta_1=\Delta_M\cup \{\alpha\in \Delta(V)$, $\alpha \bot \Delta_M$, $\psi(Z^0 \cap M_\alpha')=1\}$ and let $P_1=M_1N_1$ be the corresponding parabolic subgroup of $G$. By \ref{V.6} Proposition, the inclusion $\pi_G\subset \pi_{P_1}$ induces an isomorphism $R_M \otimes_{\cz_M} \pi_G \simeq R_M \otimes_{\cz_M} \pi_{P_1}$. But $R_M \otimes_{\cz_M}  \pi_{P_1}$ is the same as $\Ind_{P_1}^G (R_M \otimes_{\cz_{M_1}}  \ind _{M_1^0}^{M_1} V_{N_1^0})$ (\ref{V.5} Remark); if the $R_M$-module inside the induction is free, then so is $R_M \otimes_{\cz_M}  \pi_{P_1}$ (\ref{V.5} Remark). As a consequence, it is enough to prove \ref{V.2} when $\Delta_1=\Delta$. 

\medskip

\noindent\textbf{Assumption} (until \ref{V.11}): $\Delta=\Delta_M \cup \Delta(V)$, $(\Delta-\Delta_M) \bot \Delta_M$ and $\psi(Z^0 \cap M_\alpha')=1$ for $\alpha\in \Delta-\Delta_M$.

\medskip

\noindent\textbf{Notation} We put $\sigma=\ind _{M^0}^M V_{N^0}$, so $\pi_P=\Ind_P^G \sigma$. We also put $W(M)=\{w\in W_0$, $w^{-1}(\Delta_M) \subset \Phi^+\}$.   The above assumption allows us to define $\tau_\alpha$  as in  \ref{III.16} Notation.

\medskip

 By \ref{V.3} Proposition, we know that  $R_M\otimes_{\cz_M} \sigma$ is free over $R_M$, and so is $R_M \otimes_{\cz_M}  \pi_P$ (\ref{V.5} Remark). We want to deduce the same for $R_M\otimes_{\cz_M}  \pi_G$. 
For that we filter $\pi_P$ according to the double cosets $PwB$ for $w\in W(M)$  (recall that $G$ is the disjoint union of the double cosets $PwB$, $w\in W(M)$).

We consider upper sets in $W(M)$, i.e.\ subsets $A$ such that $v\in A$, $v'\in W(M)$ and $v'\geq v$ (in the Bruhat order) imply $v'\in A$. For an upper set $A$, $PAB=\bigcup_{v\in A} PvB$ is open in $G$ and we let $\pi_{P,A}$ be the subspace of functions in $\pi_P=\Ind_P^G\sigma$ with support in $PAB$; it is a $\cz_M$-submodule of $\pi_P$.

Let $A$ be non-empty upper set in $W(M)$ and choose a minimal element $w$ in $A$. Put $A'=A-\{w\}$; then $A'$ is an upper set in $W(M)$ and we have the submodule $\pi_{P,A'}$ of $\pi_{P,A}$.

Let $\bar{A}$, $\bar{A}'$ be the (open) images of $PAB$, $PA'B$ in $P\ba G$. We have an exact sequence of free $\Z$-modules
$$
0 \lgr \Cci(\bar{A}',\Z) \lgr \Cci(\bar{A},\Z) \lgr \Cci (\bar{A}-\bar{A}',\Z) \lgr 0\quad \mathrm{(cf.\ \ref{V.5}\ Remark).}
$$
Choosing a continuous section of $G\rg P\ba G$, and noting that $\bar{A}-\bar{A}'$ is the image of $PwB$ in $P\ba G$, we get from \ref{V.5} Remark that evaluating functions on $PwB$ gives an isomorphism of $\pi_{P,A}/\pi_{P,A'}$ with the $\cz_M$-module of locally constant functions $f:PwB\rg\sigma$ with $f(pg)=pf(g)$ for $p\in P$, $g\in PwB$, and with compact support in $P\ba PwB$; equivalently evaluating on $wU$ gives an isomorphism with the compactly induced representation $\ind _{w^{-1}Pw\cap U}^U {^w}\sigma$.

\medskip

\noindent\textbf{Lemma} \textit{The inclusion $\pi_{P,A}\rg \pi_P$ induces an isomorphism of $R_M\otimes_{\cz_M}  \pi_{P,A}$ onto the subspace of $R_M\otimes_{\cz_M} \pi_P =\Ind_P^G(R_M\otimes_{\cz_M} \sigma)$ consisting of functions with support in $PAB$. The sequence}
$$
0\lgr R_M\otimes_{\cz_M}  \pi_{P,A'}\lgr R_M\otimes_{\cz_M}  \pi_{P,A}\lgr R_M\otimes_{\cz_M}  (\pi_{P,A}/\pi_{P,A'}) \lgr 0
$$
\textit{is exact, and all three terms are free over $R_M$.}

\medskip

\noindent\textbf{Proof} Choosing a continuous section of $G\rg P\ba G$, $\pi_{P,A}$ appears as $\Cci(\bar{A},\Z)\otimes\sigma$, $R_M\otimes_{\cz_M} \pi_{P,A}$ as $\Cci(\bar{A},\Z) \otimes (R_M \otimes_{\cz_M} \sigma)$, so the result follows from \ref{V.5} Remark via the exact sequence $0\rg \Cci(\bar{A}',\Z)\rg \Cci(\bar{A},\Z)\rg \Cci(\bar{A}-\bar{A}',\Z) \rg 0$. $\square$

\subsection{}\label{V.8} Let $A$, $w$, $A'$ be as in \ref{V.7}, and let $Q$ be a parabolic subgroup of $G$ containing $P$. Then $\pi_Q\subset \pi_P$ and we let $\pi_{Q,A}= \pi_{P,A}\cap \pi_Q$, similarly for $A'$, so we get an inclusion of $\cz_M$-modules
$$
\pi_{Q,A}/\pi_{Q,A'} \hookrightarrow \pi_{P,A}/\pi_{P,A'}.
$$

\noindent\textbf{Notation} Set $c_{Q,w}= \Pi_{\alpha\in\Delta_Q, w^{-1}(\alpha)<0} (\tau_\alpha-1) \in \cz_M$.

\medskip

\noindent\textbf{Remarks} 1) For $\alpha\in\Delta$, $w^{-1}(\alpha)<0$ is equivalent to $s_\alpha w<w$ and it implies $\alpha\notin \Delta_M$ since $w\in W(M)$. In particular for such an $\alpha$ we have  $v_M(a_\alpha)\not=0$ by \ref{V.4} Remark.

2) By \ref{V.4} Lemma (i) (applied to $M$) $c_{Q,w}$ acts injectively on  $\sigma$ hence on $\pi_{P,A}/\pi_{P,A'}$; moreover,  $c_{Q,w}$ does not divide $0$  in $R_M$.

\medskip

\noindent\textbf{Proposition} $\pi_{Q,A}/\pi_{Q,A'}= c_{Q,w}(\pi_{P,A}/\pi_{P,A'})$ inside $\pi_{P,A}/\pi_{P,A'}$.

\medskip

Before we give the proof, we derive consequences, in particular \ref{V.2} Theorem.

\medskip

\noindent\textbf{Corollary 1} $R_M\otimes_{\cz_M} (\pi_{Q,A}/\pi_{Q,A'}) \rg R_M\otimes_{\cz_M}  (\pi_{P,A}/\pi_{P,A'})$ \textit{is injective, and $R_M\otimes_{\cz_M}  (\pi_{Q,A}/\pi_{Q,A'})$ is free over $R_M$.}

\medskip

\noindent\textbf{Proof} By the proposition, multiplication by $c_{Q,w}$ induces maps 
$$
\pi_{P,A}/\pi_{P,A'}  \twoheadrightarrow   \pi_{Q,A}/\pi_{Q,A'} \hookrightarrow \pi_{P,A}/\pi_{P,A'}.
$$
Tensoring with $R_M$  over $\cz_M$ gives 
$$
R_M\otimes_{\cz_M}  (\pi_{P,A}/\pi_{P,A'}) \twoheadrightarrow R_M\otimes_{\cz_M}  (\pi_{Q,A}/\pi_{Q,A'}) \lgr R_M\otimes_{\cz_M}  (\pi_{P,A}/\pi_{P,A'})
$$
whose composite is multiplication by $c_{Q,w}$. By the above remark 2) $c_{Q,w}$ does not divide $0$ in $R_M$; since $R_M\otimes_{\cz_M}  (\pi_{P,A}/\pi_{P,A'})$ is free over $R_M$ by \ref{V.7} Lemma, multiplication by $c_{Q,w}$ is injective on it so we get an isomorphism $R_M\otimes_{\cz_M}  (\pi_{P,A}/\pi_{P,A'}) \simeq R_M\otimes_{\cz_M}  (\pi_{Q,A}/\pi_{Q,A'})$, thus proving Corollary~1. $\square$

\medskip
 
\noindent\textbf{Corollary 2} \textit{$R_M\otimes_{\cz_M}  \pi_{Q,A}\rg R_M\otimes_{\cz_M}  \pi_{P,A}$ is injective (in particular, for $A=W(M)$, $R_M\otimes_{\cz_M}  \pi_Q \rg R_M\otimes_{\cz_M}  \pi_P$ is injective).}

\medskip

\noindent\textbf{Proof} By induction on $\# A$, $R_M \otimes_{\cz_M}  \pi_{Q,A'} \rg R_M\otimes_{\cz_M} \pi_{P,A'}$ is injective. By \ref{V.7} Lemma, $R_M \otimes_{\cz_M}  \pi_{P,A'} \rg R_M\otimes_{\cz_M}  \pi_{P,A}$ is injective and by Corollary 1, $R_M\otimes_{\cz_M}  (\pi_{Q,A}/\pi_{Q,A'}) \rg R_M\otimes_{\cz_M}  (\pi_{P,A}/\pi_{P,A'})$ is injective too. The result follows from the snake lemma applied to the commutative diagram (with exact rows)
$$
\begin{array}{cccccccccc}
 &  & R_M\otimes_{\cz_M}  \pi_{Q,A'} & \rg & R_M\otimes_{\cz_M}  \pi_{Q,A} & \rg & R_M\otimes_{\cz_M} (\pi_{Q,A}/\pi_{Q,A'}) & \rg & 0 &  \\ 
 &  & \downarrow &  & \downarrow &  & \downarrow &  &  &  \\ 
0 & \rg & R_M\otimes_{\cz_M}  \pi_{P,A'} & \rg & R_M\otimes_{\cz_M}  \pi_{P,A} & \rg & R_M\otimes_{\cz_M}  (\pi_{P,A}/\pi_{P,A'}) & \rg & 0 & \square
\end{array} 
$$

\noindent\textbf{Corollary 3}\textit{ $R_M\otimes_{\cz_M} \pi_{Q,A'} \rg R_M\otimes_{\cz_M} \pi_{Q,A}$ and $R_M\otimes_{\cz_M} \pi_{Q,A}\rg R_M\otimes_{\cz_M}\pi_Q$ are injective, and $(R_M\otimes_{\cz_M} \pi_{Q,A})/(R_M\otimes_{\cz_M} \pi_{Q,A'})\rg R_M\otimes_{\cz_M}  (\pi_{Q,A}/\pi_{Q,A'})$ is an isomorphism.}

\medskip

\noindent\textbf{Proof} In the left hand square of the previous diagram, the two vertical maps and the bottom horizontal one are injective, hence so is the top horizontal one, giving the first assertion, which immediately implies the last one. The second one follows from the first by descending induction on $\# A$. $\square$

\medskip
Now \ref{V.2} Theorem follows from the corollaries.
Indeed, by Corollary 1 and Corollary 3, $R_M\otimes_{\cz_M}\pi_Q$ is a successive extension of free modules.
Therefore $R_M\otimes_{\cz_M}\pi_Q$ is free.

\subsection{}\label{V.9} The proof of \ref{V.8} Proposition will involve an induction argument on $\dim G$. For this, a further corollary is necessary.

\medskip

\noindent\textbf{Corollary 4} \textit{Let $z\in{Z ^{+M}}\cap Z_\psi$, and assume $v_M(z) \not=0$. Then $\pi_Q \cap (\tau_z-1)\pi_P=(\tau_z-1)\pi_Q$}.

\medskip

 The proof  is given after  a  lemma.  Let $A$, $w$, $A'$ be as in \ref{V.7}, and use the notation $\pi_{P,A}$, $\pi_{Q,A}$ of \ref{V.7},~\ref{V.8}.

\medskip

\noindent\textbf{Lemma} $(\tau_z-1)\pi_{P,A}=(\tau_z-1)\pi_P \cap \pi_{P,A}$.

\medskip

\noindent\textbf{Proof} By descending induction on $\#A$, the case $A=W(M)$ being trivial. By \ref{V.4} Lemma (i), $\tau_z-1$ acts injectively on $\sigma$, hence also on $\pi_{P,A}/\pi_{P,A'}$ which is a direct sum of copies of $\sigma$ (\ref{V.5} Remark). By the snake lemma $\pi_{P,A'}/(\tau_z-1)\pi_{P,A'}$ injects into $\pi_{P,A}/(\tau_z-1)\pi_{P,A}$ i.e.\ $(\tau_z-1)\pi_{P,A}\cap \pi_{P,A'}=(\tau_z-1)\pi_{P,A'}$. The assertion $(\tau_z-1)\pi_{P,A}=(\tau_z-1)\pi_P \cap \pi_{P,A}$ then implies the similar assertion for $A'$. $\square$

\medskip \noindent\textbf{Proof of Corollary 4} 
Applying \ref{V.4} Lemma (ii) to $T= c_{Q,w}\in \cz_M$ whose support is in $\Ker v_M$ we get
$$
(\tau_z-1) \sigma \cap c_{Q,w}\sigma = (\tau_z-1) c_{Q,w}\sigma,
$$
hence
$$
(\tau_z-1) (\pi_{P,A}/\pi_{P,A'}) \cap c_{Q,w}(\pi_{P,A}/\pi_{P,A'})=(\tau_z-1) c_{Q,w}(\pi_{P,A}/\pi_{P,A'}).
$$
But by \ref{V.8} Proposition $c_{Q,w}(\pi_{P,A}/\pi_{P,A'})=\pi_{Q,A}/\pi_{Q,A'}$, so we obtain $(\tau_z-1)\pi_{P,A}\cap \pi_{Q,A} \subset (\tau_z-1)\pi_{Q,A}+\pi_{P,A'}$. By the lemma we get $(\tau_z-1) \pi_P \cap \pi_{Q,A} \subset [(\tau_z-1)\pi_{Q} \cap \pi_{P,A}]+\pi_{P,A'}$. As $(\tau_z-1)\pi_P \cap \pi_Q$ contains $(\tau_z-1)\pi_Q$, both give the same contribution to $\pi_{P,A}/\pi_{P,A'}$. Their equality now follows by induction on $\#A$. $\square$

\medskip

\subsection{}\label{V.10} We now proceed to the proof of \ref{V.8} Proposition, keeping its notation. We first deal with the (easier)  statement   that $\pi_{Q,A}/\pi_{Q,A'}$ contains $ c_{Q,w}(\pi_{P,A}/\pi_{P,A'})$.

\medskip

\noindent \textbf{Notation} Let $\Delta_w=\{\alpha\in \Delta,\, w^{-1}(\alpha)>0\}$ and let $P_w=M_wN_w$ be the corresponding parabolic subgroup of $G$; it contains $P$, and   $w$ is  in $W(M_w)$.

\medskip

\noindent\textbf{Lemma} \textit{Let $A$, $w$, $A'$ be as in \upshape{\ref{V.7}}. Then $\pi_{P_w,A} \rg \pi_{P,A}/\pi_{P,A'}$ is surjective.}

\medskip

Assume that lemma for a moment. Since $\pi_{Q\cap P_w,A}$ contains $\pi_{P_w,A}$, the map $\pi_{Q\cap P_w,A}\rg  \pi_{P,A}/\pi_{P,A'}$ is surjective as well. But by \ref{V.6} Proposition $\pi_{Q}$ contains $c_{Q,w}\, \pi_{Q\cap P_w}$, so the image of $\pi_{Q,A}$ in $\pi_{P,A}/\pi_{P,A'}$ contains $c_{Q,w}(\pi_{P,A}/\pi_{P,A'})$, i.e.\  the quotient $\pi_{Q,A}/\pi_{Q,A'}$ contains $ c_{Q,w}(\pi_{P,A}/\pi_{P,A'}) $.

\medskip

\noindent\textbf{Proof} Let $A_{\ge w}=\{v\in W(M),\, v\ge w\}$ and $A_{>w}=\{v\in W(M),v>w\}$. We use the abbreviations $ \pi_{P,\ge w}=  \pi_{P,A_{\ge w}}, \  \pi_{P,> w}=  \pi_{P,A_{> w}}$. 
Then $\pi_{P,A}\supset \pi_{P,\ge w}$ and $\pi_{P,A'} \supset \pi_{P,>w}$; moreover $\pi_{P,A'} \cap \pi_{P,\ge w}=\pi_{P,>w}$, so $\pi_{P,\ge w}/\pi_{P,>w}$ injects into $\pi_{P,A}/\pi_{P,A'}$. But evaluation on $PwB$ identifies both quotients with the same space of functions, so the injection is an isomorphism. Hence it is enough to prove the lemma for $A=A_{\ge w}$.

\medskip

\noindent\textbf{Sublemma} \textit{(i) $w^{-1} Pw\cap U= w^{-1}U w \cap U = w^{-1}P_w w\cap U$.}

\textit{(ii)} $PA_{\ge w}B =   \sqcup_{v\in W(M_w),v\ge w} P_w vB$.  

\medskip

\noindent\textbf{Proof} (i) The first equality comes from $w\in W(M)$, the second one from $w\in W(M_w)$.

(ii) By \cite[Lemma 4.20]{Abe}, $w\in W(M)$ implies $W_MA_{\ge w}=\{v\in W_0,v\ge w\}$ and similarly $w\in W(M_w)$ implies $W_{M_w}\{v\in W(M_w),v\ge w\}=\{v\in W_0,v\ge w\}$. The result follows on taking $B$-double cosets. $\square$

\medskip

To prove the lemma (for $A=A_{\ge w}$) we need to consider closely the inclusion $\pi_{P_w} \hookrightarrow \pi_P$. Both are parabolically induced from $P_w$, and the inclusion comes from the injective map $\Phi:\cz_M\otimes_{\cz_{M_w}} \ind _{M_w^0}^{M_w} V_{N_w^0} \rg \Ind_{P\cap M_w}^{M_w}\sigma$ obtained from the canonical intertwiner (III.13.1), so $\pi_{P_w}$ is simply the subspace $\Ind_{P_w}^G(\Ima \Phi)$ of $\pi_P = \Ind_{P_w}^G(\Ind_{P\cap M_w}^{M_w}\sigma)$. Seeing $\pi_P$ as induced from $P_w$, we let $\pi_{P,\ge w}'$ be the subspace of functions with support in $\bigcup_{v\in W(M_w),v\ge w} P_w vB$, and similarly $\pi_{P,> w}'$. An element $f$ of $\pi_P=\Ind_P^G\sigma$ is seen as the function $f'$ in $\Ind_{P_w}^G(\Ind_{P\cap M_w}^{M_w}\sigma)$ given by $f'(g): m \mapsto f(mg)$ for $g\in G$, $m\in M_w$. Hence by (ii) of the sublemma $\pi_{P,\ge w}=\pi_{P,\ge w}'$ and $\pi_{P,>w}\supset \pi_{P,>w}'$. By (i) of the sublemma (and \ref{V.5} Remark), choosing a continuous section of $U\rg w^{-1}Uw\cap U\ba U$ gives a $\cz_M$-linear isomorphism $\iota$ of $\pi_{P,\ge w}'/\pi_{P,>w}'$ with $\Cci(w^{-1} U w\cap U\ba U,\Z)\otimes \Ind_{P\cap M_w}^{M_w}\sigma$, a similar isomorphism of $\pi_{P,\ge w}/\pi_{P,>w}$ with $\Cci(w^{-1} Uw\cap U\ba U,\Z)\otimes \sigma$, and the quotient map $\pi_{P,\ge w}'/\pi_{P,> w}' \twoheadrightarrow \pi_{P,\ge w}/\pi_{P,> w}$ corresponds to evaluation at $1:\Ind_{P\cap M_w}^{M_w}\sigma \rg \sigma$. But $(\pi_{P_w}\cap \pi_{P,\ge w}')/(\pi_{P_w}\cap \pi_{P,\ge w}')$ is sent by $\iota$ to $\Cci(w^{-1} U w \cap U\ba U,\Z)\otimes \Ima\Phi$ so to get the surjectivity of $\pi_{P_w} \cap \pi_{P,\ge w}' \rg \pi_{P,\ge w}/\pi_{P,> w}$ it suffices to see that evaluation at $1: \Ima \Phi \rg \sigma$ is surjective. But for $x\in V_{N_w^0}$ the function in $\ind _{M_w^0}^{M_w}V_{N_w^0}$ with support $M_w^0$ and value $x$ at 1, is sent in $\Ind_{P\cap M_w}^{M_w}\sigma$ to a function with value at 1 the function in $\sigma$ with support $M^0$ and value at 1 the projection of $x$ in $V_{N^0}$; as those last functions, for varying $x$, generate $\sigma$ as a representation of $M$, and $\Ima\Phi\rg \sigma$ is $M$-equivariant, it is surjective. $\square$

\subsection{}\label{V.11} We turn to the  inclusion   $\pi_{Q,A}/\pi_{Q,A'} \subset  c_{Q,w}(\pi_{P,A}/\pi_{P,A'})$  in \ref{V.8} Proposition. We need auxiliary lemmas, where $\alpha\in \Delta-\Delta_M$ is fixed; we let $P^\alpha=M^\alpha N^\alpha$ be the parabolic subgroup corresponding to $\Delta_M\cup \{\alpha\}$ and we put $\bar{\sigma}
=\sigma/(\tau_\alpha-1)\sigma$. Note that   Hypothesis  (H) of \ref{III.15} holds with the map $\varphi: V_{N^0}\to \sigma \to \bar{\sigma}$. We also note that $\varphi \tau_\alpha = \tau_\alpha \varphi=\varphi$.

\medskip

\noindent\textbf{Lemma 1} \textit{$\bar{\sigma}$ extends to $P^\alpha$, trivially on $N$.}

\medskip

\noindent\textbf{Proof} By \ref{II.7} it suffices to prove that $\bar{\sigma}$ is trivial on $Z\cap M_\alpha'$. Since   $\alpha$ is orthogonal to $\Delta_M$  and  $\psi$ is trivial on $Z^0\cap M_\alpha'$ by assumption, that comes from the fact that $\tau_\alpha$ acts trivially on $\bar{\sigma}$ (\ref{III.17}). $\square$

\medskip

We write $^e\bar{\sigma}$ for the extension of $\Bs$ to $P^\alpha$. Inside of $\pi_P/(\tau_\alpha-1)\pi_P\simeq \Ind_P^G\Bs$ we have the subspace $\Ind_{P^\alpha}^G \,^e\Bs$, cf.\ \ref{III.22} Lemma 2.

\medskip
\noindent\textbf{Lemma 2} \textit{The image of $\pi_{P^\alpha}\rg \pi_P\rg \pi_P/(\tau_\alpha-1)\tau_P$ is contained in $\Ind_{P^\alpha}^G  {}^e\Bs$.}

\medskip
\noindent\textbf{Proof} Since $\pi_{P^\alpha} \rg \Ind_P^G \Bs$ is $\cz_M[G]$-equivariant and $\pi_{P^\alpha}$ is generated as a $\cz_M[G]$-module by $V_{P^\alpha}$ it is enough to prove that the inclusion of $\Hom_K(V_{P^\alpha}, \Ind_{P^\alpha}^G{} ^e\Bs)$ into $\Hom_K(V_{P^\alpha},\Ind_P^G\Bs)$ is an isomorphism. By Frobenius reciprocity, this means that
$$
\Hom_{M^{\alpha 0}}((V_{P^\alpha})_{N^{\alpha 0}},  {}^e \Bs) \hookrightarrow \Hom_{M^{\alpha0}} ((V_{P^\alpha})_{N^{\alpha 0}}, \Ind_{P\cap M^\alpha }^{M^\alpha }\Bs)
$$
is an isomorphism. The quotient of $\Ind_{P\cap M^\alpha}^{M^\alpha}\Bs$ by $^e\Bs$ is the  representation $^e\Bs \otimes \St_{P\cap M^\alpha }^{M^\alpha}$ and it is enough to show that $(V_{P^\alpha})_{N^{\alpha 0}}$ is not a weight of that representation. But the parameter for $(V_{P^\alpha})_{N^{\alpha 0}}$ is $(\psi,(\Delta_M \cup \{\alpha\} )\cap \Delta(V))$ and $\alpha\in \Delta(V)$ whereas by \ref{III.18} the weights of $^e\Bs \otimes \St_{P\cap M^\alpha }^{M^\alpha} = I(P\cap M^\alpha,  \Bs,P\cap M^\alpha)  $ have parameters $(\psi',I)$ where $\alpha\notin I$. $\square$

\medskip

\noindent\textbf{Lemma 3} \textit{Let $P_1 =M_1N_1$ and $P_2=M_2N_2$ be parabolic subgroups of $G$ containing $P$, and assume $\Delta_{P_2}=\Delta_{P_1}\sqcup \{\alpha\}$. Let  $A$, $w$, $A'$ be as in \upshape{\ref{V.7}}, and assume that $w^{-1}(\beta)<0$ for all   $\beta\in \Delta_{P_2}-\Delta_M$. Then}
$$
\pi_{P_2,A} \subset (\tau_\alpha-1) \pi_{P_1,A}+ \pi_{P,A'}.
$$

\noindent\textbf{Proof} Let $f\in \pi_{P_2,A}$ and let $\bar{f}$ be its image in $\Ind_P^G\Bs$. As $\pi_{P_2}\subset \pi_{P^\alpha}$, we get $\bar{f} \in \Ind_{P^\alpha}^G {}^e\Bs$ by Lemma 2. If $\bar{f}$ does not vanish on $PwB$, its support, being $P^\alpha$ invariant, contains $P s_\alpha wB$. But $w^{-1}(\alpha)<0$ means $s_\alpha w<w$ and $w$ being minimal in $A$, that contradicts $f\in \pi_{P,A}$. Hence $\bar{f}$ vanishes on $PwB$ and there exist $f_1\in \pi_P$, $f_2\in \pi_{P,A'}$ with $f=(\tau_\alpha-1) f_1+f_2$. The point is to prove that we can take $f_1$ in $\pi_{P_1,A}$. View $\pi_{P_1}$ as $\Ind_{P_1}^G \sigma_1$ with $\sigma_1=\cz_M \otimes_{\cz_{M_1}} \ind _{M_1^0}^{M_1} V_{N_1^0}$ and $\pi_P$ as $\Ind_{P_1}^G \Ind_{P\cap M_1}^{M_1}\sigma$, the inclusion $\pi_{P_1}\hookrightarrow \pi_P$ being induced by the natural intertwiner $\ind_{M_1^0}
^{M_1}V_{N_1^0} \rg \Ind_{P\cap M_1}^{M_1}\sigma$.

\medskip

\noindent\textbf{Sublemma} \textit{For $v\in A'$ we have $P_1wB\cap PvB=\emptyset$.}

\medskip

\noindent\textbf{Proof} Indeed, if $P_1wB \cap PvB\not=\emptyset$ there exists $v'$ in $W_{M_1}$ with $v'w=v$. Since $\Delta-\Delta_M$ is orthogonal to $\Delta_M$, $W_{M_1}$ is the direct product of $W_M$ and the subgroup $W_1$ generated by the $s_\beta$ for $\beta\in \Delta_{P_1}-\Delta_M$. For such a $\beta$ we have $s_\beta w<w$ and it follows (using \cite{Deo} as in \ref{IV.9} Lemma 2), by induction on length, that $v_1w\leq w$ for any   $v_1\in W_1$. Writing $v'$ as $v_2^{-1}v_1$ with $v_1\in W_1$ and $v_2\in W_M$ we get $v_1w=v_2v$. But $v_2 v\ge v$ and $v_1w\le w$ so $v\le w$ contrary to the assumption $v\in A'$. $\square$ 

\medskip

Let us pursue the proof of Lemma 3.

Since $f_2\in \pi_{P,A'}$, it follows from the sublemma that, seen as an  element of $\Ind_P^G\sigma$, it vanishes on $P_1wB$; but then, seen as an element of $\Ind_{P_1}^G \Ind_{P\cap M_1}^{M_1}\sigma$, it also vanishes on $P_1wB$. So for any $x\in P_1wB$, $f(x)=(\tau_\alpha-1)f_1(x)$ in $\Ind_{P\cap M_1}^{M_1}\sigma$. Now $\dim P_1<\dim G$ so \ref{V.8} Proposition and all its corollaries are true for $M_1$. As $v_{M_1}(a_\alpha)\not=0$ we conclude from \ref{V.9} Corollary 4 that there exists $y\in \sigma_1$ with $(\tau_\alpha-1)f_1(x)=(\tau_\alpha-1)y$. But $\tau_\alpha-1$ does not kill any element of $\Ind_{P\cap M_1}^{M_1}\sigma$, by \ref{V.4} Lemma (i), so $f_1(x)=y$ belongs to $\sigma_1$. We can choose $f_1'$ in  $ \Ind_{P_1}^G \sigma_1\cap  \pi_{P_1,A}$  with the same restriction as $f_1$ on $P_1wB$ (use \ref{V.5} Remark). Put $f_2'=f-(\tau_\alpha-1)f_1'=(\tau_\alpha-1)(f_1-f_1')+f_2$. Then  $f_1-f_1'$ vanishes on $P_1 w B$, $f_2$ vanishes on $PwB$ so $f_2'$ belongs to $\pi_{P,A'}$ and $f=(\tau_\alpha-1) f_1'+f_2'$ belongs to $(\tau_\alpha-1) \pi_{P_1,A}+\pi_{P,A'}$. $\square$

\medskip

We now finish the proof of \ref{V.8} Proposition. Let $R$ be the parabolic subgroup between $P$ and $Q$ with   $ \Delta_R-\Delta_M=\{\alpha\in \Delta_Q, w^{-1}(\alpha)<0\}$. Applying Lemma 3 successively we get $\pi_{R,A} \subset c_{R,w}\pi_{P,A}+\pi_{P,A'}$, hence the result since $\pi_{Q,A}\subset \pi_{R,A}$ and $c_{R,w}=c_{Q,w}$. $\square$

\medskip

We can get more out of that:

\medskip

\noindent\textbf{Lemma 4} \textit{Let $A$, $w$, $A'$ be as in \upshape{\ref{V.7}}. Then $\pi_{Q,A}\subset c_{Q,w} \pi_{P_w,A}+\pi_{Q,A'}$.} 

\medskip

\noindent\textbf{Proof} \ref{V.10} Lemma gives $\pi_{P,A}\subset \pi_{P_w,A}+\pi_{P,A'}$ so from \ref{V.8} Proposition we get $\pi_{Q,A} \subset c_{Q,w}\pi_{P_w,A}+\pi_{P,A'}$. But $\pi_{P_w}\subset \pi_{Q\cap P_w}$ and $c_{Q,w}\pi_{Q\cap P_w}\subset \pi_Q$ by \ref{V.6} Proposition so $c_{Q,w}\pi_{P_w,A}\subset \pi_{Q,A}$ and the result follows. $\square$

\paragraph{{\large B) Filtration theorem for $\chi\otimes_{\cz_G} \ind_K^G V$}}
\addcontentsline{toc}{section}{\ \ B) Filtration theorem for \texorpdfstring{$\chi\otimes_{\cz_G} \ind_K^G V$}{chi otimes ind\_K\^{}G V}}

\subsection{}\label{V.12} We now turn to the filtration theorem (\ref{I.6}). For that, as before, an irreducible representation $V$ of $K$ is fixed, with parameter $(\psi,\Delta(V))$, but we also fix a character $\chi$ of $\cz_G=\cz_G(V)$. We let $P=MN$ be the  parabolic subgroup with $\Delta_P=\Delta_0(\chi)$, so $P$ is the smallest parabolic subgroup of $G$ containing $B$ such that $\chi$ extends to a character -- still written $\chi$ -- of $\cz_M=\cz_M(V_{N^0})$, and that character further factors through $\cz_M \rg R_M$.

\medskip

\noindent\textbf{Notation} For a $\cz_M$-module $W$, we put $W^\chi=\chi \otimes_{\cz_M}W$.

\medskip

Recall that for each parabolic subgroup $Q$ of $G$ containing $P$, $V_Q$ denotes the  irreducible representation   of $K$ of parameter $(\psi, \Delta_Q \cap \Delta(V))$; we make the same identifications as in \ref{V.5}. In particular we get a $\cz_M[G]$-submodule $\pi_Q$ of $\pi_P=\Ind_P^G \sigma$ -- we keep writing $\sigma = \ind _{M^0}^M V_{N^0}$. Our main interest is in $\pi_G^\chi$, but its analysis goes through the $\pi_Q^\chi$, in particular $\pi_P^\chi$.

As $\sigma^\chi$ satisfies property (H) of \ref{III.15}, the maximal parabolic subgroup   of $G$ to which $\sigma^\chi$ extends, trivially on $N$, has associated set of roots $\Delta_M \sqcup \Theta_{\mathrm{max}}$ where $\Theta_{\mathrm{max}}$ is the set of $\alpha\in \Delta-\Delta_M$, orthogonal to $\Delta_M$ and such that $\psi(Z^0 \cap M_\alpha')=1$ and $\chi(\tau_\alpha)=1$ (\ref{III.17} Corollary).

\medskip

\noindent\textbf{Notation} We let $\Theta=\Theta_{\mathrm{max}}\cap \Delta(V)$, $P_e=P_{\Delta_M\sqcup \Theta}$ and write $^e\sigma^\chi$ for the extension of $\sigma^\chi$ to $P_e$, trivial on $N$. (Note that  \ref{III.22} Lemma 2  gives an identification of $\pi_P^\chi$ with $\Ind_{P_e}^G(^e\sigma^\chi \otimes \Ind_P^{P_e}1)$.)

\medskip

\noindent\textbf{Lemma} \textit{The inclusion $\pi_G\rg \pi_{P_e}$ induces an isomorphism $\pi_G^\chi \rg \pi_{P_e}^\chi$}

\medskip \noindent\textbf{Proof} It suffices to show that for $P_e \subset P_1\subset P_2$ with $\Delta_{P_2}=\Delta_{P_1}\sqcup \{\alpha\}$,  the natural map $\pi_{P_2}^\chi  \to \pi_{P_1}^\chi$ is an isomorphism. If 
$\alpha\notin \Delta(V)$ or if $\psi$ is not trivial on $Z^0 \cap M_\alpha'$ 
or  $\alpha$  not orthogonal to $\Delta_M$, then by \ref{V.6} we even have   an isomorphism $R_M\otimes_{\cz_M}\pi_{P_2}\stackrel{\sim}{\longrightarrow} R_M\otimes_{\cz_M}\pi_{P_1}$. Otherwise, $\chi (\tau_\alpha)\neq 1$ and since 
$(\tau_{\alpha}-1)\pi_{P_1}\subset \pi_{P_2}\subset \pi_{P_1}$ by \ref{V.6}, we have an isomorphism $\pi_{P_2}^\chi  \stackrel{\sim}{\longrightarrow} \pi_{P_1}^\chi$.
 $\square$

\subsection{}\label{V.13} 

\noindent\textbf{Notation}  Let $\Cd$ be the set of parabolic subgroups of $G$ between $P$ and $P_e$.

\medskip

$\bullet$ For $Q$, $Q_1$ in $\Cd$, $Q\supset Q_1$, put $c_{Q,Q_1}=\prod\limits_{\alpha\in \Delta_Q-\Delta_{Q_1}}(\tau_\alpha-1)$  (then $c_{Q,Q_1}\pi_{Q_1}\subset \pi_Q$ by \ref{V.6} Proposition).

$\bullet$ For $Q\in \Cd$,   let $\tau_Q$ be the image of $\pi_Q \stackrel{c_{P_e,Q}}{\longrightarrow} \pi_{P_e} \rg \pi_{P_e}^\chi(=\pi_G^\chi)$, and let $\rho_Q$ be the image of   $\pi_Q\hookrightarrow  \pi_P \to  \pi_P^\chi$.

\medskip
 
$\bullet$ For $Q\in \Cd$, let $Q^c\in \Cd$ be the parabolic subgroup such that  $\Delta_{Q^c}-\Delta_M=\Delta_{P_e}-\Delta_Q$, and let $\Phi_Q$, $\Psi_Q$ be the $G$-equivariant maps
$$
\begin{array}{llllllll}
\Phi _Q&: & \tau _Q& \hookrightarrow & \pi_{P_e}^\chi & \lgr & \pi_{Q^c}^\chi,&   \\ 
\Psi _Q&: & \rho_Q & \hookrightarrow & \pi_P^\chi & \stackrel{c_{Q^c,P}}{\longrightarrow} & \pi_{Q^c}^\chi. & 
\end{array} 
$$
Here, the last map is obtained from $ \pi_P  \stackrel{c_{Q^c,P}}{\longrightarrow}  \pi_{Q^c}$ by  tensoring by $\chi$.

$\bullet$ Let $I_Q$ the submodule $\Ind _{P_e}^G(^e\sigma^\chi \otimes \Ind_Q^{P_e}1)$ of $\pi_P^\chi$.
  In particular, $I_P=\rho_P=\pi_P^\chi$. Note also that  $\tau_{P_e}=\pi_{P_e}^\chi$.

\medskip
\noindent\textbf{Remark 1}   The maps $\pi_Q \stackrel{c_{P_e,Q}}{\longrightarrow} \pi_{P_e}  \hookrightarrow \pi_{Q^c}$ and $\pi_Q \hookrightarrow \pi_P  \stackrel{c_{Q^c,P}}{\longrightarrow}  \pi_{Q^c}$ are equal because $c_{Q^c,P}=c_{P_e,Q}$. Therefore $\Ima \Phi_Q=\Ima \Psi_Q$.

\medskip
\noindent\textbf{Remark 2}  For $Q$, $Q_1$ in $\Cd$, $Q\supset Q_1$, we have  $\tau_{Q_1}\subset \tau_Q$ and $\rho_{Q_1}\supset \rho_Q$.

\medskip

Our second main result in this chapter is:

\medskip

\noindent\textbf{Theorem} \textit{Let $Q\in \Cd$.}

\textit{(i)} $\rho_Q=I_Q$.

\textit{(ii)}   $\Ker \, \Psi_Q=\sum\limits_{Q_1\in \Cd, Q_1\varsupsetneq Q}\rho_{Q_1}$.

\textit{(iii)} 
$
\Ker \Phi_Q = \sum\limits_{Q_1\in \Cd, Q_1\varsubsetneq Q} \tau_{Q_1}.
$

\textit{(iv) Let $\cp \subset \Cd$; then $\tau_Q\cap \sum\limits_{Q_1\in \cp} \tau_{Q_1}=\sum\limits_{Q_1\in \cp}\tau_{Q\cap Q_1}$.}
 
 \medskip
 
 It implies \ref{I.6} Theorem 6:
 
 \medskip
 
\noindent\textbf{Corollary 1} For $Q\in \Cd$, $\tau_Q/\sum\limits_{Q_1\in \Cd, Q_1\varsubsetneq Q}\tau_{Q_1}$ is isomorphic to $I_e(P,\sigma^\chi,Q)$.

\medskip

\noindent\textbf{Proof} By Remark 1 we have that  $\tau_Q/\Ker \Phi_Q$ is isomorphic to $\rho_Q/\Ker \, \Psi_Q$. But $\rho_Q=I_Q$ by (i) so  we get by (ii) and (iii) a $G$-isomorphism between $\tau_Q/\sum\limits_{Q_1\in \Cd, Q_1\varsubsetneq Q}\tau_{Q_1}$ and $I_Q/\sum\limits_{Q_1\in \Cd,Q_1\varsupsetneq Q}I_{Q_1}$ which is $I_e(P,\sigma^\chi,Q)$. $\square$

\medskip

\noindent\textbf{Corollary 2} \textit{Enumerate the parabolic subgroups in $\Cd$ as $P=Q_1,\ldots, Q_r=P_e$, so that $i\le j$ if $Q_i\subset Q_j$. For $i=0,\ldots,r$, put $I_i= \sum\limits_{1\le j\le i} \tau_{Q_j}$. Then for $i=1,\ldots,r$, $I_i/I_{i-1}\simeq I_e(P,\sigma^\chi,Q_i)$.
}
\medskip

\noindent\textbf{Proof} For $i=1,\ldots,r$ $I_i/I_{i-1}=\tau_{Q_i}/(\tau_{Q_i}\cap \sum\limits_{1\le j<i} \tau_{Q_j})$ is also $\tau_{Q_i}/\sum\limits_{1\le j<i}\tau_{Q_i\cap Q_j}$ by (iv). The assertion follows from Corollary 1. $\square$

\medskip

\noindent\textbf{Remark 3} The proofs below are in fact valid more generally: it would suffice, for a given parabolic subgroup $P=MN$ of $G$ containing $B$, to tensor $\ind _K^GV$ with the quotient of $R_M$ in which all $\tau_\alpha-1$ for $\alpha\in \Theta$ are killed.

\medskip

Since we consider only parabolic subgroups in $\Cd$, and all the representations we consider are parabolically induced from analogously defined representations of the Levi quotient of $P_e$, it is enough to prove the theorem when $P_e=G$, i.e.\ $\Delta=\Delta_M\sqcup \Theta$, which we assume from now on.

\subsection{}\label{V.14} Under that assumption $\Delta=\Delta_M \sqcup \Theta$, we prove \ref{V.13} Theorem in a succession of lemmas.

We fix $Q\in \Cd$ and let $M_Q$ be its Levi subgroup containing $M$.

\medskip

\noindent\textbf{Lemma 1} $\rho_Q \subset I_Q$.

\medskip

\noindent\textbf{Proof} Equality is clear when $Q=P$, so we assume $Q\varsupsetneq P$. For each $\alpha\in \Delta_Q-\Delta_P$, let $P^\alpha$ be as in \ref{V.11}. By \ref{V.11} Lemma 2, $\rho_{P^\alpha}$ is included in $I_{P^\alpha}$ so a fortiori $\rho_Q\subset I_{P^\alpha}$. But the subgroup of $G$ generated by the $P^\alpha$'s for $\alpha\in \Delta_Q-\Delta_P$ is $Q$, so  $\cap_{\alpha\in \Delta_Q-\Delta_P} I_{P^\alpha}=I_Q$, and $\rho_Q\subset I_Q$. $\square$

\medskip

To prove equality in Lemma 1, we resort to filtration arguments. In the following $A$, $w$, $A'$ are as in \ref{V.7} and $\pi_{P,A}$, $\pi_{Q,A}$ as in \ref{V.7}, \ref{V.8}.

\medskip

\noindent\textbf{Remark 1} We can also filter $\pi_P^\chi$ by support yielding $(\pi_P^\chi)_A \subset \pi_P^\chi$. But from \ref{V.7} Lemma we get, after tensoring with $\chi:R_M\rg C$, that $\pi_{P,A}\rg \pi_P^\chi$ induces an isomorphism $(\pi_{P,A})^\chi \simeq  (\pi_{P}^\chi)_A$. We let $\pi_{P,A}^\chi$ denote $(\pi_{P}^\chi)_A$.

\medskip

We put $\rho_{Q,A}=\rho_Q \cap \pi_{P,A}^\chi$, $I_{Q,A}=I_Q\cap \pi_{P,A}^\chi$, so $\rho_{Q,A}=\rho_Q\cap I_{Q,A}$.

\medskip

\noindent\textbf{Remark 2}  By \ref{V.8} Corollary~2 and Corollary~3, 
$$0\to R_M \otimes_{\cz_M} \pi_{Q,A}\to R_M \otimes_{\cz_M} \pi_Q\to R_M \otimes_{\cz_M} (\pi_Q /\pi_{Q,A})\to 0$$
is an exact sequence of free $ R_M$-modules  (an extension of free  $R_M$-modules is free and by induction $R_M \otimes_{\cz_M} \pi_{Q,A}$ and $R_M \otimes_{\cz_M} \pi_{Q}$ are free $R_M$-modules). Therefore  the map $(\pi_{Q,A})^\chi \rg \pi_Q^\chi$ is injective.

\medskip

\noindent\textbf{Lemma 2} \textit{(i) If $w\notin W(M_Q)$ then $I_{Q,A}=I_{Q,A'}$, and $\rho_{Q,A}=\rho_{Q,A'}$.}

\textit{(ii) If $w\in W(M_Q)$ the maps $\pi_Q\rg \rho_Q \rg I_Q\rg \pi_P^\chi$ induces isomorphisms }
$$
(\pi_{Q,A})^\chi/(\pi_{Q,A'})^\chi \simeq \rho_{Q,A}/\rho_{Q,A'} \simeq  I_{Q,A}/I_{Q,A'} \simeq \pi_{P,A}^\chi/\pi_{P,A'}^\chi.
$$

\textit{(iii) $\rho_{Q,A}$ is the image of $\pi_{Q,A}$ in $\pi_P^\chi$.}

\medskip

\noindent\textbf{Note} $w\in W(M_Q)$ means that for $\alpha\in \Delta_Q$, $w^{-1}(\alpha)>0$; it is equivalent to  $c_{Q,w}=1$ (\ref{V.8}).

\medskip

\noindent\textbf{Proof} (i) Let $f\in I_{Q,A}-I_{Q,A'}$; then $f$ is not identically $0$ on $PwB$, but its support is   left $Q$-equivariant, so for any $v\in W_{M_Q}$, $f$ is not identically $0$ on $PvwB$.
 If $w\notin W(M_Q)$ we can choose $v\in W_{M_Q}$ so that $vw<w$. That implies $vw\notin A$ by minimality of $w$, a contradiction. So $I_{Q,A}=I_{Q,A'}$ and $\rho_{Q,A}=\rho_{Q,A'}$ follows by intersecting with $\rho_Q$.

(ii) Let $w\in W(M_Q)$. Then $c_{Q,w}=1$ and \ref{V.8} Proposition gives that  the map $\pi_{Q,A} \rg \pi_{P,A}$ induces an isomorphism $\pi_{Q,A}/\pi_{Q,A'}\simeq \pi_{P,A}/\pi_{P,A'}$. Tensoring with $\chi$ gives an isomorphism of  $(\pi_{Q,A})^\chi/(\pi_{Q,A'})^\chi$ onto $(\pi_{P,A})^\chi/(\pi_{P,A'})^\chi$ which is $\pi_{P,A}^\chi/\pi_{P,A'}^\chi$ by Remark 1; since the image of that isomorphism is contained in $\rho_{Q,A}/\rho_{Q,A'}$, itself contained in $I_{Q,A}/I_{Q,A'}$, we get (ii).

(iii) We prove it by descending induction on $\#A$, the case $A=W(M)$ being true by definition of $\rho_Q$. We assume that the result is true for $A$ and prove it for $A'$. By \ref{V.11} Lemma 4 we have
$$
\pi_{Q,A}\subset c_{Q,w} \pi_{P_w,A}+\pi_{Q,A'}.
$$
If $w\notin W(M_Q)$ then $\chi(c_{Q,w})=0$. Hence $\pi_{Q,A}$ and $\pi_{Q,A'}$ have the same image in $\pi_P^{\chi}$, which is $\rho_{Q,A}$ by induction and $\rho_{Q,A'}$ by (i). If $w\in W(M_Q)$ we use the isomorphism $(\pi_{Q,A})^\chi/(\pi_{Q,A'})^\chi \simeq \rho_{Q,A}/\rho_{Q,A'}$ in (ii). Since $(\pi_{Q,A})^\chi \rg \rho_{Q,A}$ is surjective by induction, $(\pi_{Q,A'})^\chi \rg \rho_{Q,A'}$ has to be surjective too. $\square$

\medskip

\noindent\textbf{Lemma 3} $\rho_Q=I_Q$.

\medskip

\noindent\textbf{Proof} By induction on $\#A$: if $\rho_{Q,A'}=I_{Q,A'}$, then Lemma 2 (i), (ii),  and Lemma 1 give $\rho_{Q,A}=I_{Q,A}$. $\square$

\medskip

\noindent\textbf{Lemma 4} \textit{For $Q_1\in \Cd$, $Q_1\varsupsetneq Q$, $\Ker \, \Psi_Q$ contains $\rho_{Q_1}$.}

\medskip
\noindent\textbf{Proof} It enough to show that the composite map $\pi_{Q_1}^\chi \rg \pi_Q^\chi\rg \rho_Q \stackrel{\Psi_Q}{\longrightarrow} \pi_{Q^c}^{\chi}$ is $0$. But it factors as $\pi_{Q_1}^\chi \rg \pi_Q^\chi \stackrel{c_{G,Q}}{\longrightarrow} \pi_{G}^{\chi} \rg \pi_{Q^c}^\chi$ since $c_{Q^c,P}=c_{G,Q}$. From $c_{G,Q}=c_{Q_1,Q}c_{G,Q_1}$ we get $c_{G,Q}\pi_{Q_1}^\chi=c_{Q_1,Q} c_{G,Q_1}\pi_{Q_1}^\chi \subset c_{Q_1,Q}\pi_G^\chi$ which is $0$ since \hbox{$\chi(c_{Q_1,Q})=0$.} $\square$

\medskip
\noindent\textbf{Lemma 5} $\Ker \, \Psi_Q \subset \sum\limits_{Q_1\in \Cd,Q_1\varsupsetneq Q} \rho_{Q_1}$.

\medskip
\noindent\textbf{Proof} We show  by induction on $\#A$ that 
\begin{equation}\tag{$*$} 
\Ker \, \Psi_Q \cap \pi_{P,A}^\chi \subset \sum\limits_{Q_1\in \Cd,Q_1\varsupsetneq Q} \rho_{Q_1}.
\end{equation}
We assume that ($*$)  is true for $A'$ and prove it for $A$. Note that $\Ker \, \Psi_Q \subset \rho_Q \subset \pi_P^\chi$ so $\Ker \, \Psi_Q \cap \pi_{P,A}^\chi = \Ker \, \Psi_Q \cap \rho_{Q,A}$. If $w\notin W(M_Q)$ then $\rho_{Q,A}=\rho_{Q,A'}$ by Lemma~2 (i), so the result is immediate. Assume $w\in W(M_Q)$. On $\rho_{Q,A}/\rho_{Q,A'}$, $\Psi_Q$ induces $\bar{\Psi}_Q: \rho_{Q,A}/\rho_{Q,A'} \rg \pi_{P,A}^\chi/\pi_{P,A'}^\chi \stackrel{c_{Q^c,P}}{\longrightarrow} (\pi_{Q^c,A})^\chi/(\pi_{Q^c,A'})^\chi$. By Lemma 2(ii), the first map is an isomorphism, so we focus  on the second map.

\medskip
\noindent\textbf{Notation} Put $d_w^Q = \prod\limits_{\alpha\in \Delta-\Delta_Q, w^{-1}(\alpha)>0}(\tau_\alpha-1)$, so that  $c_{Q^c,P}= d_w^Q c_{Q^c,w}$ because $\Delta-\Delta_Q=\Delta_{Q^c}-\Delta_M$.

\medskip

By \ref{V.8} Proposition and the remark before it, $c_{Q^c,w}$ gives an isomorphism $\pi_{P,A}/\pi_{P,A'} \stackrel{\sim}{\longrightarrow} \pi_{Q^c,A}/\pi_{Q^c,A'}$. If $d_w^Q=1$ then $\bar{\Psi}_Q$ is injective and $\Ker \, \Psi_Q \cap \pi_{P,A}^\chi=\Ker \, \Psi_Q \cap \pi_{P,A'}^\chi$ so ($*$) follows from the induction hypothesis. Let $d_w^Q\not= 1$, choose $\alpha\in \Delta-\Delta_Q$ with $w^{-1}(\alpha)>0$ and let $Q^\alpha$ be the parabolic subgroup of $G$ corresponding to $\Delta_Q \cup \{\alpha\}$. Then $w\in W(M_{Q^\alpha})$ and Lemma 2 (ii) gives the isomorphism
$$
\rho_{Q^\alpha,A}/\rho_{Q^\alpha,A'} \stackrel{\sim}{\longrightarrow} \pi_{P,A}^\chi/\pi_{P,A'}^\chi.
$$
Let $f\in \Ker \, \Psi_Q \cap \pi_{P,A}^\chi$, and choose $f'\in \rho_{Q^\alpha,A}$ with $f-f' \in \pi_{P,A'}^\chi$. As $f'\in \Ker \, \Psi_Q$ by Lemma~4, $f-f'\in \Ker \, \Psi_Q$ so $f-f'$ belongs to $\sum\limits_{Q_1\in \Cd, Q_1\varsupsetneq Q} \rho_{Q_1}$ by induction; as $f'$ also belongs to that space, the result follows. $\square$

\medskip

\subsection{}\label{V.15} We have proved (i) and (ii) in \ref{V.13} Theorem, and now we turn to  part (iii).
Describing $\Ker \Phi_Q$ is analogous to describing $\Ker \, \Psi_Q$. We let $A$, $w$, $A'$ be as before, and let $\tau_{Q,A}\subset \tau_Q$ be the image of $\pi_{Q,A}$ (or $(\pi_{Q,A})^\chi)$ in $\pi_G^\chi=\tau_G$, via the map   $\pi_Q^\chi \stackrel{c_{G,Q}}{\longrightarrow} \pi_{G}^\chi$.  We observe that  
 $\tau_{Q,A'}\subset \tau_{Q,A}$ and $\tau_{Q_1,A}\subset \tau_{Q,A}$ if $Q_1 \subset Q$ in $ \Cd$. We note also that by \ref{V.14} Remark 2 we have $(\pi_{G,A})^\chi=\tau_{G,A}\subset \pi_G^\chi$.
 
\medskip
\noindent\textbf{Lemma 6} \textit{(i)\ If for some $\alpha\in \Delta-\Delta_Q$, $w^{-1}(\alpha)>0$ then $\tau_{Q,A}=\tau_{Q,A'}$. Otherwise, the natural maps $(\pi_{Q,A})^\chi/(\pi_{Q,A'})^\chi \twoheadrightarrow \tau_{Q,A}/\tau_{Q,A'}\rg \tau_{G,A}/\tau_{G,A'}$ are isomorphisms.
}

\textit{(ii)} $\tau_{Q,A} = \tau_{G,A} \cap \tau_Q$.

\medskip
\noindent\textbf{Proof} (i) Let $\phi\in\pi_{Q,A}$. With $P_w$ as in \ref{V.10}, \ref{V.11} Lemma 4 implies that we can write $\phi=c_{Q,w}\phi_w+\phi'$ with $\phi_w\in\pi_{P_w,A}$ and $\phi'\in\pi_{Q,A'}$. Since $d_w^Qc_{G,w} = c_{G,Q}c_{Q,w} $  we get $c_{G,Q}\phi=d_w^Q(c_{G,w}\phi_w)+c_{G,Q}\phi'$. But 
$c_{G,w}=c_{G,P_w}$ so $c_{G,w}\phi_w$ belongs to $\pi_G$ by \ref{V.6} Proposition. In the first case of (i), 
$\chi(d_w^Q)=0$, so $\phi$ has the same image as $\phi'$  in $\tau_Q$; this implies $\tau_{Q,A}=\tau_{Q,A'}$. Let us assume we are in the second case of (i), so $d_w^Q=1$. Consider the natural inclusions
$$
\pi_{G,A}/\pi_{G,A'} \hookrightarrow \pi_{Q,A}/\pi_{Q,A'} \hookrightarrow \pi_{P,A}/\pi_{P,A'}.
$$
By \ref{V.8} Proposition, the first space is $c_{G,w}(\pi_{P,A}/\pi_{P,A'})$ and the second is $c_{Q,w}(\pi_{P,A}/\pi_{P,A'})$. Consequently, $c_{G,Q}(\pi_{Q,A}/\pi_{Q,A'}) = \pi_{G,A}/\pi_{G,A'}$ since $d_w^Q =1$. Thus $c_{G,Q}$ induces a surjective map of $\pi_{Q,A}/\pi_{Q,A'}$ onto $\pi_{G,A}/\pi_{G,A'}$.  But by \ref{V.4} Lemma (i) (applied to $M$) $c_{G,Q}$ acts injectively on $\sigma$ hence on $\pi_{P,A}/\pi_{P,A'}$, so we actually get an isomorphism. Tensoring with $\chi$ we get an isomorphism $(\pi_{Q,A})^\chi/(\pi_{Q,A'})^\chi \to \tau_{G,A}/\tau_{G,A'}$; but this factors as in the statement of (i), so (i) follows again.

(ii) We proceed by descending induction on $\#A$, the case $A=W(M)$ being obvious. The containment $\tau_{Q,A'} \subset \tau_{G,A'} \cap \tau_Q$ is clear, and we have $\tau_{Q,A}=\tau_{G,A} \cap \tau_Q$ by induction. In the first case of (i) $\tau_{Q,A'}=\tau_{Q,A} =\tau_{G,A}\cap \tau_Q \supset \tau_{G,A'} \cap \tau_Q$ so $\tau_{Q,A'}= \tau_{G,A'}\cap \tau_Q$. In the second case of (i), $\tau_{Q,A}/\tau_{Q,A'} \rg \tau_{G,A}/\tau_{G,A'}$ is an isomorphism; as moreover $\tau_{G,A'}\cap \tau_Q \subset \tau_{Q,A}$ by induction, 
the result follows. $\square$

\medskip
\noindent\textbf{Lemma 7} \textit{For $Q_1\in \Cd$, $Q_1\varsubsetneq Q$, then $\tau_{Q_1} \subset \Ker \Phi_Q$.}

\medskip
\noindent\textbf{Proof} Let $P_1$ be the parabolic subgroup corresponding to $\Delta_{Q_1} \sqcup (\Delta-\Delta_Q) = \Delta_{Q_1}\cup \Delta_{Q^c}$. Since $Q_1\varsubsetneq Q$, we get $P_1\varsubsetneq G$. We have $c_{P_1,Q_1}\pi_{Q_1} \subset \pi_{P_1}\subset \pi_{Q^c}$ so $c_{G,Q_1}\pi_{Q_1} \subset  c_{G,P_1} \pi_{Q^c}$. As $\chi(c_{G,P_1})=0$ the image of $\pi_{Q_1} \stackrel{c_{Q^c,Q_1}}{\longrightarrow} \pi_{Q^c} \rg \pi_{Q^c}^\chi$ is 0; but that image is $\Phi_Q(\tau_{Q_1})$.~$\square$

\medskip
\noindent\textbf{Lemma 8 } $\Ker \Phi_Q \subset \sum\limits_{Q_1\in \Cd,Q_1\varsubsetneq Q}\tau_{Q_1}$.

\medskip
\noindent\textbf{Proof} We prove that $\Ker \,\Phi_Q\cap \tau_{G,A}$ is contained in the right-hand side, by induction on $\#A$. In the first case of Lemma 6 (i), $\tau_{Q,A}=\tau_{Q,A'}$,  so $\tau_{G,A}\cap \tau_Q=\tau_{G,A'} \cap \tau_Q$ by Lemma 6 (ii). Consequently, $\Ker \,\Phi_Q \cap \tau_{G,A}=\Ker \Phi_Q \cap \tau_{G,A'}$ and we are done. So we assume that for all $\alpha\in \Delta-\Delta_Q=\Delta_{Q^c}-\Delta_P$ we have $w^{-1}(\alpha)<0$. 
\medskip
On $\tau_{Q,A}/\tau_{Q,A'}$, $\Phi_Q$ induces $\bar{\Phi}_Q: \tau_{Q,A}/\tau_{Q,A'} \rg (\pi_{G,A})^\chi/(\pi_{G,A'})^\chi{\longrightarrow} (\pi_{Q^c,A})^\chi/(\pi_{Q^c,A'})^\chi$, where the first map is an isomorphism by Lemma 6 (i), and the second comes, upon tensoring with $\chi$, from the inclusion of $\pi_{G,A}/\pi_{G,A'}$ into $\pi_{Q^c,A}/\pi_{Q^c,A'}$. By \ref{V.8} Proposition, we have, inside $\pi_{P,A}/\pi_{P,A'}$, $\pi_{G,A}/\pi_{G,A'}=c_{G,w}(\pi_{P,A}/\pi_{P,A'})$, and $\pi_{Q^c,A}/\pi_{Q^c,A'}=c_{Q^c,w}(\pi_{P,A}/\pi_{P,A'})$. If for all $\alpha\in \Delta-\Delta_{Q^c}$ we have $w^{-1}(\alpha)>0$, then $c_{G,w}=c_{Q^c,w}$, and $\pi_{G,A}/\pi_{G,A'}=\pi_{Q^c,A}/\pi_{Q^c,A'}$; thus $\Ker \,\Phi_Q \cap \tau_{G,A}=\Ker \Phi_Q \cap \tau_{G,A'}$, so we conclude by induction. In the opposite case, choose $\alpha\in \Delta-\Delta_{Q^c}=\Delta_{Q}-\Delta_P$ with $w^{-1}(\alpha)<0$, and let $Q_\alpha$ correspond to $\Delta_{Q}-\{\alpha\}$. Then  $\tau_{Q_\alpha,A} /\tau_{Q_\alpha,A'} \rg \tau_{G,A}/\tau_{G,A'}$ is an isomorphism by Lemma~6~(i). If $f\in \Ker \Phi_Q \cap \tau_{Q,A}$, there is $f'\in \tau_{Q_\alpha,A}$ with $f-f'\in \tau_{G,A'}$. As $\tau_{Q_\alpha} \subset \tau_Q$ we have $f-f'\in \tau_{G,A'} \cap \tau_Q=\tau_{Q,A'}$ by Lemma 6 (ii). 
  Lemma~7 gives $\Phi_Q(f')=0$, so $\Phi_Q(f-f')=0$ and by induction $f-f'$ belongs to the right-hand side of Lemma~8; since $f'\in \tau_{Q_\alpha}$ also belongs to that space, so does $f$.~$\square$

\subsection{}\label{V.16} It remains to prove (iv) of \ref{V.13} Theorem.

\medskip
\noindent\textbf{Lemma 9} \textit{Let $\cp \subset \Cd$. Then $\Big(\sum\limits_{Q_1\in \cp} \tau_{Q_1}\Big) \cap \tau_{G,A}= \sum\limits_{Q_1\in \cp} \tau_{Q_1,A}$.}

\medskip
\noindent\textbf{Proof} The containment $\supset$ is clear; we prove the other direction by descending induction on $\#A$. Let $\cp^-=\{Q_1\in \cp  \mid  w^{-1}(\alpha) <0 $    for any  $\alpha\in \Delta-\Delta_{Q_1}  \} $. If $\cp^-$ is empty then $\tau_{Q_1,A}=\tau_{Q_1,A'}$ for any $Q_1\in \cp$ (Lemma~6 (i)), and we have nothing to prove. Assume $\cp^-$ is not empty, and put $Q_\cap =\bigcap\limits_{Q_1\in \cp^-}Q_1$. Then for $\alpha\in \Delta -\Delta_{Q_\cap}$ we have $w^{-1}(\alpha)<0$ so by Lemma 6 (i) the map $\tau_{Q_\cap,A}\rg \tau_{G,A}/\tau_{G,A'}$ is surjective. For $Q_1\in \cp$ let $f_{Q_1}\in \tau_{Q_1}$ be chosen so that $\sum\limits_{Q_1\in \cp} f_{Q_1}\in \tau_{G,A'}$; by the inductive hypothesis we may assume that   all $f_{Q_1}\in \tau_{Q_1,A}$. For $Q_1\in \cp - \cp^-$, we even have $f_{Q_1}\in \tau_{Q_1,A'}$ by Lemma 6 (i). Fix $Q_2\in \cp^-$; for $Q_1\in \cp^-$, $Q_1\not= Q_2$ choose $f_{Q_1}' \in \tau_{Q_\cap,A}$ with $f_{Q_1}-f_{Q_1}' \in\tau_{G,A'}$. Since $\tau_{Q_\cap,A} \subset \tau_{Q_1,A}$, $f_{Q_1}-f_{Q_1}'$ belongs to $\tau_{G,A'} \cap \tau_{Q_1}=\tau_{Q_1,A'}$. So $\sum\limits_{Q_1\in \cp} f_{Q_1}$ appears as $f_{Q_2}+ \sum\limits_{Q_1\in \cp^-,Q_1\not=Q_2}f_{Q_1}'$ plus terms in $\sum\limits_{Q_1\in \cp} \tau_{Q_1,A'}$. But for $Q_1\in \cp^-$, $Q_1 \not= Q_2$, $f_{Q_1}'$ belongs to $\tau_{Q_\cap,A} \subset \tau_{Q_2,A}$ so $f_{Q_2} + \sum\limits_{Q_1\in \cp^-,Q_1\not=Q_2} f_{Q_1}'$ belongs to $\tau_{Q_2}\cap \tau_{G,A'}= \tau_{Q_2,A'} \subset \sum\limits_{Q_1\in \cp}\tau_{Q_1,A'}$.~$\square$

\medskip 
We finally prove (iv) of \ref{V.13} Theorem. Fix $Q\in   \Cd$ and let $\cp\subset \Cd$. It is clear that
$$
\Big(\sum_{Q_1\in \cp} \tau_{Q_1}\Big) \cap \tau_Q  \supset \sum_{Q_1\in \cp} (\tau_{Q_1} \cap \tau_Q) \supset \sum_{Q_1\in \cp} \tau_{Q_1\cap Q}.
$$
We prove now 
$$
\Big(\sum_{Q_1\in \cp} \tau_{Q_1}\Big) \cap \tau_{Q,A} \subset \sum_{Q_1\in \cp} \tau_{Q_1\cap Q}\ \mathrm{by\ induction\ on\ }\#A. 
$$
If there is $\alpha\in \Delta-\Delta_Q$ with $w^{-1}(\alpha) >0$ then $\tau_{Q,A}=\tau_{Q,A'}$ (Lemma~6 (i)) and there is nothing to prove, so we assume the contrary. By Lemma~9
$$
\Big(\sum_{Q_1\in \cp} \tau_{Q_1}\Big)\cap \tau_{Q,A} =\Big(\sum_{Q_1\in \cp} \tau_{Q_1,A}\Big) \cap \tau_{Q,A}.
$$
Let  $\cp^- \subset \cp$ be the same subset as in the proof of Lemma~9. If $\cp^-$ is empty, then $\tau_{Q_1,A}=\tau_{Q_1,A'}$ for any $Q_1$ in $\cp$. Hence
$$
\Big(\sum_{Q_1\in \cp} \tau_{Q_1}\Big) \cap \tau_{Q,A} = \Big(\sum_{Q_1\in \cp} \tau_{Q_1,A'}\Big) \cap \tau_{Q,A} = \Big(\sum_{Q_1\in \cp} \tau_{Q_1,A'}\Big) \cap \tau_{Q,A'},
$$
 and the result follows from Lemma 9 and the induction hypothesis. Now assume $\cp^-\not=\emptyset$, and write $Q_\cap =Q\cap \bigcap\limits_{Q_1\in \cp^-}Q_1$; then for $\alpha\in \Delta-\Delta_{Q_\cap}$, $w^{-1}(\alpha)<0$ and again $\tau_{Q_\cap,A}\rg \tau_{G,A}/\tau_{G,A'}$ is surjective. For $Q_1\in \cp$ let $f_{Q_1}\in \tau_{Q_1}$ be chosen so that $\sum\limits_{Q_1\in \cp} f_{Q_1} \in \tau_{Q,A}$. 
By Lemma 9 we may assume $f_{Q_1}\in \tau_{Q,A}$. 
For $Q_1\in \cp^-$, choose $f_{Q_1}'\in \tau_{Q_\cap,A}$ with $f_{Q_1}-f_{Q_1}' \in \tau_{G,A'}$ (then $f_{Q_1}-f_{Q_1}' \in \tau_{Q_1,A'})$. Write 
$$
\sum_{Q_1\in \cp} f_{Q_1} = \sum_{Q_1\in \cp^-} (f_{Q_1} -f_{Q_1}') + \sum_{Q_1\in \cp-\cp^-} f_{Q_1} + \sum_{Q_1\in \cp^-} f_{Q_1}'.
$$
 We examine the right hand side. The last term belongs to $\tau_{Q_\cap,A} \subset \tau_{Q,A}$, so the sum of the first two belongs to $\tau_Q$. As each summand in those two terms indexed by $Q_1$ is  in $\tau_{Q_1,A'}$, their sum belongs to $(\sum\limits_{Q_1\in \cp} \tau_{Q_1,A'}) \cap \tau_Q$, which is in $\sum\limits_{Q_1\in \cp} \tau_{Q_1\cap Q}$ by the induction hypothesis. But for $Q_1\in \cp^-$, $f_{Q_1}'\in \tau_{Q_\cap,A}$, and $\tau_{Q_\cap} \subset \tau_{Q_1\cap Q}$ since $Q_\cap \subset Q_1 \cap Q$. Thus the third term also belongs to $\sum\limits_{Q_1\in \cp^-} \tau_{Q_1\cap Q}$. $\square$

\section{Consequences of the classification}\label{VI}

\subsection{}\label{VI.1}  
We recall from \ref{I.3} that a  representation of $G$ is supercuspidal if it is irreducible, admissible, and does not appear as a  subquotient of a parabolically induced representation $\Ind_P^G\sigma$, where $P$ is a  proper parabolic subgroup of $G$ and $\sigma$ an irreducible admissible representation of the Levi quotient of~$P$.  

 It is well known \cite{BL1,Br} that there exist irreducible admissible supercuspidal representations when $G=\GL_2(\mathbb Q_p)$, therefore  the following proposition shows that  we cannot drop the condition that $\sigma$ be irreducible admissible in the definition of supercuspidality, unlike for representations of $G$ over a field of characteristic different from $p$.
 
\medskip 
 \noindent \textbf{Proposition} \textit{Any irreducible   \repr\ $\pi$ of $G$ is a subquotient of $\Ind_B^G \sigma$ for some  \repr\ $\sigma $ of $Z$.}

\medskip 
\noindent \textbf{Proof} The smoothness of $\pi$ implies that $\pi$ has a weight $V$. The irreducibility of $\pi$ implies that $\pi$ is a quotient of $\ind_K^G V$.  The representation $\ind_K^G V$ embeds in $\Ind_B^G(\ind_{Z^0}^Z V_{U^0})$ by the intertwiner $\mathcal I$ of  \ref{III.13}. $\square$

\medskip
 
\subsection{}\label{VI.2} We derive the desired consequences of \ref{I.5} Theorem~4. Mostly we follow the pattern of \cite{He2}.  

 We now prove \ref{I.5} Theorem 5, which we recall.

\vskip2mm

\noindent \textbf{Theorem} \textit{Let $\pi$ be an irreducible admissible \repr\ of $G$. Then $\pi$ is supercuspidal if and only if $\pi$ is supersingular.
}

\vskip2mm

As observed in the introduction, this theorem shows that the notion of supersingularity, for an irreducible admissible representation of $G$, is independent of the choices of $\gs, \bb, K$.

\medskip

\noindent\textbf{Proof} Let $\pi$ be supercuspidal. By  \ref{I.5} Theorem 4, there is a supersingular $B$-triple $(P,\sigma,Q)$ such that $\pi \simeq I(P,\sigma,Q)$. By \ref{III.24} Proposition, $I(P,\sigma,Q)$ is a component of $\Ind_P^G\sigma$, so $P=G$ and $\pi\simeq \sigma$ is supersingular.

Let $\pi$ be supersingular. Assume it occurs as a subquotient of $\Ind_P^G\sigma$ for a parabolic subgroup $P$ of $G$ and an irreducible admissible \repr\ $\sigma$ of the Levi quotient $M$ of $P$; we may and do assume that $P$ contains $B$.  By \ref{I.5} Theorem~4,  \ref{III.24} Proposition, and transitivity of parabolic induction, we may assume that $\sigma$ is supersingular. By \ref{III.24} Proposition, $\pi$ is isomorphic to some $I(P,\sigma,Q)$ and  \ref{I.5} Theorem~4 implies that $P=G$, so that $\pi$ is indeed  supercuspidal. $\square$

\medskip Theorems 1 to 3 in Section \ref{I.3} are now rather immediate. They follow from \ref{I.5} Theorem~4 and the following elementary observations:

\begin{enumerate}
\item Any triple is $G$-conjugate to a $B$-triple.
\item A $B$-triple is supersingular if and only if it is supercuspidal (by the theorem).
\item $I(P,\sigma,Q)\simeq I(P',\sigma',Q')$ if the triples $(P,\sigma,Q)$, $(P',\sigma',Q')$ are $G$-conjugate.
\end{enumerate} 

\subsection{}\label{VI.3} We also have the desired consequence about supercuspidal support.

\vskip2mm

\noindent \textbf{Proposition} \textit{Let $\pi$  be an irreducible admissible \repr\ of $G$. Then there is a parabolic subgroup $P$ of $G$ and a supercuspidal \repr\ $\sigma$ of the Levi quotient of $P$ such that $\pi$ is a subquotient of $\Ind_P^G\sigma$. If $P_1$ is a parabolic subgroup of $G$ and $\sigma_1$ a supercuspidal \repr\ of the Levi quotient of $P_1$ such that $\pi$ is a subquotient of $\Ind_{P_1}^G \sigma_1$, then there is $g$ in $G$ such that $P_1=gPg^{-1}$ and that $\sigma_1$ is equivalent to $x\mapsto\sigma(g^{-1}xg)$.
}

\vskip2mm

\noindent \textbf{Proof}  By \ref{I.3} Theorem 3, $\pi$ has the form $I(P,\sigma,Q)$ for some \superc\ triple $(P,\sigma,Q)$ and the first assertion comes from \ref{III.24} Proposition. The uniqueness assertion is derived in the same way from  \ref{I.3} Theorem 2. $\square$

\vskip2mm We say that the \textbf{supercuspidal support} of $\pi$ is the class of $(P,\sigma)$ for the equivalence relation appearing in the proposition.

\subsection{}\label{VI.4} We give one more consequence mentioned in the introduction.

\medskip

\noindent\textbf{Proposition}  \textit{Let $(P,\sigma,Q)$ be a $B$-triple. Assume that $\sigma$ is a supercuspidal (or equivalently, supersingular) representation of $M$. Then $I(P,\sigma,Q)$ is finite-dimensional if and only if $P=B$ and $Q=G$.}

\medskip

\noindent\textbf{Proof} As $Z$ is compact mod centre, any irreducible representation $\tau$ of $Z$ is finite dimensional \cite{Hn,Vig2} and consequently supercuspidal. If $P(\tau)=G$ then $I(B,\tau,G)={}^e\tau$ is finite dimensional. Conversely, let $\pi$ be a finite-dimensional irreducible representation of $G$. Then its kernel is an open normal subgroup of $G$. Considering $\iota:G^{\is}\rg G$ as in Chapter~\ref{II}, $\Ker(\sigma\circ\iota)$ is an open normal subgroup of $G^{\is}$ which by \ref{II.3} Proposition  has to be $G^{\is}$ itself. Thus $\pi$ is trivial on $G'$ and since $G=ZG'$, $\pi$ restricts to an irreducible (supercuspidal) representation $\tau$ of $Z$; we have $P(\tau)=G$ and ${}^e\tau=\pi$, $\pi=I(B,\tau,G)$. $\square$

\medskip 

\subsection{}\label{VI.5} It is worth noting that our results recover the classifications obtained previously in special cases. Keep the notation of Chapter~\ref{III}. When $\bg$ is split, then for $\alpha\in \Delta$, $Z\cap M_\alpha'$ is simply the image in $Z=S$ of the coroot $\alpha^\vee$, so our classification is the same as that of \cite{Abe}; it also gives the classification of \cite{He2} for $\bg=\bgl_n$. 
Other special cases are worth mentioning: groups of semisimple rank 1 and inner forms of~$\GL_n$.
Of course if $\bg$ has relative rank 0, all irreducible representations of $G$ are finite dimensional and supercuspidal, and our classification theorem says nothing.
If $\bg$ has relative semisimple rank 1, the classification is rather simple (see also \cite{BL1, BL2, Abd, Che, Ko, Ly2}). An irreducible admissible representation $\pi$ of $G$ falls into one (and only one) of the following cases:

\medskip

\noindent 1) $\pi$ is supercuspidal (hence infinite dimensional), i.e.\    $\pi\simeq I(G,\pi,G)$.

\noindent 2) $\pi$ is finite dimensional; it is then trivial on $G'$ and restricts to an irreducible representation $\tau$ of $Z$, trivial on $Z\cap G'$, and $\pi\simeq I(B,\tau,G)$.

\noindent 3) $\pi\simeq \sigma\otimes \St_B^G$ where $\sigma$ is as in 2), i.e.\  $\pi\simeq I(B,\sigma|_Z,B)$.

\noindent 4) $\pi \simeq I(B,\tau,B)$ where $\tau$ is an irreducible representation of $Z$ (hence finite dimensional and supercuspidal) which is not trivial on $Z\cap G'$.

\subsection{}\label{VI.6}Let us briefly consider the case of inner forms of general linear groups. Thus $\bg=\GL_{n/D}$ where $D$ is a central division algebra  of finite degree over $F$. We take for $S$ the diagonal subgroup $(F^\times)^n$ (so that $Z$ is the diagonal subgroup $(D^\times)^n)$, and for $B$ the upper triangular subgroup. We can take $K=\GL_n(\Oo_D)$ where $\Oo_D$ is the ring of integers of $D$; all other special parahoric subgroups of $G$ are conjugate to~$K$.

A parabolic subgroup $P$ of $G$ containing $B$ is an upper triangular block subgroup, and the corresponding Levi subgroup $M$ is the  block diagonal subgroup. If the successive blocks down the diagonal have size $n_1,\ldots,n_r$, then $M$ appears as $M_1\times\cdots\times M_r$, $M_i=\GL_{n_i}(D)$ and an irreducible admissible representation of $M$ factors as a tensor product $\pi_1 \otimes \cdots \otimes \pi_r$, where $\pi_i$ is an irreducible admissible representation of $M_i$ for $i=1,\ldots,r$ determined up to isomorphism. (Conversely such a tensor product is an irreducible admissible representation of $M$: the reader can devise a proof as suggested in \cite{He2}, perhaps using \cite[7.10 Lemma]{HV2}.) Note that the group $G'$ is the kernel of the non-commutative determinant $\det: G\rg F^\times$. Parameters for the irreducible admissible representations of $G$ can then be described in a way entirely parallel to the case $D=F$ obtained in \cite{He2}. (The cases of $\GL_n(D)$ where $n \le 3$ are treated in T.\ Ly's Ph.D.\ thesis~\cite{Ly2}, \cite[Chapter 3]{Ly3}.)

We simply state the results, leaving to the reader the translation from our classification in this paper.

For $i=1,\ldots,r$ let $\pi_i$ be a representation of $M_i$ which is either supercuspidal or of the form $\chi_i\circ \det$ for some character $\chi_i: F^\times \rg C^\times$; if for two consecutive indices $i$, $i+1$ we have $\pi_i=\chi_i\circ\det$ and $\pi_{i+1}=\chi_{i+1}\circ\det$, assume $\chi_i\not= \chi_{i+1}$. 

For each index $i$ such that $\pi_i=\chi_i\circ \det$, choose an upper (block) triangular parabolic subgroup $Q_i$ of $M_i$, and put $\sigma_i=(\chi_i\circ \det) \otimes \St_{Q_i}^{M_i}$; for other indices $i$ put $\sigma_i=\pi_i$. Then $\Ind_{P}^G(\sigma_1\otimes\cdots \otimes \sigma_r)$ is irreducible and admissible. Conversely any irreducible admissible representation of $G$ has such a shape, where the integers $n_1,\ldots,n_r$, the parabolic subgroups $Q_i$ of $M_i=\GL_{n_i}(D)$, and the isomorphism  classes of the $\pi_i$, are determined.

\end{document}